\documentclass[12pt]{article} 

\usepackage{latexsym, 
amscd, 
graphicx, amssymb}

\newtheorem{thm}{Theorem}[section]

\newtheorem{prop}[thm]{Proposition}
\newtheorem{cor}[thm]{Corollary}
\newtheorem{lemma}[thm]{Lemma}

\newcommand{\halmos}{\rule{1ex}{1.4ex}}

\newcommand{\text}[1]{\mbox{\rm #1}}

\newcommand{\nn}{\nonumber \\}

 \newcommand{\pf}{{\it Proof.}\hspace{2ex}}
 
 \newcommand{\epfv}{\hspace*{\fill}\mbox{$\halmos$}\vspace{1em}}

\newcommand{\wt}{\mbox{\rm wt}\ }

\newcommand{\tr}{\mbox{\rm Tr}}
\newcommand{\A}{\mathcal{A}}
\newcommand{\Y}{\mathcal{Y}}
\newcommand{\C}{\mathbb{C}}
\newcommand{\Z}{\mathbb{Z}}

\newcommand{\Q}{\mathbb{Q}}
\newcommand{\N}{\mathbb{N}}

\newcommand{\F}{\mathcal{F}}

\title{ {\bf Vertex operator algebras and
the Verlinde conjecture} }
\date{}
\author{Yi-Zhi Huang}

\begin{document}

\bibliographystyle{alpha}
\maketitle

\begin{abstract} 
We prove the Verlinde conjecture in the following 
general form: 
Let $V$ be a simple vertex operator algebra satisfying the 
following conditions: (i) $V_{(n)}=0$ for $n<0$, 
$V_{(0)}=\mathbb{C}\mathbf{1}$ and 
$V'$ is isomorphic to $V$ as a $V$-module. (ii)
Every $\mathbb{N}$-gradable weak $V$-module is completely 
reducible. (iii) $V$ is $C_{2}$-cofinite. (In the presence 
of Condition (i), Conditions (ii) and (iii) are equivalent to
a single condition, namely, that every weak 
$V$-module is completely 
reducible.)
Then the matrices formed by the fusion rules 
among the irreducible $V$-modules are diagonalized by 
the matrix given by the action of the modular transformation
$\tau\mapsto -1/\tau$ on the space of  characters
of irreducible $V$-modules. Using this result, we obtain the 
Verlinde formula for the fusion rules. We also 
prove that the matrix associated to the 
modular transformation
$\tau\mapsto -1/\tau$ is symmetric. 
\end{abstract}

\renewcommand{\theequation}{\thesection.\arabic{equation}}
\renewcommand{\thethm}{\thesection.\arabic{thm}}
\setcounter{equation}{0}
\setcounter{thm}{0}
\setcounter{section}{-1}

\section{Introduction}

In the present paper, we formulate and prove a general version
of the Verlinde conjecture and prove the Verlinde formula
for fusion rules using the representation theory of vertex 
operator algebras. 

The Verlinde conjecture \cite{V} in conformal field theory states that
the action of the modular transformation $\tau \mapsto -1/\tau$ on the
space of characters of a rational conformal field theory
diagonalizes the fusion rules. Except for some particular examples
(see below), the general Verlinde conjecture has been an open problem
for twenty years.
In \cite{MS1}, Moore and Seiberg showed 
on a physical level of rigor that
this conjecture follows from the axioms for rational conformal
field theories (see \cite{K}, \cite{S1}, \cite{S2} and \cite{S3} for 
axioms for conformal field theories and 
\cite{MS2} for axioms
and assumptions for rational conformal field theories on a physical
level of rigor).  Axioms for rational conformal field theories 
are in fact much stronger than statements such as the Verlinde conjecture.
The work \cite{MS1} \cite{MS2} of Moore and Seiberg
greatly advanced our understanding of the structure of conformal field
theories.  However, it is a very hard problem to actually {\it construct}
theories
satisfying these axioms {\it mathematically} and therefore the existence of
rational conformal field theories is a very strong assumption.
In fact, the construction of full rational conformal field theories,
especially the hard part of verifying the axioms, 
need, among many other things, the Verlinde conjecture (see \cite{HK1}
and \cite{HK2}). 

Assuming that the axioms for
higher-genus rational conformal field theories are 
satisfied, the Verlinde 
conjecture leads to a Verlinde formula \cite{V} for the dimensions of
the spaces of conformal blocks on higher-genus Riemann surfaces. 
In the special case of 
the conformal field theories associated to affine Lie algebras
(the Wess-Zumino-Novikov-Witten models), 
this Verlinde formula gives a surprising formula 
for the dimensions of the spaces of sections of the 
"generalized theta divisors"  and has given rise to a great deal of 
excitement and new mathematics. See the works \cite{TUY} by 
Tsuchiya-Ueno-Yamada, \cite{BL} by Beauville-Laszlo, 
\cite{F} by Faltings and \cite{KNR} by Kumar-Narasimhan-Ramanathan
for details and proofs of this particular case of the Verlinde formulas.

In this
paper, using the results on the duality and modular invariance of
genus-zero and genus-one correlation functions, especially those
obtained recently in \cite{H8} and \cite{H9}, we formulate and prove a
general version of the Verlinde conjecture in the framework of the
theory of vertex operator algebras, which were first introduced and
studied in mathematics by Borcherds \cite{B} and
Frenkel-Lepowsky-Meurman \cite{FLM}.  Our theorem assumes only that the
vertex operator algebra we consider satisfies certain natural grading,
finiteness and reductivity properties (see below).  This result is part
of a program to mathematically construct conformal field theories in the
sense of Kontsevich \cite{K} 
and Segal \cite{S1} \cite{S2} \cite{S3} from representations
of vertex operator algebras. It will be a main step towards a
mathematical proof of the Verlinde formula for the dimensions of
conformal blocks on higher-genus Riemann surfaces for conformal field
theories other than the Wess-Zumino-Novikov-Witten models. 

To state our main result, we need to discuss briefly
matrices formed by fusion rules and modular transformations. 
Let $V$ be a simple vertex operator algebra
satisfying the following conditions: (i) $V_{(n)}=0$ for $n<0$,
$V_{(0)}=\mathbb{C}\mathbf{1}$ and $V'$ is isomorphic to $V$ as a
$V$-module. (ii) Every $\mathbb{N}$-gradable weak $V$-module is completely
reducible. (iii) $V$ is $C_{2}$-cofinite (that is, $\dim V/C_{2}(V)<\infty$
where $C_{2}(V)$ is the subspace of $V$ spanned by $u_{-2}v$ for $u, v\in V$. 
By results of Li \cite{L} and 
Abe-Buhl-Dong \cite{ABD}, 
Conditions (ii) and (iii) are equivalent to a 
single condition that every weak $V$-module is completely 
reducible.
Let $\A$ be the set of 
equivalence classes of irreducible $V$-modules. For $a\in \A$,
we choose a representative $W^{a}$ in $a$  
such that $W^{e}=V$ where $e$ is the equivalence class
containing $V$. 

For $a_{1}, a_{2}, 
a_{3}\in \A$, let $N_{a_{1}a_{2}}^{a_{3}}
=N_{W^{a_{1}}W^{a_{2}}}^{W^{a_{3}}}$ be the corresponding
fusion rules, that is, the dimensions of the spaces of 
intertwining operators of types ${W^{a_{3}}\choose 
W^{a_{1}}W^{a_{2}}}$ (see \cite{FHL}). 
For $a\in \A$, let $\mathcal{N}(a)$ be the matrix 
whose entries are $N_{aa_{1}}^{a_{2}}$ for $a_{1}, a_{2}\in \A$,
that is, 
$$\mathcal{N}(a)=(N_{aa_{1}}^{a_{2}}).$$

For $a\in \A$, we have the characters 
$\tr_{W^{a}}q_{\tau}^{L(0)-\frac{c}{24}}$,
also called shifted graded dimensions, 
where $q_{\tau}=e^{2\pi \tau}$ and $\tau\in \mathbb{H}$.
In \cite{Z}, Zhu proved under certain conditions
that the maps given by
$$u\mapsto \tr_{W^{a}}Y_{W^{a}}(e^{2\pi zL(0)}u, e^{2\pi z})
q_{\tau}^{L(0)-\frac{c}{24}}$$
for $u\in V$, where $a\in \mathcal{A}$,
are linearly independent and 
there exist
$S_{a_{1}}^{a_{2}}\in \C$ for $a_{1}, a_{2}\in \A$ such that
\begin{eqnarray*}
\lefteqn{\tr_{W^{a_{1}}}Y_{W^{a_{1}}}\left(
e^{-\frac{2\pi z}{\tau}L(0)}\left(-\frac{1}{\tau}\right)u, 
e^{-\frac{2\pi z}{\tau}}\right)
q_{-\frac{1}{\tau}}^{L(0)-\frac{c}{24}}}\nn
&&=\sum_{a_{2}\in \A}S_{a_{1}}^{a_{2}}
\tr_{W^{a_{2}}}Y_{W^{a_{2}}}(e^{2\pi zL(0)}u, e^{2\pi z})
q_{\tau}^{L(0)-\frac{c}{24}}.
\end{eqnarray*}
in \cite{DLM}, Dong, Li and Mason
improved Zhu's results above by showing that they also hold under
the conditions (slightly weaker than what) we assume in our paper. 
In partiuclar, when $u=\mathbf{1}$, we have 
$$\tr_{W^{a_{1}}}q_{-\frac{1}{\tau}}^{L(0)-\frac{c}{24}}
=\sum_{a_{2}\in \A}S_{a_{1}}^{a_{2}}
\tr_{W^{a_{2}}}q_{\tau}^{L(0)-\frac{c}{24}}.$$

Now we can state the main result of the present paper:
The matrix 
$S=(S_{a_{1}}^{a_{2}})$ diagonalizes 
$\mathcal{N}(a)$ for $a\in \A$.  See Theorem
\ref{main} for a more complete and precise statement. 
Using this result, we obtain the Verlinde formula for 
the fusion rules for such a vertex operator algebra. 
We also prove that the matrix $S$
is symmetric. 

Note that fusion rules and modular transformations of characters were
already defined mathematically for suitable vertex operator algebras
by Frenkel-Huang-Lepowsky in \cite{FHL} using intertwining operators
and by Zhu in \cite{Z} using his modular invariance theorem,
respectively. The $C_{2}$-cofiniteness condition was also introduced
in \cite{Z}.  These are already enough for the mathematical 
formulation of the general version of the Verlinde conjecture proved
in this paper. Further results on intertwining operators and 
modular invariance
were obtained in \cite{HL1}--\cite{HL4} by Huang-Lepowsky, in
\cite{H1}, \cite{H3} and \cite{H7} by the author, in \cite{DLM} by
Dong-Li-Mason and in \cite{M} by Miyamoto. However, the proof of the 
Verlinde conjecture in the present paper requires much stronger 
results than these.  We need the duality and modular
invariance properties for genus-zero and genus-one {\it multi-point}
correlation functions constructed from intertwining operators for a
vertex operator algebra satisfying the conditions above. These
properties have been proved recently in \cite{H8} and \cite{H9} by the
author so that the proof of the Verlinde conjecture in the present
paper has now become possible.

The main content of the present paper is to establish mathematically certain
formulas needed in the proof of the Verlinde conjecture. Most of the
formulas were first obtained on a physical level of rigor by Moore and
Seiberg \cite{MS1} \cite{MS2} using the (assumed) 
axioms for rational conformal
field theories. In the present paper, our formulations and proofs are
based on the results obtained in the theory of vertex operator
algebras. We use only those results which have been established
mathematically and we do not assume that all axioms for rational
conformal field theories hold. In \cite{MS1} and \cite{MS2}, Moore and
Seiberg discussed many of the subtle technical details using examples
such as the minimal models. Many of these discussions cannot be
generalized to general cases. For example, in \cite{MS1} and
\cite{MS2}, spaces of chiral vertex operators (intertwining operators)
are identified with tensor products of spaces of lowest weight vectors
in modules.  This identification allows them to give an $S_{3}$ action
easily on the direct sum of these spaces and the formulas obtained in
\cite{MS1} and \cite{MS2} depend heavily on this action. It is known
that in general spaces of intertwining operators cannot be identified
with tensor products of spaces of lowest weight vectors and thus, even if we
assume that all axioms for rational conformal field theories hold, the
method based on particular examples in \cite{MS1} and \cite{MS2}
cannot be directly adopted to establish the formulas we need. In the
present paper, we define this action of $S_{3}$ using the
skew-symmetry for intertwining operators and the contragredient
intertwining operators in \cite{FHL} and \cite{HL2} and prove all the
results and formulas needed. There are other examples similar to this
one in the present paper. In this sense, even if we assume that all
axioms for rational conformal field theories hold, the formulas and
results stated in \cite{MS1} and \cite{MS2} are actually conjectures
and in the present paper we give mathematical proofs.

The results of the present paper have been used in \cite{H10}
to prove the 
rigidity and modularity of the tensor category of modules 
for a vertex operator algebra satisfying the conditions 
mentioned above. They have also been used in the constructions
of genus-zero and genus-one full conformal field theories 
in \cite{HK1} and \cite{HK2}.
The results of the present paper have been announced 
in \cite{H9.3} and \cite{H9.7}. See also \cite{Le} for an exposition.

We assume that the reader is familiar with the basic theory of vertex
operator algebras as presented in \cite{FLM}, \cite{FHL} and
\cite{LL}. We also assume that the reader has some basic knowledge in
the theories of intertwining operators, tensor products,
composition-invertible formal series and the Virasoro algebra, and the
modular invariance, as developed in, for example, \cite{DLM},
\cite{HL1}--\cite{HL4}, \cite{H1}--\cite{H9}, \cite{M} and \cite{Z}.

The present paper is organized as follows: In Section 1, we state our
basic assumptions and we discuss intertwining operators and genus-zero
correlation functions constructed from them. An action of $S_{3}$ on
the direct sum of spaces of intertwining operators among irreducible
modules are also given in this section. In section 2, we discuss
geometrically-modified intertwining operators and genus-one
correlation functions constructed from these operators. The modular
transformation associated to $\tau\mapsto -1/\tau$ is recalled in this
section.  Section 3 is devoted to the proof of three formulas for
braiding and fusing matrices. Using all the results obtained in
Sections 1, 2 and 3, we prove two formulas derived first by Moore and
Seiberg from the axioms for rational conformal field theories in
Section 4.  Finally in Section 5, we prove the Verlinde conjecture,
the Verlinde formula for fusion rules and that the matrix associated
to the modular transformation $\tau\mapsto -1/\tau$ is symmetric.

\paragraph{Notations} In this paper, $i$ is either 
$\sqrt{-1}$ or an index, and it should be easy to tell which is
which. The symbols $\mathbb{N}$, $\mathbb{Z}_{+}$, $\mathbb{Z}$,
$\mathbb{Q}$, $\mathbb{C}$, $\C^{\times}$ and $\mathbb{H}$ denote the
nonnegative integers, positive integers, integers, rational numbers,
complex numbers, nonzero complex numbers and the upper half plane,
respectively. We shall use $x,  y, ...$ to denote commuting formal
variables and $z, z_{1}, z_{2}, ...$ to denote complex numbers or
complex variables.

\paragraph{Acknowledgment}
This research is supported in part by NSF grants DMS-0070800 and 
DMS-0401302. I am
grateful to J. Lepowsky,  
L. Kong and the referee for their comments.

\renewcommand{\theequation}{\thesection.\arabic{equation}}
\renewcommand{\thethm}{\thesection.\arabic{thm}}
\setcounter{equation}{0}
\setcounter{thm}{0}

\section{Intertwining operators and genus-zero correlation functions}

Let $V$ be a simple vertex operator algebra and $C_{2}(V)$ the subspace of 
$V$ spanned by $u_{-2}v$ for $u, v\in V$. In the present paper,
we shall always assume that $V$ satisfies the 
following conditions:

\begin{enumerate}

\item $V_{(n)}=0$ for $n<0$, $V_{(0)}=\mathbb{C}\mathbf{1}$ and 
$V'$ is isomorphic to $V$ as a $V$-module.

\item Every $\mathbb{N}$-gradable weak $V$-module is completely 
reducible.

\item $V$ is $C_{2}$-cofinite, that is, $\dim V/C_{2}(V)<\infty$. 

\end{enumerate}

Note that when $V_{(n)}=0$ for $n<0$, $V_{(0)}=\mathbb{C}\mathbf{1}$,
$V'$ is isomorphic to $V$ as a $V$-module if for any irreducible
$V$-module not isomorphic to $V$, $W_{(0)}=0$. So the results of the
present paper still hold if Condition 1 above is replaced by the
condition $V_{(n)}=0$ for $n<0$, $V_{(0)}=\mathbb{C}\mathbf{1}$, and
$W_{(0)}=0$ for any irreducible $V$-module not isomorphic to $V$.
As we have mentioned in the abstract and introduction, 
by results of Li \cite{L} and 
Abe-Buhl-Dong \cite{ABD}, 
Conditions (ii) and (iii) are equivalent to a 
single condition that every weak $V$-module is completely 
reducible.

From \cite{DLM}, we know that there are only finitely many
inequivalent irreducible $V$-modules.  Let $\mathcal{A}$ be the set of
equivalence classes of irreducible $V$-modules. We denote the
equivalence class containing $V$ by $e$. Note that the contragredient
module of an irreducible module is also irreducible (see \cite{FHL}).
So we have a map
\begin{eqnarray*}
^{\prime}:
\mathcal{A}&\to& \mathcal{A}\\
a&\mapsto& a'.
\end{eqnarray*}
For $a\in \mathcal{A}$, if $a\ne a'$, we choose 
representatives $W^{a}$ and $W^{a'}$ of $a$ and $a'$ such that 
$(W^{a})'=W^{a'}$ and, after we identify 
$(W^{a})''$ with $W^{a}$, we also have $(W^{a'})'=W^{a}$. 
If $a=a'\ne e$, we choose any nondegenerate symmetric
bilinear invariant form on $W^{a}$ and using this form, we 
can identify $(W^{a})'$ with $W^{a'}$ (and $(W^{a'})'$ with $W^{a}$).
For $a=e$, we choose the nondegenerate symmetric bilinear invariant form
$(\cdot, \cdot)$ normalized by $(\mathbf{1}, \mathbf{1})$=1.
Since the 
results of the present paper involve only elements of $\mathcal{A}$, 
not representatives of these elements, it is convenient
to identify $V$-modules and
their double contragredient modules and to identify 
$(W^{a})'$ with $W^{a'}$ and $(W^{a'})'$ with $W^{a}$
using the chosen 
nondegenerate bilinear invariant forms. 
After these identifications, we see that 
we can find a representative $W^{a}$ of $a$ for each $a\in \mathcal{A}$
such that $W^{e}=V$ and $(W^{a})'=W^{a'}$. 
In this paper, for simplicity, 
we fix such a choice. 
From \cite{AM} and \cite{DLM}, we know that irreducible $V$-modules 
are in fact graded by rational numbers. 
Thus for $a\in \A$, 
there exists $h_{a}\in \Q$ such that 
$W^{a}=\coprod_{n\in h_{a}+\N}W_{(n)}$ and $W_{(h_{a})}\ne 0$. 

For $a_{1}, a_{2}, a_{3}\in \mathcal{A}$, 
let $\mathcal{V}_{a_{1}a_{2}}^{a_{3}}$ 
be the space of intertwining operators of 
type ${W^{a_{3}}\choose W^{a_{1}}W^{a_{2}}}$. For any 
$\Y\in \mathcal{V}_{a_{1}a_{2}}^{a_{3}}$, we know from \cite{FHL}
that 
for $w_{a_{1}}\in W^{a_{1}}$ and $w_{a_{2}}\in W^{a_{2}}$
\begin{equation}\label{y-delta}
\Y(w_{a_{1}}, x)w_{a_{2}}\in x^{\Delta(\Y)}W^{a_{3}}[[x, x^{-1}]],
\end{equation}
where 
$$\Delta(\mathcal{Y})=h_{a_{3}}-h_{a_{1}}-h_{a_{2}}.$$

In the present paper, we shall use the following conventions: For $z\in 
\C^{\times}$, $\log z=\log |z|+i\arg z$, where $0\le \arg z<2\pi$. 
For $z\in \C^{\times}$ and $r\in \C$, $z^{r}=e^{r\log z}$. 
For  $n\in \Z$ and $s\in \C$, 
$(e^{n\pi i})^{s}=e^{sn\pi i}$.  Similarly,
for $z\in \C^{\times}$ and $O$ an linear operator, $z^{O}=e^{(\log z) O}$,
if $e^{(\log z) O}$ is well-defined. For $n\in \Z$ and 
$O$ a linear operator, $(e^{n\pi i})^{O}=e^{n\pi iO}$, if $e^{n\pi iO}$ 
is well
defined. For $a_{1}, a_{2}, a_{3}\in \A$,
$z\in \C^{\times}$, $r\in \C$, $\Y\in \mathcal{V}_{a_{1}a_{2}}^{a_{3}}$, 
$w_{a_{1}}\in W^{a_{1}}$ and $w_{a_{2}}\in W^{a_{2}}$, 
$\Y(w_{a_{1}}, z)w_{a_{2}}$
is $\Y(w_{a_{1}}, x)w_{a_{2}}|_{x^{n}=e^{n\log z}, \;
n\in \C}$
and 
$\Y(w_{a_{1}}, e^{r}x)w_{a_{2}}$ is $\Y(w_{a_{1}}, y)w_{a_{2}}
|_{y^{n}=e^{nr}x^{n}, \;
n\in \C}$ by the conventions above.

We also know that the fusion rules 
$N_{a_{1}a_{2}}^{a_{3}}
=N_{W^{a_{1}}W^{a_{2}}}^{W^{a_{3}}}$ for 
$a_{1}, a_{2}, a_{3}\in \mathcal{A}$ are all finite (see \cite{GN},
\cite{L},
\cite{AN},
\cite{H8}). 
For $a_{1}, a_{2}, a_{3}\in \mathcal{A}$, we have isomorphisms
$\Omega_{-r}: \mathcal{V}_{a_{1}a_{2}}^{a_{3}} \to 
\mathcal{V}_{a_{2}a_{1}}^{a_{3}}$ and 
$A_{-r}: \mathcal{V}_{a_{1}a_{2}}^{a_{3}} \to 
\mathcal{V}_{a_{1}a'_{3}}^{a'_{2}}$ for $r\in \Z$, defined
in Section 7 of \cite{HL2}.
Using these isomorphisms, we define a left action of 
the symmetric group $S_{3}$ on 
$$\mathcal{V}=\coprod_{a_{1}, a_{2}, a_{3}\in 
\mathcal{A}}\mathcal{V}_{a_{1}a_{2}}^{a_{3}}$$
as follows:
For $a_{1}, a_{2}, a_{3}\in \A$,
$\mathcal{Y}\in \mathcal{V}_{a_{1}a_{2}}^{a_{3}}$, 
we define 
\begin{eqnarray*}
\sigma_{12}(\mathcal{Y})&=&e^{\pi i \Delta(\mathcal{Y})}
\Omega_{-1}(\mathcal{Y})\\
&=&e^{-\pi i \Delta(\mathcal{Y})}\Omega_{0}(\mathcal{Y}),\\
\sigma_{23}(\mathcal{Y})&=&e^{\pi i h_{a_{1}}}
A_{-1}(\mathcal{Y})\\
&=&e^{-\pi i h_{a_{1}}}A_{0}(\mathcal{Y}).
\end{eqnarray*}

\begin{prop}
The actions $\sigma_{12}$ and $\sigma_{23}$ of $(12)$ and 
$(23)$ on $\mathcal{V}$ defined above
generate an action of $S_{3}$ on $\mathcal{V}$. 
\end{prop}
\pf
We need only prove the relations which must be satisfied by 
$\sigma_{12}$ and $\sigma_{23}$. Propositions 7.1 and 7.3
in \cite{HL2} are actually equivalent to $\sigma_{12}^{2}=1$ and 
$\sigma_{23}^{2}=1$, respectively. So the only relation we need to 
prove is
\begin{equation}\label{s3-sym}
\sigma_{12}\sigma_{23}\sigma_{12}\sigma_{23}=\sigma_{23}\sigma_{12}.
\end{equation}

Let $\mathcal{Y}$ be an element of $\mathcal{V}_{a_{1}a_{2}}^{a_{3}}$.
Then for $w_{a_{1}}\in W^{a_{1}}$, $w_{a_{2}}\in W_{a_{2}}$ and 
$w'_{a_{3}}\in W_{a'_{3}}$, 
\begin{eqnarray}\label{s3-sym-r}
\lefteqn{\langle \sigma_{23}(\sigma_{12}(\Y))(w_{a_{2}},
x)w'_{a_{3}}, w_{a_{1}}\rangle}\nn
&&=e^{\pi ih_{a_{2}}}\langle  w'_{a_{3}}, \sigma_{12}(\Y)(e^{xL(1)}
(e^{-\pi i}x^{-2})^{L(0)}w_{a_{2}},
x^{-1})w_{a_{1}}\rangle\nn
&&=e^{\pi i(h_{a_{3}}-h_{a_{1}})}\langle  w'_{a_{3}}, e^{x^{-1}L(-1)}
\Y(w_{a_{1}}, e^{-\pi i}x^{-1})e^{xL(1)}
(e^{-\pi i}x^{-2})^{L(0)}w_{a_{2}}\rangle.\nn
&&
\end{eqnarray}
On the other hand,
\begin{eqnarray}\label{s3-sym-l}
\lefteqn{\langle 
\sigma_{12}(\sigma_{23}(\sigma_{12}(\sigma_{23}(\Y))))(w_{a_{2}},
x)w'_{a_{3}}, w_{a_{1}}\rangle}\nn
&&=e^{\pi i (h_{a_{1}}-h_{a_{2}}-h_{a_{3}})}
\langle e^{xL(-1)} \sigma_{23}(\sigma_{12}(\sigma_{23}(\Y)))(w'_{a_{3}},
e^{-\pi i}x)w_{a_{2}}, w_{a_{1}}\rangle\nn
&&=e^{\pi i (h_{a_{1}}-h_{a_{2}})}\cdot\nn
&&\quad\quad\cdot
\langle w_{a_{2}}, \sigma_{12}(\sigma_{23}(\Y))(e^{-xL(1)}
(e^{-\pi i}(e^{-\pi i}x)^{-2})^{L(0)}w'_{a_{3}},
e^{\pi i}x^{-1})e^{xL(1)}w_{a_{1}}\rangle\nn
&&=e^{-\pi i h_{a_{3}}}
\langle w_{a_{2}}, e^{-x^{-1}L(-1)}\sigma_{23}(\Y)(e^{xL(1)}w_{a_{1}},
x^{-1})e^{-xL(1)}
(e^{\pi i}x^{-2})^{L(0)}w'_{a_{3}}\rangle\nn
&&=e^{\pi i (h_{a_{1}}-h_{a_{3}})}
\langle \Y(e^{x^{-1}L(1)}(e^{-\pi i}(x^{-1})^{-2})^{L(0)}
e^{xL(1)}w_{a_{1}},
x)e^{-x^{-1}L(1)}w_{a_{2}}, \nn
&&\quad\quad\quad\quad\quad\quad\quad\quad\quad\quad\quad\quad
\quad\quad\quad\quad\quad\quad\quad\quad\quad
e^{-xL(1)}
(e^{\pi i}x^{-2})^{L(0)}w'_{a_{3}}\rangle\nn
&&=e^{\pi i (h_{a_{1}}-h_{a_{3}})}
\langle  w'_{a_{3}}, (e^{\pi i}x^{-2})^{L(0)}e^{-xL(-1)}\cdot\nn
&&\quad\quad\quad\quad\quad\quad\quad\quad\cdot
\Y(e^{x^{-1}L(1)}(e^{-\pi i}x^{2})^{L(0)}
e^{xL(1)}w_{a_{1}},
x)e^{-x^{-1}L(1)}w_{a_{2}}\rangle\nn
&&=e^{\pi i (h_{a_{1}}-h_{a_{3}})}
\langle  w'_{a_{3}}, e^{x^{-1}L(-1)}(e^{\pi i}x^{-2})^{L(0)}\cdot\nn
&&\quad\quad\quad\quad\quad\quad\quad\quad\cdot
\Y(e^{x^{-1}L(1)}(e^{-\pi i}x^{2})^{L(0)}
e^{xL(1)}w_{a_{1}},
x)e^{-x^{-1}L(1)}w_{a_{2}}\rangle\nn
&&=e^{\pi i (h_{a_{1}}-h_{a_{3}})}
\langle  w'_{a_{3}}, e^{x^{-1}L(-1)}\cdot\nn
&&\quad\quad\quad\quad\quad\quad\quad\quad\cdot
\Y((e^{\pi i}x^{-2})^{L(0)}
e^{x^{-1}L(1)}(e^{-\pi i}x^{2})^{L(0)}
e^{xL(1)}w_{a_{1}},
e^{\pi i}x^{-1})\cdot\nn
&&\quad\quad\quad\quad\quad\quad\quad\quad\quad\quad\quad\quad\quad\quad
\quad\quad\quad\quad\quad\cdot
(e^{\pi i}x^{-2})^{L(0)}e^{-x^{-1}L(1)}w_{a_{2}}\rangle\nn
&&=e^{\pi i (h_{a_{3}}-h_{a_{1}})}
\langle  w'_{a_{3}}, e^{x^{-1}L(-1)}
\Y(w_{a_{1}},
e^{-\pi i}x^{-1})e^{xL(1)}(e^{-\pi i}x^{-2})^{L(0)}w_{a_{2}}\rangle.
\end{eqnarray}
From (\ref{s3-sym-r}) and (\ref{s3-sym-l}), we obtain
$$\langle \sigma_{23}(\sigma_{12}(\Y))(w_{a_{2}},
x)w'_{a_{3}}, w_{a_{1}}\rangle=\langle 
\sigma_{12}(\sigma_{23}(\sigma_{12}(\sigma_{23}(\Y))))(w_{a_{2}},
x)w'_{a_{3}}, w_{a_{1}}\rangle.$$
Since $a_{1}, a_{2}, a_{3}$, $w_{a_{1}}, w_{a_{2}}, w_{a_{3}}'$ and $\Y$ are
arbitrary, we see that (\ref{s3-sym}) holds.
\epfv

For $p=1, \dots, 6$ and $a_{1}, a_{2}, 
a_{3}\in \mathcal{A}$, let 
$\Y_{a_{1}a_{2}; i}^{a_{3}; (p)}$, $i=1, \dots, N_{a_{1}a_{2}}^{a_{3}}$, be
bases of $\mathcal{V}_{a_{1}a_{2}}^{a_{3}}$.
From Thereom 3.9 in \cite{H8}, 
the space $\coprod_{a\in \mathcal{A}}W^{a}$
has a natural structure of an intertwining operator algebra
in the sense of \cite{H5} and \cite{H7}.
In particular, we have the associativity of intertwining operators
(see (14.55) in \cite{H1}, Condition (vii) in Definition 3.1 
in \cite{H5} and Axiom 3 in Definition 2.1 in \cite{H7}).
Thus there exist 
$$F(\Y_{a_{1}a_{5}; i}^{a_{4}; (1)}\otimes \Y_{a_{2}a_{3}; j}^{a_{5}; (2)}; 
\Y_{a_{6}a_{3}; l}^{a_{4}; (3)}
\otimes \Y_{a_{1}a_{2}; k}^{a_{6}; (4)}) \in \C$$
for $a_{1}, \dots, a_{6}\in \mathcal{A}$, $i=1, \dots, 
N_{a_{1}a_{5}}^{a_{4}}$, $j=1, \dots, 
N_{a_{2}a_{3}}^{a_{5}}$, $k=1, \dots, 
N_{a_{1}a_{2}}^{a_{6}}$, $l=1, \dots, 
N_{a_{6}a_{3}}^{a_{4}}$, 
such that 
\begin{eqnarray}\label{assoc}
\lefteqn{\langle w_{a'_{4}}, \Y_{a_{1}a_{5}; i}^{a_{4}; (1)}(w_{a_{1}}, z_{1})
\Y_{a_{2}a_{3}; j}^{a_{5}; (2)}(w_{a_{2}}, z_{2})w_{a_{3}}\rangle}\nn
&&=\sum_{a_{6}\in \A}\sum_{k=1}^{N_{a_{1}a_{2}}^{a_{6}}}
\sum_{l=1}^{N_{a_{6}a_{3}}^{a_{4}}}
F(\Y_{a_{1}a_{5}; i}^{a_{4}; (1)}\otimes \Y_{a_{2}a_{3}; j}^{a_{5}; (2)}; 
\Y_{a_{6}a_{3}; l}^{a_{4}; (3)}
\otimes \Y_{a_{1}a_{2}; k}^{a_{6}; (4)})\cdot\nn
&&\quad\quad\quad\quad\cdot 
\langle w_{a'_{4}}, 
\Y_{a_{6}a_{3}; l}^{a_{4}; (3)}(\Y_{a_{1}a_{2}; k}^{a_{6}; (4)}(w_{a_{1}}, z_{1}-z_{2})
w_{a_{2}}, z_{2})w_{a_{3}}\rangle
\end{eqnarray}
when $|z_{1}|>|z_{2}|>|z_{1}-z_{2}|>0$, 
for $a_{1}, \dots, a_{5}\in \A$, $w_{a_{1}}\in W^{a_{1}}$, 
$w_{a_{2}}\in W^{a_{2}}$, $w_{a_{3}}\in W^{a_{3}}$,  
$w_{a'_{4}}\in W^{a'_{4}}=(W^{a_{4}})'$, $i=1, \dots, 
N_{a_{1}a_{5}}^{a_{4}}$ and $j=1, \dots, 
N_{a_{2}a_{3}}^{a_{5}}$. Note that here we have used our convention
on the choices of values of intertwining operators. The numbers 
$$F(\Y_{a_{1}a_{5}; i}^{a_{4}; (1)}\otimes \Y_{a_{2}a_{3}; j}^{a_{5}; (2)}; 
\Y_{a_{6}a_{3}; l}^{a_{4}; (3)}
\otimes \Y_{a_{1}a_{2}; k}^{a_{6}; (4)})$$
are matrix elements of the fusing isomorphism.
In \cite{H7}, when these four bases are chosen to be the same, 
these matrix elements are denoted as $F_{a_{5}; a_{6}}^{i, j; k, l}(a_{1}, 
a_{2}, a_{3}; a_{4})$.
In the present paper, since we want to emphasize their dependence on the bases
and since we shall need matrix elements of the fusing 
isomorphism under different bases, we
have to instead use the notation above. 

The fusing isomorphism is invertible. Thus for any bases as above, 
there exist 
$$F^{-1}(\Y_{a_{6}a_{3}; l}^{a_{4}; (1)}
\otimes \Y_{a_{1}a_{2}; k}^{a_{6}; (2)};
\Y_{a_{1}a_{5}; i}^{a_{4}; (3)}\otimes \Y_{a_{2}a_{3}; j}^{a_{5}; (4)}) \in \C$$
for $a_{1}, \dots, a_{6}\in \mathcal{A}$, $i=1, \dots, 
N_{a_{1}a_{5}}^{a_{4}}$, $j=1, \dots, 
N_{a_{2}a_{3}}^{a_{5}}$, $k=1, \dots, 
N_{a_{1}a_{2}}^{a_{6}}$, $l=1, \dots, 
N_{a_{6}a_{3}}^{a_{4}}$, 
such that 
\begin{eqnarray}\label{assoc-inv}
\lefteqn{\langle w_{a'_{4}}, 
\Y_{a_{6}a_{3}; l}^{a_{4}; (1)}(\Y_{a_{1}a_{2}; k}^{a_{6}; (2)}(w_{a_{1}}, z_{1}-z_{2})
w_{a_{2}}, z_{2})w_{a_{3}}\rangle}\nn
&&=\sum_{a_{5}\in \A}\sum_{i=1}^{N_{a_{1}a_{5}}^{a_{4}}}
\sum_{j=1}^{N_{a_{2}a_{3}}^{a_{5}}}
F^{-1}(\Y_{a_{6}a_{3}; l}^{a_{4}; (1)}
\otimes \Y_{a_{1}a_{2}; k}^{a_{6}; (2)};
\Y_{a_{1}a_{5}; i}^{a_{4}; (3)}\otimes \Y_{a_{2}a_{3}; j}^{a_{5}; (4)})\cdot\nn
&&\quad\quad\quad\quad\cdot 
\langle w_{a'_{4}}, \Y_{a_{1}a_{5}; i}^{a_{4}; (3)}(w_{a_{1}}, z_{1})
\Y_{a_{2}a_{3}; j}^{a_{5}; (4)}(w_{a_{2}}, z_{2})w_{a_{3}}\rangle
\end{eqnarray}
when $|z_{1}|>|z_{2}|>|z_{1}-z_{2}|>0$,
for $a_{1}, \dots, a_{4}, a_{6}\in \A$, $w_{a_{1}}\in W^{a_{1}}$, 
$w_{a_{2}}\in W^{a_{2}}$, $w_{a_{3}}\in W^{a_{3}}$,  
$w_{a'_{4}}\in W^{a'_{4}}=(W^{a_{4}})'$, $k=1, \dots, 
N_{a_{1}a_{2}}^{a_{6}}$ and $l=1, \dots, 
N_{a_{6}a_{3}}^{a_{4}}$. These numbers are matrix elements of the inverse
of the fusing isomorphism.

By Lemma 4.1 in \cite{H7} or by the differential equations
given by Theorem 1.4 in \cite{H8}, we know that 
\begin{eqnarray}
&\langle w_{a'_{4}}, \Y_{a_{1}a_{5}; i}^{a_{4}; (1)}(w_{a_{1}}, z_{1})
\Y_{a_{2}a_{3}; j}^{a_{5}; (2)}(w_{a_{2}}, z_{2})w_{a_{3}}\rangle, 
&\label{product}\\
&\langle w_{a'_{4}}, 
\Y_{a_{6}a_{3}; l}^{a_{4}; (3)}
(\Y_{a_{1}a_{2}; k}^{a_{6}; (4)}(w_{a_{1}}, z_{1}-z_{2})
w_{a_{2}}, z_{2})w_{a_{3}}\rangle
\rangle, &\label{iterate}\\
&\langle w_{a'_{4}}, 
\Y_{a_{2}a_{6}; k}^{a_{4}; (5)}(
w_{a_{2}}, z_{2})\Y_{a_{1}a_{3}; l}^{a_{6}; (6)}(w_{a_{1}}, z_{1})
w_{a_{3}}\rangle&\label{product2}
\end{eqnarray}
can all be analytically extended to multi-valued analytic 
functions on 
$$M^{2}=\{(z_{1}, z_{2})\in \C^{2}\;|\; z_{1}, z_{2}\ne 0, z_{1}\ne z_{2}\}.$$
We can lift  multi-valued analytic functions on $M^{2}$ to single-valued 
analytic functions on the universal covering $\widetilde{M}^{2}$ of 
$M^{2}$. However, these single-valued liftings are not unique. To 
obtain a unique lifting, we have to give single-valued branches
of the multi-valued analytic functions in simply-connected
subsets of $M^{2}$.

Consider the following three regions:
\begin{eqnarray}
M^{2}(|z_{1}|>|z_{2}|>0, \; 0\le \arg z_{1}, \arg z_{2}<2\pi);\label{re1}\\
M^{2}(|z_{2}|>|z_{1}|>0, \; 0\le \arg z_{1}, \arg z_{2}<2\pi);\label{re2}\\
M^{2}(|z_{2}|>|z_{1}-z_{2}|>0, \; 0\le \arg z_{2}, \arg (z_{1}-z_{2})<2\pi).
\label{re3}
\end{eqnarray}
These are, respectively, the region $|z_{1}|>|z_{2}|>0$ with cuts 
along the real lines in the $z_{1}$- and $z_{2}$-planes; the
region $|z_{2}|>|z_{1}|>0$ with cuts 
along the real lines in the $z_{1}$- and $z_{2}$-planes;
the region $|z_{2}|>|z_{1}-z_{2}|>0$ with cuts 
along the real lines in the $z_{2}$- and $z_{1}-z_{2}$-planes. 
The  multi-valued analytic extensions of (\ref{product}) and
(\ref{iterate}) have single-valued branches 
(\ref{product}) and (\ref{iterate}), respectively, 
on (\ref{re1}) and (\ref{re3}), respectively. 
Thus we have the corresponding single-valued
analytic functions on $\widetilde{M}^{2}$. We denote these 
analytic functions by 
$$E(\langle w_{a'_{4}}, \Y_{a_{1}a_{5}; i}^{a_{4}; (1)}(w_{a_{1}}, x_{1})
\Y_{a_{2}a_{3}; j}^{a_{5}; (2)}(w_{a_{2}}, x_{2})w_{a_{3}}\rangle)$$
and
$$E(\langle w_{a'_{4}}, 
\Y_{a_{6}a_{3}; l}^{a_{4}; (3)}(\Y_{a_{1}a_{2}; k}^{a_{6}; (4)}(w_{a_{1}}, z_{1}-z_{2})
w_{a_{2}}, z_{2})w_{a_{3}}\rangle),$$
respectively.

Then from (\ref{assoc}) and (\ref{assoc-inv}), we 
immediately obtain
\begin{eqnarray}\label{assoc2}
\lefteqn{E(\langle w_{a'_{4}}, \Y_{a_{1}a_{5}; i}^{a_{4}; (1)}(w_{a_{1}}, z_{1})
\Y_{a_{2}a_{3}; j}^{a_{5}; (2)}(w_{a_{2}}, z_{2})w_{a_{3}}\rangle)}\nn
&&=\sum_{a_{6}\in \A}\sum_{k=1}^{N_{a_{1}a_{2}}^{a_{6}}}
\sum_{l=1}^{N_{a_{6}a_{3}}^{a_{4}}}
F(\Y_{a_{1}a_{5}; i}^{a_{4}; (1)}\otimes \Y_{a_{2}a_{3}; j}^{a_{5}; (2)}; 
\Y_{a_{6}a_{3}; l}^{a_{4}; (3)}
\otimes \Y_{a_{1}a_{2}; k}^{a_{6}; (4)})\cdot\nn
&&\quad\quad\quad\quad\cdot 
E(\langle w_{a'_{4}}, 
\Y_{a_{6}a_{3}; l}^{a_{4}; (3)}(\Y_{a_{1}a_{2}; k}^{a_{6}; (4)}(w_{a_{1}}, z_{1}-z_{2})
w_{a_{2}}, z_{2})w_{a_{3}}\rangle)
\end{eqnarray}
and
\begin{eqnarray}\label{assoc-inv2}
\lefteqn{E(\langle w_{a'_{4}}, 
\Y_{a_{6}a_{3}; l}^{a_{4}; (1)}(\Y_{a_{1}a_{2}; k}^{a_{6}; (2)}(w_{a_{1}}, z_{1}-z_{2})
w_{a_{2}}, z_{2})w_{a_{3}}\rangle)}\nn
&&=\sum_{a_{5}\in \A}\sum_{i=1}^{N_{a_{1}a_{5}}^{a_{4}}}
\sum_{j=1}^{N_{a_{2}a_{3}}^{a_{5}}}
F^{-1}(\Y_{a_{6}a_{3}; l}^{a_{4}; (1)}
\otimes \Y_{a_{1}a_{2}; k}^{a_{6}; (2)};
\Y_{a_{1}a_{5}; i}^{a_{4}; (3)}\otimes \Y_{a_{2}a_{3}; j}^{a_{5}; (4)})\cdot\nn
&&\quad\quad\quad\quad\cdot 
E(\langle w_{a'_{4}}, \Y_{a_{1}a_{5}; i}^{a_{4}; (3)}(w_{a_{1}}, z_{1})
\Y_{a_{2}a_{3}; j}^{a_{5}; (4)}(w_{a_{2}}, z_{2})w_{a_{3}}\rangle).
\end{eqnarray}

From a single-valued analytic function $f$ on $\widetilde{M}^{2}$,
we can construct other single-valued analytic functions on $\widetilde{M}^{2}$
such that they and $f$ correspond to the same multi-valued analytic
functions on $M^{2}$. 
One particularly interesting construction is given by a braiding
operation. 
To give this construction, we need only  give a single-valued 
branch of the multi-valued analytic function on $M^{2}$ corresponding to 
$f$
on a simply-connected region in $M^{2}$.
We can view the region (\ref{re2}) as a simply-connected 
region in $\widetilde{M}^{2}$. Thus $f$ gives a single-valued function
on this region. For any $r\in \Z$, we now give a single-valued analytic 
functions on the region (\ref{re1}) in the 
following way: Consider the path 
$$t \mapsto \left(\frac{3}{2}
-\frac{e^{(2r+1)\pi i t}}{2}, \frac{3}{2}
+\frac{e^{(2r+1)\pi i t}}{2}\right)$$ 
from the point $(1, 2)$ in the region (\ref{re2})
to the point $(2, 1)$ in the region (\ref{re1}).
This path and the restriction of $f$ to the region (\ref{re2})
gives a unique single-valued analytic function on the 
region (\ref{re1}). This single-valued analytic function
on this simply-connected region determines uniquely a
single-valued analytic function $B^{(r)}(f)$ on $\widetilde{M}^{2}$
such that $f$ and $B^{(r)}(f)$ correspond to the same 
multi-valued analytic function on $M^{2}$.

We now consider (\ref{product2}). From the discussion above, we 
have a single-valued analytic function 
$$E(\langle w_{a'_{4}}, 
\Y_{a_{2}a_{6}; k}^{a_{4}; (5)}(
w_{a_{2}}, z_{2})\Y_{a_{1}a_{3}; l}^{a_{6}; (6)}(w_{a_{1}}, z_{1})
w_{a_{3}}\rangle)$$
on $\widetilde{M}^{2}$. 
Apply our construction above, we obtain another 
single-valued analytic function 
$$B^{(r)}(E(\langle w_{a'_{4}}, 
\Y_{a_{2}a_{6}; k}^{a_{4}; (3)}(
w_{a_{2}}, z_{2})\Y_{a_{1}a_{3}; l}^{a_{6}; (4)}(w_{a_{1}}, z_{1})
w_{a_{3}}\rangle))$$
on $\widetilde{M}^{2}$. 
In terms of these functions, the commutativity for intertwining 
operators (see Theorem 3.1 in \cite{H3} and Proposition 2.2
\cite{H7}) can be written as 
\begin{eqnarray*}
\lefteqn{B^{(r)}(E(\langle w_{a'_{4}}, \Y_{a_{1}a_{5}; i}^{a_{4}; (1)}(w_{a_{1}}, z_{2})
\Y_{a_{2}a_{3}; j}^{a_{5}; (2)}(w_{a_{2}}, z_{1})w_{a_{3}}\rangle))}\nn
&&=\sum_{a_{6}\in \A}
\sum_{k=1}^{N_{a_{2}a_{6}}^{a_{4}}}\sum_{l=1}^{N_{a_{1}a_{3}}^{a_{6}}}
B^{(r)}(\Y_{a_{1}a_{5}; i}^{a_{4}; (1)}\otimes \Y_{a_{2}a_{3}; j}^{a_{5}; (2)};
\Y_{a_{2}a_{6}; k}^{a_{4}; (3)}\otimes \Y_{a_{1}a_{3}; l}^{a_{6}; (4)})\cdot\nn
&&\quad\quad\quad\quad\cdot 
E(\langle w_{a'_{4}}, 
\Y_{a_{2}a_{6}; k}^{a_{4}; (3)}(
w_{a_{2}}, z_{1})\Y_{a_{1}a_{3}; l}^{a_{6}; (4)}(w_{a_{1}}, z_{2})
w_{a_{3}}\rangle),
\end{eqnarray*}
where  the numbers 
$$B^{(r)}(\Y_{a_{1}a_{5}; i}^{a_{4}; (1)}\otimes \Y_{a_{2}a_{3}; j}^{a_{5}; (2)};
\Y_{a_{2}a_{6}; k}^{a_{4}; (3)}\otimes \Y_{a_{1}a_{3}; l}^{a_{6}; (4)})$$
are the matrix elements of the braiding isomorphism. 

From the construction above, it is easy to see that the square 
$(B^{(r)})^{2}$ of $B^{(r)}$ is actually 
the map sending a single-valued analytic function $f$ on $\widetilde{M}^{2}$
to another one $(B^{(r)})^{2}(f)$ 
corresponding to the same multi-valued analytic functions 
on $M^{2}$ in the following way: Consider the path 
$$t \mapsto \left(\frac{3}{2}
+\frac{e^{2(2r+1)\pi i t}}{2}, \frac{3}{2}
-\frac{e^{2(2r+1)\pi i t}}{2}\right)$$ 
from the point $(2, 1)$ to itself.
This path and the restriction of $f$ to the region (\ref{re1})
gives another single-valued analytic function on the same
region. This new single-valued analytic function on the region
(\ref{re1}) gives a single-valued analytic 
function $(B^{(r)})^{2}(f)$ on $\widetilde{M}^{2}$. From this
description of $(B^{(r)})^{2}(f)$, we see 
that $(B^{(r)})^{2}(f)$ is in fact the 
monodromy of the multi-valued function corresponding to 
$f$ given by $\log (z_{1}-z_{2})\mapsto \log (z_{1}-z_{2}) +2(2r+1)\pi i$.

We shall also use similar notations to denote the matrix elements of
the square $(B^{(r)})^{2}$ of $B^{(r)}$ under the bases above as 
$$(B^{(r)})^{2}(\Y_{a_{1}a_{5}; i}^{a_{4}; (1)}\otimes \Y_{a_{2}a_{3}; j}^{a_{5}; (2)};
\Y_{a_{1}a_{6}; i}^{a_{4}; (3)}\otimes \Y_{a_{2}a_{3}; j}^{a_{6}; (4)}),$$
that is,
\begin{eqnarray}\label{b-r-2}
\lefteqn{(B^{(r)})^{2}(E(\langle w_{a'_{4}}, \Y_{a_{1}a_{5}; i}^{a_{4}; (1)}(w_{a_{1}}, z_{1})
\Y_{a_{2}a_{3}; j}^{a_{5}; (2)}(w_{a_{2}}, z_{2})w_{a_{3}}\rangle))}\nn
&&=\sum_{a_{6}\in \A}
\sum_{k=1}^{N_{a_{2}a_{6}}^{a_{4}}}\sum_{l=1}^{N_{a_{1}a_{3}}^{a_{6}}}
(B^{(r)})^{2}(\Y_{a_{1}a_{5}; i}^{a_{4}; (1)}\otimes \Y_{a_{2}a_{3}; j}^{a_{5}; (2)};
\Y_{a_{1}a_{6}; i}^{a_{4}; (3)}\otimes \Y_{a_{2}a_{3}; j}^{a_{6}; (4)})\cdot\nn
&&\quad\quad\quad\quad\cdot 
E(\langle w_{a'_{4}}, \Y_{a_{1}a_{6}; i}^{a_{4}; (1)}(w_{a_{1}}, z_{1})
\Y_{a_{2}a_{3}; j}^{a_{6}; (2)}(w_{a_{2}}, z_{2})w_{a_{3}}\rangle).
\end{eqnarray}

As in the case of correlation functions constructed from 
products and iterates of intertwining operators above, if we have a
series $\varphi$ which is absolutely convergent in a region in $M^{2}$ and 
can be analytically extended to a multi-valued analytic function on $M^{2}$,
then we obtain a unique single-valued analytic 
function $E(\varphi)$ on $\widetilde{M}^{2}$. For example, we know that 
for $z\in \C$,
$$\langle w_{a'_{4}}, \Y_{a_{1}a_{5}; i}^{a_{4}; (1)}(w_{a_{1}}, z_{1})
e^{zL(-1)}\Y_{a_{2}a_{3}; j}^{a_{5}; (2)}(w_{a_{2}}, z_{2}-z)
e^{-zL(-1)}w_{a_{3}}\rangle$$
is absolutely convergent in the region given by $|z_{1}|>|z_{2}|>0$ and
$|z_{2}-z|>|z|$ and it can be 
analytically extended to a multi-valued analytic function on $M^{2}$. 
Thus we have a single-valued analytic function
$$E(\langle w_{a'_{4}}, \Y_{a_{1}a_{5}; i}^{a_{4}; (1)}(w_{a_{1}}, z_{1})
e^{zL(-1)}\Y_{a_{2}a_{3}; j}^{a_{5}; (2)}(w_{a_{2}}, z_{2}-z)
e^{-zL(-1)}w_{a_{3}}\rangle)$$
on $\widetilde{M}^{2}$. All the other properties of 
intertwining operators, for example,
the skew-symmetry, the contragredient intertwining operators and so on,
can all be expressed using equalities for such single-valued  analytic functions 
on $\widetilde{M}^{2}$. 

We shall need the following result:

\begin{prop}\label{independence0}
For $a_{1}, a_{2}, a_{3}, a_{4}\in \A$, the maps from 
$W^{a'_{4}}\otimes W^{a_{1}}\otimes W^{a_{2}}\otimes W^{a_{3}}$
to the space of single-valued analytic functions on $\widetilde{M}^{2}$
given by 
$$w_{a'_{4}}\otimes w_{a_{1}} \otimes w_{a_{2}}\otimes w_{a_{3}}
\mapsto
E(\langle w_{a'_{4}}, \Y_{a_{1}a_{5}; i}^{a_{4}; (1)}(w_{a_{1}}, z_{1})
\Y_{a_{2}a_{3}; j}^{a_{5}; (2)}(w_{a_{2}}, z_{2})w_{a_{3}}\rangle),$$
$a_{5}\in \A$, $i=1, \dots, N_{a_{1}a_{5}}^{a_{4}}$, 
$j=1, \dots, N_{a_{2}a_{3}}^{a_{5}}$, are linearly independent. 
Similarly, for $a_{1}, a_{2}, a_{3}, a_{4}\in \A$, the maps from 
$W^{a'_{4}}\otimes W^{a_{1}}\otimes W^{a_{2}}\otimes W^{a_{3}}$
to the space of single-valued analytic functions on $\widetilde{M}^{2}$
given by 
$$w_{a'_{4}}\otimes w_{a_{1}} \otimes w_{a_{2}}\otimes w_{a_{3}}
\mapsto
E(\langle w_{a'_{4}}, \Y_{a_{6}a_{3}; k}^{a_{4}; (3)}(
\Y_{a_{1}a_{2}; l}^{a_{6}; (4)}(w_{a_{1}}, z_{1}-z_{2})
w_{a_{2}}, z_{2})w_{a_{3}}\rangle),$$
$a_{6}\in \A$, $k=1, \dots, N_{a_{6}a_{3}}^{a_{4}}$,
$l=1, \dots, N_{a_{1}a_{2}}^{a_{6}}$, are linearly independent. 
\end{prop}
\pf
We prove only the linear independence of 
the maps obtained
from products of intertwining operators. For iterates, the proof 
is similar.

Since analytic extensions are unique, we need only prove that 
the maps from 
$W^{a'_{4}}\otimes W^{a_{1}}\otimes W^{a_{2}}\otimes W^{a_{3}}$
to the space of the single-valued analytic functions 
on the region (\ref{re1}) given by 
$$w_{a'_{4}}\otimes w_{a_{1}} \otimes w_{a_{2}}\otimes w_{a_{3}}
\mapsto
\langle w_{a'_{4}}, \Y_{a_{1}a_{5}; i}^{a_{4}; (1)}(w_{a_{1}}, z_{1})
\Y_{a_{2}a_{3}; j}^{a_{5}; (2)}(w_{a_{2}}, z_{2})w_{a_{3}},
\rangle, $$
$a_{5}\in \A$, $i=1, \dots, N_{a_{1}a_{5}}^{a_{4}}$, 
$j=1, \dots, N_{a_{2}a_{3}}^{a_{5}}$, are linearly independent.
Assume that 
\begin{eqnarray}\label{independence}
\sum_{a_{5}\in \A}\sum_{i=1}^{N_{a_{1}a_{5}}^{a_{4}}}
\sum_{j=1}^{N_{a_{2}a_{3}}^{a_{5}}}\lambda_{a_{5}, i, j} 
\langle w_{a'_{4}}, \Y_{a_{1}a_{5}; i}^{a_{4}; (1)}(w_{a_{1}}, z_{1})
\Y_{a_{2}a_{3}; j}^{a_{5}; (2)}(w_{a_{2}}, z_{2})w_{a_{3}}\rangle=0
\end{eqnarray}
Since (\ref{independence}) holds for all $z_{1}$ and $z_{2}$ 
satisfying $|z_{1}|>|z_{2}|>0$, we obtain the following 
equation in formal variables:
\begin{eqnarray}\label{independence1}
\sum_{a_{5}\in \A}\sum_{i=1}^{N_{a_{1}a_{5}}^{a_{4}}}
\sum_{j=1}^{N_{a_{2}a_{3}}^{a_{5}}}\lambda_{a_{5}, i, j} 
\langle w_{a'_{4}}, \Y_{a_{1}a_{5}; i}^{a_{4}; (1)}(w_{a_{1}}, x_{1})
\Y_{a_{2}a_{3}; j}^{a_{5}; (2)}(w_{a_{2}}, x_{2})w_{a_{3}}\rangle=0
\end{eqnarray}
We want to show that $\lambda_{a_{5}, i, j}=0$ for 
$a_{5}\in \A$, $i=1, \dots, N_{a_{1}a_{5}}^{a_{4}}$ and 
$j=1, \dots, N_{a_{2}a_{3}}^{a_{5}}$. 

From the tensor product theory in \cite{HL1} and \cite{HL4}, we know that 
the tensor product module $W^{a_{2}}\boxtimes_{P(z_{2})} W^{a_{3}}$ is isomorphic to 
$\oplus_{a_{5}\in \A}N_{a_{2}a_{3}}^{a_{5}}W^{a_{5}}$. 
For $a_{5}\in \A$ and $j=1, \dots, N_{a_{2}a_{3}}^{a_{5}}$,
let $\pi_{a_{5}; j}$ be the projections from 
$\oplus_{a_{5}\in \A}N_{a_{2}a_{3}}^{a_{5}}W^{a_{5}}$ to 
the $j$-th copy of $W^{a_{5}}$. 
Let $f: W^{a_{2}}\boxtimes_{P(z_{2})} W^{a_{3}}
\to \oplus_{a_{5}\in \A}N_{a_{2}a_{3}}^{a_{5}}W^{a_{5}}$ be the 
isomorphism such that 
$$\overline{\pi}_{a_{5}; j}(\overline{f}(w_{a_{2}}\boxtimes_{P(z_{2})}w_{a_{3}}))
=\Y_{a_{2}a_{3}; j}^{a_{5}; (2)}(w_{a_{2}}, z_{2})w_{a_{3}}$$
for $w_{a_{2}}\in W^{a_{2}}$ and $w_{a_{3}}\in W^{a_{3}}$, where
$$\overline{\pi}_{a_{5}; j}: 
\overline{\oplus_{a_{5}\in \A}N_{a_{2}a_{3}}^{a_{5}}W^{a_{5}}}
\to \overline{W}^{a_{5}}$$ 
and 
$$\overline{f}: \overline{W^{a_{2}}\boxtimes_{P(z_{2})} W^{a_{3}}}
\to \overline{\oplus_{a_{5}\in \A}N_{a_{2}a_{3}}^{a_{5}}W^{a_{5}}}$$ 
are the natural 
extensions of $\pi_{a_{5}; j}$ and $f$.
By the universal property for the tensor product module, 
such $f$ indeed exists. Let $\Y_{2}$ be the 
intertwining operator corresponding to the 
intertwining map $\boxtimes_{P(z_{2})}: W^{a_{2}}\otimes W^{a_{3}}
\to W^{a_{2}}\boxtimes_{P(z_{2})} W^{a_{3}}$ (see \cite{HL1}
and \cite{HL4}).
Then we have
\begin{equation}\label{y-2}
\pi_{a_{5}; j}(f(\Y_{2}(w_{a_{2}}, x)w_{a_{3}}))
=\Y_{a_{2}a_{3}; j}^{a_{5}; (2)}(w_{a_{2}}, x)w_{a_{3}}
\end{equation}
for $w_{a_{2}}\in W^{a_{2}}$ and $w_{a_{3}}\in W^{a_{3}}$.
Let $\Y_{1}$ be the intertwining operator 
of type ${W^{a_{4}}\choose W^{a_{1}}\; (W^{a_{2}}\boxtimes_{P(z_{2})} W^{a_{3}})}$
given by 
\begin{equation}\label{y-1}
\Y_{1}(w_{a_{1}}, x)w=\sum_{a_{5}\in \A}\sum_{i=1}^{N_{a_{1}a_{5}}^{a_{4}}}
\sum_{j=1}^{N_{a_{2}a_{3}}^{a_{5}}}\lambda_{a_{5}, i, j}
\Y_{a_{1}a_{5}; i}^{a_{4}; (1)}(w_{a_{1}}, x)\pi_{a_{5}; j}(f(w))
\end{equation}
for $w_{a_{1}}\in W^{a_{1}}$ and 
$w\in W^{a_{2}}\boxtimes_{P(z_{2})} W^{a_{3}}$.

By (\ref{y-2}) and (\ref{y-1}), 
the left-hand side of (\ref{independence1}) is equal to 
\begin{eqnarray*}
&{\displaystyle \sum_{a_{5}\in \A}\sum_{i=1}^{N_{a_{1}a_{5}}^{a_{4}}}
\sum_{j=1}^{N_{a_{2}a_{3}}^{a_{5}}}\lambda_{a_{5}, i, j}
 \langle w_{a'_{4}}, \Y_{a_{1}a_{5}; i}^{a_{4}; (1)}(w_{a_{1}}, x_{1})
\pi_{a_{5}; j}(f(\Y_{2}(w_{a_{2}}, x_{2})w_{a_{3}}))
\rangle}&\nn
&\quad\quad\quad\quad\quad\quad =\langle w_{a'_{4}}, \Y_{1}(w_{a_{1}}, x_{1})
\Y_{2}(w_{a_{2}}, x_{2})w_{a_{3}}
\rangle.&
\end{eqnarray*}
Thus we have 
\begin{equation}\label{independence2}
\langle w_{a'_{4}}, \Y_{1}(w_{a_{1}}, x_{1})
\Y_{2}(w_{a_{2}}, x_{2})w_{a_{3}}
\rangle=0
\end{equation}
for $w_{a'_{4}}\in W^{a'_{4}}$, 
$w_{a_{1}}\in W^{a_{1}}$, $w_{a_{2}}\in W^{a_{2}}$, $w_{a_{3}}\in W^{a_{3}}$.
Since the homogeneous components of $w_{a_{2}}\boxtimes_{P(z_{2})}w_{a_{3}}$ 
or equivalently the $x_{2}$-coefficients of $\Y_{2}(w_{a_{2}}, x_{2})w_{a_{3}}$
for $w_{a_{2}}\in W^{a_{2}}$ and $w_{a_{3}}\in W^{a_{3}}$ span 
$W^{a_{2}}\boxtimes_{P(z_{2})} W^{a_{3}}$ (see \cite{H1}), 
we obtain from (\ref{independence2}) that 
$$\langle w_{a'_{4}}, \Y(w_{a_{1}}, x_{1})w
\rangle=0.$$
for $w_{a'_{4}}\in W^{a'_{4}}$, 
$w_{a_{1}}\in W^{a_{1}}$ and $w\in W^{a_{2}}\boxtimes_{P(z_{2})} W^{a_{3}}$.

Now take $w$ to be an element of  $W^{a_{2}}\boxtimes_{P(z_{2})} W^{a_{3}}$
such that $f(w)$ is in the $j$-th copy of $W^{a_{5}}$ in 
$\oplus_{a_{5}\in \A}N_{a_{2}a_{3}}^{a_{5}}W^{a_{5}}$, that is, take $w$
such that 
$\pi_{a_{5}; j}(f(w))=f(w)$, $\pi_{a_{5}; m}(f(w))=0$ for $m\ne j$
and $\pi_{a; m}(f(w))=0$ for $a \ne a_{5}$. Then we have 
$$\sum_{i=1}^{N_{a_{1}a_{5}}^{a_{4}}}
\lambda_{a_{5}, i, j}
\langle w_{a'_{4}}, \Y_{a_{1}a_{5}; i}^{a_{4}; (1)}(w_{a_{1}}, x_{1})f(w)
\rangle=0.$$
Since $w_{a'_{4}}$, $w_{a_{1}}$, $w$ are arbitrary elements of 
$W^{a'_{4}}$, $W^{a_{1}}$ and the $j$-th copy of $W^{a_{5}}$ in 
$\oplus_{a_{5}\in \A}N_{a_{2}a_{3}}^{a_{5}}W^{a_{5}}$, respectively, 
we obtain
$$\sum_{i=1}^{N_{a_{1}a_{5}}^{a_{4}}}
\lambda_{a_{5}, i, j}\Y_{a_{1}a_{5}; i}^{a_{4}; (1)}=0.$$
Since $\Y_{a_{1}a_{5}; i}^{a_{4}; (1)}$ for 
$i=1, \dots, N_{a_{1}a_{5}}^{a_{4}}$
are linearly 
independent, we obtain
$\lambda_{a_{5}, i, j}=0$ for  
$a_{5}\in \A$, $i=1, \dots, N_{a_{1}a_{5}}^{a_{4}}$ and 
$j=1, \dots, N_{a_{2}a_{3}}^{a_{5}}$. 
\epfv

Since
\begin{equation}\label{3-product}
\langle w_{a'_{2}}, \Y_{a_{2}e; 1}^{a_{2}}(w^{1}_{a_{2}}, z_{1})
\Y_{a'_{3}a_{3}; 1}^{e}(w_{a'_{3}}, z_{2})
\Y_{a_{1}a_{2}; i}^{a_{3}}(w_{a_{1}}, z_{3})w^{2}_{a_{2}}\rangle
\end{equation}
and 
\begin{equation}\label{3-iterate}
\langle w_{a'_{2}}, 
\Y_{a_{4}a_{3}; k}^{a_{2}}(\Y_{a_{2}a'_{3}; j}^{a_{4}}(w^{1}_{a_{2}}, 
z_{1}-z_{2})w_{a'_{3}}, z_{2})
\Y_{a_{1}a_{2}; i}^{a_{3}}(w_{a_{1}}, z_{3})w^{2}_{a_{2}}\rangle
\end{equation}
satisfy a system of differential equations of regular singular
points with coefficients in 
$$\C[z_{1}, z_{2}^{-1}, z_{2}, z_{2}^{-1}, z_{3}, z_{3}^{-1}, 
(z_{1}-z_{2})^{-1}, (z_{1}-z_{3})^{-1}, (z_{2}-z_{3})^{-1}]$$
(see \cite{H8}),
they are absolutely convergent in the region given by
$|z_{1}|>|z_{2}|>|z_{3}|>0$ and in the region given by
$|z_{2}|>|z_{1}-z_{2}|>0$, $|z_{2}|>|z_{3}|>0$, $|z_{2}-z_{3}|>|z_{1}-z_{2}|>0$,
respectively.  Using also these differential equations, we see that 
(\ref{3-product})
and (\ref{3-iterate})
can also be analytically extended to multi-valued analytic functions
on 
$$M^{3}=\{(z_{1}, z_{2}, z_{3})\in \C^{3}\;|\; 
z_{1}, z_{2}, z_{3}\ne 0, z_{1}\ne z_{2}, z_{2}\ne z_{3}, z_{1}\ne z_{3}\}.$$
We can also choose single valued branches of these multi-valued analytic functions
on the region given by 
$|z_{1}|>|z_{2}|>|z_{3}|>0$ with cuts along the real lines in the 
$z_{1}$-, $z_{2}$- and $z_{3}$-planes and on the region 
$|z_{2}|>|z_{1}-z_{2}|>0$, $|z_{2}|>|z_{3}|>0$, $|z_{2}-z_{3}|>|z_{1}-z_{2}|>0$
with cuts along the real lines in the 
$z_{2}$-, $z_{1}-z_{2}$- and $z_{3}$-planes. 

Similarly to the case of two intertwining 
operators above, we also have the single-valued analytic functions
$$E(\langle w_{a'_{2}}, \Y_{a_{2}e; 1}^{a_{2}}(w^{1}_{a_{2}}, z_{1})
\Y_{a'_{3}a_{3}; 1}^{e}(w_{a'_{3}}, z_{2})
\Y_{a_{1}a_{2}; i}^{a_{3}}(w_{a_{1}}, z_{3})w^{2}_{a_{2}}\rangle)$$
and 
$$E(\langle w_{a'_{2}}, 
\Y_{a_{4}a_{3}; k}^{a_{2}}(\Y_{a_{2}a'_{3}; j}^{a_{4}}(w^{1}_{a_{2}}, 
z_{1}-z_{2})w_{a'_{3}}, z_{2})
\Y_{a_{1}a_{2}; i}^{a_{3}}(w_{a_{1}}, z_{3})w^{2}_{a_{2}}\rangle)$$
on the universal covering $\widetilde{M}^{3}$ of $M^{3}$. 
The associativity can also be written using these single valued analytic
functions.
Similarly we have other single valued analytic functions
on  $\widetilde{M}^{3}$, for example, 
\begin{eqnarray*}
&E(\langle w_{a'_{2}}, 
\Y_{a_{4}a_{3}; k}^{a_{2}}(\Y_{a_{2}a'_{3}; j}^{a_{4}}(w^{1}_{a_{2}}, 
z_{1}-z_{2})w_{a'_{3}}, z_{2})
\Y_{a_{1}a_{2}; i}^{a_{3}}(w_{a_{1}}, z_{3})w^{2}_{a_{2}}\rangle),&\\
&B^{(r)}_{12}(E(\langle w_{a'_{2}}, \Y_{a_{2}e; 1}^{a_{2}}(w^{1}_{a_{2}}, z_{1})
\Y_{a'_{3}a_{3}; 1}^{e}(w_{a'_{3}}, z_{2})
\Y_{a_{1}a_{2}; i}^{a_{3}}(w_{a_{1}}, z_{3})w^{2}_{a_{2}}\rangle)),&\\
&B^{(r)}_{23}(E(\langle w_{a'_{2}}, \Y_{a_{2}e; 1}^{a_{2}}(w^{1}_{a_{2}}, z_{1})
\Y_{a'_{3}a_{3}; 1}^{e}(w_{a'_{3}}, z_{2})
\Y_{a_{1}a_{2}; i}^{a_{3}}(w_{a_{1}}, z_{3})w^{2}_{a_{2}}\rangle))&
\end{eqnarray*}
and so on, where the subscripts $12$ and $23$ in 
$B^{(r)}_{12}$ and $B^{(r)}_{23}$, respectively, mean that they 
corresponding to braiding isomorphisms which braid the first two and 
the last two intertwining operators, respectively. 

Similarly to the case of two intertwining operators, we also 
have the following result: 

\begin{prop}\label{independence-3p}
For $a_{1}, a_{2}, a_{3}, a_{4}, a_{5}\in \A$, the maps from 
$W^{a'_{5}}\otimes W^{a_{1}}\otimes W^{a_{2}}\otimes W^{a_{3}}\otimes W^{a_{4}}$
to the space of single-valued analytic functions on $\widetilde{M}^{3}$
given by 
\begin{eqnarray*}
\lefteqn{w_{a'_{5}}\otimes w_{a_{1}} \otimes w_{a_{2}}\otimes w_{a_{3}}\otimes w_{a_{4}}
\mapsto}\nn
&&\quad\quad E(\langle w_{a'_{4}}, \Y_{a_{6}a_{4}; k}^{a_{5}; (1)}(
\Y_{a_{7}a_{3}; l}^{a_{6}; (2)}(
\Y_{a_{1}a_{2}; m}^{a_{7}; (3)}(w_{a_{1}}, z_{1}-z_{2})w_{a_{2}}, 
z_{2}-z_{3})w_{a_{3}}, z_{3})w_{a_{4}}\rangle)
\end{eqnarray*}
for $a_{6}, a_{7}\in \A$, $k=1, \dots, N_{a_{6}a_{4}}^{a_{5}}$,
$l=1, \dots, N_{a_{7}a_{3}}^{a_{6}}$ and $m=1, \dots, N_{a_{1}a_{2}}^{a_{7}}$
are linearly independent. 
Similarly, $a_{1}, a_{2}, a_{3}, a_{4}, a_{5}\in \A$, the maps from 
$W^{a'_{5}}\otimes W^{a_{1}}\otimes W^{a_{2}}\otimes W^{a_{3}}\otimes W^{a_{4}}$
to the space of single-valued analytic functions on $\widetilde{M}^{3}$
given by the lifting to the universal covering of 
analytic extensions of products or combinations of 
products and iterates of bases of intertwining operators are linearly 
independent. 
\end{prop}

The proof of this proposition is similar to the proof of Proposition 
\ref{independence0} and is omitted.

\renewcommand{\theequation}{\thesection.\arabic{equation}}
\renewcommand{\thethm}{\thesection.\arabic{thm}}
\setcounter{equation}{0}
\setcounter{thm}{0}

\section{Geometrically-modified intertwining operators and 
genus-one correlation functions}

We now discuss genus-one correlation functions. In the present paper,
we need only one or two point functions. Genus-one correlation functions 
are defined using what we called geometrically-modified intertwining 
operators (see \cite{H9}). We shall first discuss these operators and prove
some important properties needed in this paper. 

Let $A_{j}$, $j\in \mathbb{Z}_{+}$, be the 
complex numbers defined 
by 
\begin{eqnarray*}
\frac{1}{2\pi i}\log(1+2\pi i y)=\left(\exp\left(\sum_{j\in \mathbb{Z}_{+}}
A_{j}y^{j+1}\frac{\partial}{\partial y}\right)\right)y.
\end{eqnarray*}
For any $V$-module 
$W$, let
\begin{equation}\label{ux}
\mathcal{U}(x)=(2\pi ix)^{L(0)}e^{-L^{+}(A)}\in (\mbox{\rm End}\;W)\{x\}
\end{equation}
where $L_{+}(A)=\sum_{j\in \mathbb{Z}_{+}}
A_{j}L(j)$ and $L(j)$ for $j\in \mathbb{Z}$ are the Virasoro operators
on $W$. Given an intertwining operator $\Y$ of type ${W^{a_{3}}\choose
W^{a_{1}}W^{a_{2}}}$, the map $W^{a_{1}}\otimes W^{a_{2}}\to
W^{a_{3}}\{x\}$ defined  by $w_{a_{1}}\otimes w_{a_{2}}\mapsto 
\Y(\mathcal{U}(x)w_{a_{1}}, x)w_{a_{2}}$ is the corresponding 
geometrically-modified intertwining operator. We need the following 
property of geometrically-modified intertwining operators:

\begin{prop}
Let $\Y$ be an intertwining operator of type ${W^{a_{3}}\choose
W^{a_{1}}W^{a_{2}}}$. Then for $w_{a_{1}}\in W^{a_{1}}$, $w_{a_{1}}\in W^{a_{1}}$
and $w_{a_{3}}'\in W^{a'_{3}}$, we have 
\begin{eqnarray}\label{contrag-m}
\langle w_{a_{3}}', \sigma_{23}(\Y)(\mathcal{U}(x)w_{a_{1}}, x)w_{a_{2}}\rangle
&=&e^{\pi i h_{a_{1}}}\langle \Y(\mathcal{U}(x^{-1})e^{-\pi iL(0)}w_{a_{1}}, 
x^{-1})w_{a_{3}}', w_{a_{2}}\rangle\nn
&=&e^{-\pi i h_{a_{1}}}\langle \Y(\mathcal{U}(x^{-1})e^{\pi iL(0)}w_{a_{1}}, 
x^{-1})w_{a_{3}}', w_{a_{2}}\rangle.\nn
&&
\end{eqnarray}
\end{prop}
\pf
We first prove the first equality in (\ref{contrag-m}).
By the definition of $\sigma_{23}(\Y)$, we have
\begin{eqnarray}\label{contrag}
\lefteqn{\langle w_{a_{3}}', \sigma_{23}(\Y)(\mathcal{U}(x)w_{a_{1}}, 
x)w_{a_{2}}\rangle}\nn
&&=e^{\pi i h_{a_{1}}}\langle \Y(e^{xL(1)}(e^{-\pi i}x^{-2})^{L(0)}
\mathcal{U}(x)w_{a_{1}}, x^{-1})w_{a_{3}}', w_{a_{2}}\rangle.
\end{eqnarray}
We now calculate $e^{xL(1)}(e^{-\pi i}x^{-2})^{L(0)}
\mathcal{U}(x)w$. The equality
$$\frac{1}{2\pi i}\log \left(1-\frac{x}{y^{-1}+x}\right)
=-\frac{1}{2\pi i}\log (1+xy)$$
gives 
\begin{eqnarray*}
\lefteqn{e^{\displaystyle -xy^{2}\frac{\partial}{\partial y}}
\; (e^{-\pi i}x^{-2})^{\displaystyle -y\frac{\partial}{\partial y}}\; 
(2\pi i x)^{\displaystyle -y\frac{\partial}{\partial y}}\; e^{\displaystyle 
\sum_{j\in \Z_{+}}A_{j}y^{j+1}
\frac{\partial}{\partial y}} \; y}\nn
&&\quad\quad
=(2\pi i x^{-1})^{\displaystyle -y\frac{\partial}{\partial y}}\; 
e^{\displaystyle \sum_{j\in \Z_{+}}A_{j}y^{j+1}
\frac{\partial}{\partial y}} \; e^{\displaystyle \pi i y\frac{\partial}{\partial y}}\;y.
\end{eqnarray*}
Using the method developed and results obtained 
in Chapters 4 and 5 of \cite{H4}, we obtain
$$e^{xL(1)}(e^{-\pi i}x^{-2})^{L(0)}
(2\pi i x)^{L(0)}e^{-L_{+}(A)}
=(2\pi i x^{-1})^{L(0)}e^{-L_{+}(A)} e^{-\pi i L(0)},$$
that is,
$$e^{xL(1)}(e^{-\pi i}x^{-2})^{L(0)}\mathcal{U}(x)
= \mathcal{U}(x^{-1})e^{-\pi i L(0)},$$
Thus the right-hand side of (\ref{contrag}) is equal to
$$e^{\pi i h_{a_{1}}} \langle \Y(\mathcal{U}(x^{-1})e^{-\pi i L(0)}w_{a_{1}}, 
x^{-1})w_{a_{3}}', w_{a_{2}}\rangle,$$
proving the first equality in (\ref{contrag-m}).

The second equality in (\ref{contrag-m}) follows immediately
from the first equality and the equality
\begin{eqnarray*}
e^{-\pi L(0)}w_{a_{1}}&=&e^{\pi L(0)}e^{-2\pi L(0)}w_{a_{1}}\nn
&=&e^{-2\pi h_{a_{1}}}e^{\pi L(0)}w_{a_{1}}. 
\end{eqnarray*}
\epfv

As in the preceding section, For $a_{1}, a_{2}, 
a_{3}\in \mathcal{A}$ and $p=1, 2, 3, 4, 5, 6$, let 
$\Y_{a_{1}a_{2}; i}^{a_{3}; (p)}$, $i=1, \dots, N_{a_{1}a_{2}}^{a_{3}}$,
be bases of $\mathcal{V}_{a_{1}a_{2}}^{a_{3}}$. Let $q_{\tau}=e^{2\pi i\tau}$
for $\tau \in \mathbb{H}$. 
We consider $q_{\tau}$-traces of geometrically-modified 
intertwining operators of the following form:
\begin{equation}\label{1-trace}
\tr_{W^{a_{2}}}
\Y_{a_{1}a_{2}; i}^{a_{2}; (1)}(\mathcal{U}(e^{2\pi iz})w_{a_{1}}, e^{2\pi i z})
q_{\tau}^{L(0)-\frac{c}{24}}
\end{equation}
for $a_{1}, a_{2}\in \A$, $i=1, \dots, N_{a_{1}a_{2}}^{a_{2}}$. 
In \cite{M} and \cite{H9}, it was shown that these $q$-traces are independent of 
$z$, are absolutely convergent
when $0<|q_{\tau}|<1$ and can be analytically extended to 
analytic functions of $\tau$ in the upper-half plane. 
We shall denote the analytic extension of (\ref{1-trace})
by 
$$E(\tr_{W^{a_{2}}}
\Y_{a_{1}a_{2}; i}^{a_{2}; (1)}(\mathcal{U}(e^{2\pi iz})w_{a_{1}}, e^{2\pi i z})
q_{\tau}^{L(0)-\frac{c}{24}}).$$
These are genus-one one-point correlation functions. 
In \cite{M} and \cite{H9},  
the following modular invariance property is also proved:
For 
$$\left(\begin{array}{cc}a&b\\ c&d\end{array}\right)\in SL(2, \Z),$$
let $\tau'=\frac{a\tau+b}{c\tau+d}$. Then for fixed $a_{1}\in 
\mathcal{A}$, there exist unique
$A_{a_{2}}^{a_{3}}\in \C$ for $a_{2}, a_{3}\in \mathcal{A}$ 
such that 
\begin{eqnarray*}
\lefteqn{\tr_{W^{a_{2}}}
\Y_{a_{1}a_{2}; i}^{a_{2}; (1)}
\left(\mathcal{U}(e^{2\pi i\frac{z}{c\tau+d}})
\left(\frac{1}{c\tau+d}\right)^{L(0)}
w_{a_{1}}, e^{2\pi i \frac{z}{c\tau+d}}\right)
q_{\tau'}^{L(0)-\frac{c}{24}}}\nn
&&=\sum_{j=1}^{N_{a_{1}a_{3}}^{a_{3}}}
\sum_{a_{3}\in \A}A_{a_{2}; i}^{a_{3}; j}
\tr_{W^{a_{3}}}
\Y_{a_{1}a_{3}; j}^{a_{3}; (1)}(\mathcal{U}(e^{2\pi iz})w_{a_{1}}, e^{2\pi i z})
q_{\tau}^{L(0)-\frac{c}{24}}
\end{eqnarray*}
for $w_{a_{1}}\in W^{a_{1}}$ and $i=1, \dots, N_{a_{1}a_{2}}^{a_{2}}$. 
This modular invariance property can be written in terms of 
the analytic extensions as follows: 
\begin{eqnarray*}
\lefteqn{E\left(\tr_{W^{a_{2}}}
\Y_{a_{1}a_{2}; i}^{a_{2}; (1)}\left(\mathcal{U}(e^{\frac{2\pi iz}{c\tau+d}})
\left(\frac{1}{c\tau+d}\right)^{L(0)}
w_{a_{1}}, e^{\frac{2\pi i z}{c\tau+d}}\right)
q_{\tau'}^{L(0)-\frac{c}{24}}\right)}\nn
&&=\sum_{j=1}^{N_{a_{1}a_{3}}^{a_{3}}}\sum_{a_{3}\in \A}A_{a_{2}}^{a_{3}}
E(\tr_{W^{a_{3}}}
\Y_{a_{1}a_{3}; i}^{a_{3}; (1)}(\mathcal{U}(e^{2\pi iz})w_{a_{1}}, e^{2\pi i z})
q_{\tau}^{L(0)-\frac{c}{24}}).
\end{eqnarray*}

We consider
the space $\F^{e}_{1; 1}$ spanned by linear maps of the form
\begin{eqnarray*}
\Psi_{a}: V&\to &
\mathbb{G}^{e}_{1; 1}\nn
u&\mapsto & 
\Psi_{a}(u; \tau)
\end{eqnarray*}
for $a\in \A$, where 
$$\Psi_{a}(u; \tau)=E(\tr_{W^{a}}
\Y_{ea; 1}^{a; (1)}(\mathcal{U}(e^{2\pi iz})u, e^{2\pi i z})
q_{\tau}^{L(0)-\frac{c}{24}})$$
and $\mathbb{G}^{e}_{1; 1}$ is  the space spanned by
functions 
of $\tau$ of the form 
$\Psi_{a}(u; \tau)$.

We also need genus-one two-point correlation functions. 
They are constructed from the $q_{\tau}$-traces of products or iterates
of two geometrically-modified intertwining operators as 
follows: For $a_{1}, a_{2}, a_{3}, a_{4}\in \A$, 
$i=1, \dots$, $N_{a_{1}a_{3}}^{a_{4}}$, $j=1, \dots, 
N_{a_{2}a_{4}}^{a_{3}}$, consider 
\begin{equation}\label{1-product}
\tr_{W^{a_{4}}}
\Y_{a_{1}a_{3}; i}^{a_{4}; (1)}(\mathcal{U}(e^{2\pi iz_{1}})w_{a_{1}}, 
e^{2\pi i z_{1}})\Y_{a_{2}a_{4}; j}^{a_{3}; (2)}
(\mathcal{U}(e^{2\pi iz_{2}})w_{a_{2}},  e^{2\pi i z_{2}})
q_{\tau}^{L(0)-\frac{c}{24}}.
\end{equation}
In \cite{H9}, it was shown that these $q_{\tau}$-traces
are absolutely convergent when
$0<|q_{\tau}|<|e^{2\pi i z_{1}}|<|e^{2\pi i z_{2}}|<1$
and can be analytically extended to multi-valued 
analytic functions on 
$$M_{1}^{2}=\{(z_{1}, z_{2}, \tau)\in \C^{3}\;|\; 
z_{1}\ne z_{2}+p\tau +q \; \mbox{\rm for}\; 
p, q\in \Z,\; \tau\in \mathbb{H}\}.$$
(In fact, these multi-valued analytic functions depend only on $z_{1}-z_{2}$
and $\tau$ by the $L(-1)$-derivative property. See \cite{H9}.)
These multi-valued analytic functions on $M_{1}^{2}$
are genus-one two-point correlation functions. 
They can be lifted to 
single-valued analytic functions on the universal 
covering $\widetilde{M}_{1}^{2}$. As in the case of 
genus-zero correlation functions, these liftings are 
not unique. To obtain a unique single-valued 
analytic function on $\widetilde{M}_{1}^{2}$ from 
a multi-valued analytic functions $f$ on $M_{1}^{2}$,
we have to give a single valued branch of $f$ on some simply connected
region in $M_{1}^{2}$ or values of $f$ 
on some simply-connected subset of $M_{1}^{2}$. 
In fact, (\ref{1-product}) gives a single valued 
branch of its analytic extension in the region given by
$0<|q_{\tau}|<|e^{2\pi i z_{1}}|<|e^{2\pi i z_{2}}|<1$. 
We denote this single-valued analytic function on 
$\widetilde{M}_{1}^{2}$ by 
\begin{equation}\label{1-product-e}
E(\tr_{W^{a_{4}}}
\Y_{a_{1}a_{3}; i}^{a_{4}; (1)}(\mathcal{U}(e^{2\pi iz_{1}})w_{a_{1}}, 
e^{2\pi i z_{1}})\Y_{a_{2}a_{4}; j}^{a_{3}; (2)}
(\mathcal{U}(e^{2\pi iz_{2}})w_{a_{2}},  e^{2\pi i z_{2}})
q_{\tau}^{L(0)-\frac{c}{24}}).
\end{equation}

Similarly, these genus-one two-point correlation functions can
also be constructed from $q_{\tau}$-traces of iterates of 
two intertwining operators as follows:
For $a_{1}, a_{2}, a_{5}, a_{6}\in \A$, 
$k=1, \dots, N_{a_{3}a_{4}}^{a_{4}}$, $l=1, \dots, 
N_{a_{1}a_{2}}^{a_{3}}$, consider 
\begin{equation}\label{1-iterate}
\tr_{W^{a_{6}}}
\Y_{a_{5}a_{6}; k}^{a_{6}; (3)}(\mathcal{U}(e^{2\pi iz_{1}})
\Y_{a_{1}a_{2}; l}^{a_{5}; (4)}
(w_{a_{1}},  z_{1}-z_{2})w_{a_{2}}, 
e^{2\pi i z_{2}})
q_{\tau}^{L(0)-\frac{c}{24}}.
\end{equation}
In \cite{H9}, it was shown that these $q_{\tau}$-traces
are absolutely convergent when
$0<|q_{\tau}|<|e^{2\pi i z_{2}}|<1$ and
$0<|e^{2\pi i(z_{1}-z_{2})}-1|<1$
and can be analytically extended to multi-valued 
analytic functions on $M_{1}^{2}$.
(Again, these multi-valued analytic functions depend only on $z_{1}-z_{2}$
and $\tau$ by the $L(-1)$-derivative property.)
These multi-valued analytic functions on $M_{1}^{2}$
are also genus-one two-point correlation functions. 
The multi-valued analytic extension of (\ref{1-iterate})
can also be lifted
uniquely to a 
single-valued analytic function on $\widetilde{M}_{1}^{2}$
using the single-valued branch (\ref{1-iterate}). 
We denote it by 
\begin{equation}\label{1-iterate-e}
E(\tr_{W^{a_{6}}}
\Y_{a_{5}a_{6}; k}^{a_{6}; (3)}(\mathcal{U}(e^{2\pi iz_{1}})
\Y_{a_{1}a_{2}; l}^{a_{5}; (4)}
(w_{a_{1}},  z_{1}-z_{2})w_{a_{2}}, 
e^{2\pi i z_{2}})
q_{\tau}^{L(0)-\frac{c}{24}}).
\end{equation}

In \cite{H9}, an associativity property for geometrically-modified
intertwining operators is
proved. This associativity together with the convergence 
property of $q_{\tau}$-traces of gives
\begin{eqnarray}\label{1-assoc}
\lefteqn{E(\tr_{W^{a_{4}}}
\Y_{a_{1}a_{3}; i}^{a_{4}; (1)}(\mathcal{U}(e^{2\pi iz_{1}})w_{a_{1}}, 
e^{2\pi i z_{1}})\Y_{a_{2}a_{4}; j}^{a_{3}; (2)}
(\mathcal{U}(e^{2\pi iz_{2}})w_{a_{2}},  e^{2\pi i z_{2}})
q_{\tau}^{L(0)-\frac{c}{24}})}\nn
&&=\sum_{a_{5}\in \A}\sum_{l=1}^{N_{a_{1}a_{2}}^{a_{5}}}
\sum_{k=1}^{N_{a_{5}a_{4}}^{a_{4}}}
F(\Y_{a_{1}a_{3}; i}^{a_{4}; (1)}\otimes \Y_{a_{2}a_{4}; j}^{a_{3}; (2)}; 
\Y_{a_{5}a_{4}; k}^{a_{4}; (3)}
\otimes \Y_{a_{1}a_{2}; l}^{a_{5}; (4)})\cdot\nn
&&\quad\quad\quad\cdot 
E(\tr_{W^{a_{4}}}
\Y_{a_{5}a_{4}; k}^{a_{4}; (3)}(\mathcal{U}(e^{2\pi iz_{1}})
\Y_{a_{1}a_{2}; l}^{a_{5}; (4)}
(w_{a_{1}},  z_{1}-z_{2})w_{a_{2}}, 
e^{2\pi i z_{2}})
q_{\tau}^{L(0)-\frac{c}{24}})
\end{eqnarray}
and
\begin{eqnarray}\label{1-assoc-inv}
\lefteqn{E(\tr_{W^{a_{4}}}
\Y_{a_{3}a_{4}; k}^{a_{4}; (3)}(\mathcal{U}(e^{2\pi iz_{1}})
\Y_{a_{1}a_{2}; l}^{a_{3}; (4)}
(w_{a_{1}},  z_{1}-z_{2})w_{a_{2}}, 
e^{2\pi i z_{2}})
q_{\tau}^{L(0)-\frac{c}{24}}}\nn
&&=\sum_{a_{5}\in \A}\sum_{i=1}^{N_{a_{1}a_{2}}^{a_{6}}}
\sum_{j=1}^{N_{a_{6}a_{3}}^{a_{4}}}
F^{-1}(\Y_{a_{3}a_{4}; k}^{a_{4}; (3)}
\otimes \Y_{a_{1}a_{2}; l}^{a_{3}; (4)}; 
\Y_{a_{1}a_{5}; i}^{a_{4}; (1)}\otimes \Y_{a_{2}a_{4}; j}^{a_{5}; (2)})\cdot\nn
&&\quad\quad\quad\cdot 
E(\tr_{W^{a_{4}}}
\Y_{a_{1}a_{5}; i}^{a_{4}; (1)}(\mathcal{U}(e^{2\pi iz_{1}})w_{a_{1}}, 
e^{2\pi i z_{1}})\cdot\nn
&&\quad\quad\quad\quad\quad\quad\quad\quad\quad\quad\quad\quad
\cdot \Y_{a_{2}a_{4}; j}^{a_{5}; (2)}
(\mathcal{U}(e^{2\pi iz_{2}})w_{a_{2}},  e^{2\pi i z_{2}})
q_{\tau}^{L(0)-\frac{c}{24}}).
\end{eqnarray}
In particular, the space  of all single-valued analytic functions 
on  $\widetilde{M}_{1}^{2}$ spanned by functions of the form
(\ref{1-product-e}) and the space of all single-valued analytic functions 
on  $\widetilde{M}_{1}^{2}$ spanned by functions of the form
(\ref{1-iterate-e}) are the same. We shall denote this space by
$\mathbb{G}^{e}_{1; 2}$. 

We need the following result: 

\begin{prop}\label{1-independence0}
For $a_{1}, a_{2}\in \A$, the maps from $W^{a_{1}}\otimes W^{a_{2}}$
to $\mathbb{G}^{e}_{1; 2}$ given by 
\begin{eqnarray*}
\lefteqn{w_{a_{1}}\otimes w_{a_{2}}\mapsto}\nn 
&& E(\tr_{W^{a_{4}}}
\Y_{a_{1}a_{3}; i}^{a_{4}; (1)}(\mathcal{U}(e^{2\pi iz_{1}})w_{a_{1}}, 
e^{2\pi i z_{1}})\Y_{a_{2}a_{4}; j}^{a_{3}; (2)}
(\mathcal{U}(e^{2\pi iz_{2}})w_{a_{2}},  e^{2\pi i z_{2}})
q_{\tau}^{L(0)-\frac{c}{24}}),
\end{eqnarray*}
$a_{3}, a_{4}\in \A$, $i=1, \dots, N_{a_{1}a_{3}}^{a_{4}}$,
$j=1, \dots, N_{a_{2}a_{4}}^{a_{3}}$,
are linearly independent. Similarly, for $a_{1}, a_{2}\in \A$, 
the maps from $W^{a_{1}}\otimes W^{a_{2}}$
to $\mathbb{G}^{e}_{1; 2}$ given by 
$$w_{a_{1}}\otimes w_{a_{2}}\mapsto E(\tr_{W^{a_{4}}}
\Y_{a_{3}a_{4}; k}^{a_{4}; (3)}(\mathcal{U}(e^{2\pi iz_{2}})
\Y_{a_{1}a_{2}; l}^{a_{3}; (4)}
(w_{a_{1}},  z_{1}-z_{2})w_{a_{2}}, 
e^{2\pi i z_{2}})
q_{\tau}^{L(0)-\frac{c}{24}}),$$
$a_{3}, a_{4}\in \A$, $k=1, \dots, N_{a_{3}a_{4}}^{a_{3}}$,
$l=1, \dots, N_{a_{1}a_{2}}^{a_{3}}$,
are linearly independent.
\end{prop}
\pf
We prove only the linear independence of 
the maps obtained
from $q_{\tau}$-traces of iterates of intertwining operators. 
For the linear independence
of the maps obtained from traces of products of intertwining 
operators, the proof is similar.

Since analytic extensions are unique, we need only prove that 
the maps given by 
$$w_{a_{1}} \otimes w_{a_{2}}
\mapsto \tr_{W^{a_{4}}}
\Y_{a_{3}a_{4}; k}^{a_{4}; (3)}(\mathcal{U}(e^{2\pi iz_{1}})
\Y_{a_{1}a_{2}; l}^{a_{3}; (4)}
(w_{a_{1}},  z_{1}-z_{2})w_{a_{2}}, 
e^{2\pi i z_{2}})
q_{\tau}^{L(0)-\frac{c}{24}},$$
$a_{3}, a_{4}\in \A$, $k=1, \dots, N_{a_{3}a_{4}}^{a_{3}}$,
$l=1, \dots, N_{a_{1}a_{2}}^{a_{3}}$,
are linearly independent. Assume 
\begin{eqnarray}\label{1-independence}
\lefteqn{\sum_{a_{3}, a_{4}\in \A}\sum_{k=1}^{N_{a_{3}a_{4}}^{a_{4}}}
\sum_{l=1}^{N_{a_{1}a_{4}}^{a_{2}}}\lambda_{a_{3}, a_{4}, k, l}\;\cdot}\nn
&&\quad\quad\quad\quad \cdot \tr_{W^{a_{4}}}
\Y_{a_{3}a_{4}; k}^{a_{4}; (3)}(\mathcal{U}(e^{2\pi iz_{2}})
\Y_{a_{1}a_{2}; l}^{a_{3}; (4)}
(w_{a_{1}},  z_{1}-z_{2})w_{a_{2}}, 
e^{2\pi i z_{2}})
q_{\tau}^{L(0)-\frac{c}{24}}\nn
&&=0
\end{eqnarray}
for $w_{a_{1}}\in W^{a_{1}}$ and $w_{a_{2}}\in W^{a_{2}}$. 
Since (\ref{1-independence}) holds for all $z_{1}$ and $z_{2}$
satisfying
$0<|q_{\tau}|<|e^{2\pi i z_{2}}|<1$ and
$0<|e^{2\pi i(z_{1}-z_{2})}-1|<1$, we obtain the following
equation in which the variable $z_{1}-z_{2}$ is replaced by a 
formal variable $x_{0}$:
\begin{eqnarray}\label{1-independence1}
\lefteqn{\sum_{a_{3}, a_{4}\in \A}\sum_{k=1}^{N_{a_{3}a_{4}}^{a_{4}}}
\sum_{l=1}^{N_{a_{1}a_{4}}^{a_{2}}}\lambda_{a_{3}, a_{4}, k, l}\;\cdot}\nn
&&\quad\quad\quad\quad \cdot \tr_{W^{a_{4}}}
\Y_{a_{3}a_{4}; k}^{a_{4}; (3)}(\mathcal{U}(e^{2\pi iz_{2}})
\Y_{a_{1}a_{2}; l}^{a_{3}; (4)}
(w_{a_{1}},  x_{0})w_{a_{2}}, 
e^{2\pi i z_{2}})
q_{\tau}^{L(0)-\frac{c}{24}}\nn
&&=0
\end{eqnarray}
for $w_{a_{1}}\in W^{a_{1}}$ and $w_{a_{2}}\in W^{a_{2}}$. 
We want to show that $\lambda_{a_{3}, a_{4}, k, l}=0$ for 
$a_{3}, a_{4}\in \A$, $k=1, \dots, N_{a_{3}a_{4}}^{a_{3}}$ 
and $l=1, \dots, N_{a_{1}a_{2}}^{a_{3}}$.

As in the proof of Proposition \ref{independence0}, there 
exists an isomorphism 
$$f: W^{a_{1}}\boxtimes_{P(z_{1}-z_{2})} W^{a_{2}}\to 
\oplus_{a_{3}\in \A}N_{a_{1}a_{2}}^{a_{3}}W^{a_{3}}$$
such that for $w_{a_{1}}\in W^{a_{1}}$ and $w_{a_{2}}\in W^{a_{2}}$,
$$\overline{\pi}_{a_{3}; l}(\overline{f}(w_{a_{1}}
\boxtimes_{P(z_{z_{1}-z_{2}})}w_{a_{2}}))
=\Y_{a_{1}a_{2}; l}^{a_{3}; (4)}(w_{a_{1}}, z_{1}-z_{2})w_{a_{2}},$$
where  $\pi_{a_{3}; l}$  is the projections from 
$\oplus_{a_{3}\in \A}N_{a_{1}a_{2}}^{a_{3}}W^{a_{3}}$ to 
the $l$-th copy of $W^{a_{3}}$ and 
$$\overline{\pi}_{a_{3}; l}: 
\overline{\oplus_{a_{3}\in \A}N_{a_{1}a_{2}}^{a_{3}}W^{a_{3}}}\to 
\overline{W}^{a_{3}}$$
and 
$$\overline{f}: \overline{W^{a_{1}}\boxtimes_{P(z_{1}-z_{2})} W^{a_{2}}}
\to \overline{\oplus_{a_{3}\in \A}N_{a_{1}a_{2}}^{a_{3}}W^{a_{3}}}$$
are the natural extension of $\pi_{a_{3}; l}$ and $f$, respectively. 
Let $\Y_{2}$ be the 
intertwining operator corresponding to the 
intertwining map $\boxtimes_{P(z_{2})}: W^{a_{1}}\otimes W^{a_{2}}
\to W^{a_{1}}\boxtimes_{P(z_{1}-z_{2})} W^{a_{2}}$ (see \cite{HL1}
and \cite{HL4}).
Then we have
\begin{equation}\label{y-2-1}
\pi_{a_{3}; l}(f(\Y_{2}(w_{a_{1}}, x)w_{a_{2}}))
=\Y_{a_{1}a_{2}; j}^{a_{3}; (4)}(w_{a_{1}}, x)w_{a_{2}}
\end{equation}
for $w_{a_{1}}\in W^{a_{1}}$ and $w_{a_{2}}\in W^{a_{2}}$.
For $a_{4}\in \A$, let $\Y_{a_{4}}$ be the 
intertwining operator of type ${W^{a_{4}}\choose 
(W^{a_{1}}\boxtimes_{P(z_{1}-z_{2})} W^{a_{2}})\; 
W^{a_{4}}}$
given by 
$$\Y_{a_{4}}(w, x)w_{a_{4}}
=\sum_{a_{3}\in \A}\sum_{k=1}^{N_{a_{3}a_{4}}^{a_{4}}}
\sum_{l=1}^{N_{a_{1}a_{4}}^{a_{2}}}\lambda_{a_{3}, a_{4}, k, l}
\Y_{a_{3}a_{4}; k}^{a_{4}; (3)}(
\pi_{a_{3}; l}(f(w)), x)w_{a_{4}}$$
for $w\in W^{a_{1}}\boxtimes_{P(z_{1}-z_{2})} W^{a_{2}}$
and $w_{a_{4}}\in W^{a_{4}}$. 
Then the left-hand side of (\ref{1-independence1}) is equal to 
\begin{eqnarray}\label{1-independence1.5}
\lefteqn{\sum_{a_{3}, a_{4}\in \A}\sum_{k=1}^{N_{a_{3}a_{4}}^{a_{4}}}
\sum_{l=1}^{N_{a_{1}a_{2}}^{a_{3}}}\lambda_{a_{3}, a_{4}, k, l}\;\cdot}\nn
&&\quad\quad\quad\cdot \tr_{W^{a_{4}}}
\Y_{a_{3}a_{4}; k}^{a_{4}; (3)}(\mathcal{U}(e^{2\pi iz_{2}})
\pi_{a_{3}; l}(f(\Y_{2}(w_{a_{1}}, x_{0})w_{a_{2}})), 
e^{2\pi i z_{2}})
q_{\tau}^{L(0)-\frac{c}{24}}\nn
&&=\sum_{a_{4}\in \A}\tr_{W^{a_{4}}}
\Y_{a_{4}}(\mathcal{U}(e^{2\pi iz_{2}})
\Y_{2}(w_{a_{1}}, x_{0})w_{a_{2}}, e^{2\pi i z_{2}})
q_{\tau}^{L(0)-\frac{c}{24}}
\end{eqnarray}
By (\ref{1-independence1}) and (\ref{1-independence1.5}), we have 
\begin{equation}\label{1-independence2}
\sum_{a_{4}\in \A}\tr_{W^{a_{4}}}
\Y_{a_{4}}(\mathcal{U}(e^{2\pi iz_{2}})
\Y_{2}(w_{a_{1}}, x_{0})w_{a_{2}}, e^{2\pi i z_{2}})
q_{\tau}^{L(0)-\frac{c}{24}}=0
\end{equation}
for $w_{a_{1}}\in W^{a_{1}}$ and $w_{a_{2}}\in W^{a_{2}}$. 
Since the coefficients of $\Y_{2}(w_{a_{1}}, x_{0})w_{a_{2}}$ for 
$w_{a_{1}}\in W^{a_{1}}$ and $w_{a_{2}}\in W^{a_{2}}$
span $W^{a_{1}}\boxtimes_{P(z_{1}-z_{2})} W^{a_{2}}$, 
we obtain from (\ref{1-independence2}) that 
\begin{equation}\label{1-independence3}
\sum_{a_{4}\in \A}\tr_{W^{a_{4}}}
\Y_{a_{4}}(\mathcal{U}(e^{2\pi iz_{2}})w, e^{2\pi i z_{2}})
q_{\tau}^{L(0)-\frac{c}{24}}=0
\end{equation}
for $w\in W^{a_{1}}\boxtimes_{P(z_{1}-z_{2})} W^{a_{2}}$
and $w_{a_{3}}\in W^{a_{3}}$. 

Since $W^{a_{4}}$ for $a_{4}\in \A$
are irreducible $V$-modules, we have 
$T(W^{a_{4}})=W^{a_{4}}_{(h_{a_{4}})}$, where
$$T(W^{a_{4}})=\{w\in W^{a_{4}}\;|\; u_{n}w=0, u\in V, \wt u-n-1<0\}$$ 
(see \cite{H9}).
Since $\tau\in \mathbb{H}$ is arbitrary,  from (\ref{1-independence3}),
we have
\begin{equation}\label{1-independence3.5}
\sum_{a_{4}\in \A}
\tr_{T(W^{a_{4}})}
o_{\Y_{a_{4}}}(\mathcal{U}(1)w)=0
\end{equation}
for $w\in W^{a_{1}}\boxtimes_{P(z_{1}-z_{2})} W^{a_{2}}$,
where 
$$o_{\Y_{a_{4}}}(\tilde{w})=(\Y_{a_{4}})_{\wt \hat{w}-1}(\hat{w})$$
for homogeneous $\tilde{w}\in W^{a_{1}}\boxtimes_{P(z_{1}-z_{2})} W^{a_{2}}$
(see Chapter 6 of \cite{H9} for more details).
Since $W^{a_{4}}$ for $a_{4}\in \A$
are inequivalent irreducible $V$-modules, 
$T(W^{a_{4}})$ for $a_{4}\in \A$ 
are inequivalent irreducible $\tilde{A}(V)$-modules
by Proposition 6.5 in \cite{H9}. 
Thus $\tr_{T(W^{a_{4}})}$ for $a_{4}\in \A$ are linearly 
independent. So from (\ref{1-independence3.5}), we obtain 
\begin{equation}\label{1-independence3.6}
o_{\Y_{a_{4}}}(\mathcal{U}(1)w)=0
\end{equation}
for $w\in W^{a_{1}}\boxtimes_{P(z_{1}-z_{2})} W^{a_{2}}$.
Thus $\rho(\Y_{a_{4}})=0$ where 
$$\rho(\Y_{a_{4}}): \tilde{A}(W^{a_{1}}\boxtimes_{P(z_{1}-z_{2})} W^{a_{2}})
\otimes_{\tilde{A}(V)} T(W^{a_{4}})\to T(W^{a_{4}})$$
is given by 
$$\rho(\Y_{a_{4}})((w+\tilde{O}(W^{a_{1}}\boxtimes_{P(z_{1}-z_{2})} W^{a_{2}}))
\otimes w_{a_{4}})=o_{\Y_{a_{4}}}(\mathcal{U}(1)w)w_{a_{4}}$$
for $w\in W^{a_{1}}\boxtimes_{P(z_{1}-z_{2})} W^{a_{2}}$ and 
$w_{a_{4}}\in T(W^{a_{4}})$ (see chapter 6 of \cite{H9}). 
By Theorem 6.9 in \cite{H9}, $\rho(\Y_{a_{4}})=0$ is equivalent to 
$\Y_{a_{4}}=0$ for $a_{4}\in \A$, that is 
\begin{equation}\label{1-independence5}
\Y_{a_{4}}(w, x)=0
\end{equation}
for $a_{4}\in \A$ and $w\in W^{a_{1}}
\boxtimes_{P(z_{1}-z_{2})} W^{a_{2}}$.

For $a_{3}, a_{4}\in \A$ and $l=1, \dots, N_{a_{1}a_{2}}^{a_{3}}$,
take $w$ to be an element of the tensor product module $W^{a_{1}}
\boxtimes_{P(z_{1}-z_{2})} W^{a_{2}}$ such that 
$f(w)$ is in the $l$-th copy of $W^{a_{3}}$ in 
$\oplus_{a_{3}\in \A}N_{a_{1}a_{2}}^{a_{3}}W^{a_{3}}$, that
is, $\pi_{a_{3}; l}(f(w))=f(w)$, 
$\pi_{a_{3}; m}(f(w))=0$ for $m\ne l$ and 
$\pi_{a; m}(f(w))=0$ for $a\ne a_{3}$. Then by (\ref{1-independence5})
and the definition of 
the intertwining operator $\Y_{a_{4}}$, we have 
$$\sum_{k=1}^{N_{a_{3}a_{4}}^{a_{4}}}
\lambda_{a_{3}, a_{4}, k, l}
\Y_{a_{3}a_{4}; k}^{a_{4}; (3)}(f(w), x)=0.$$
Since $f(w)$ is an arbitrary element of 
the $l$-th copy of $W^{a_{3}}$ in 
$\oplus_{a_{3}\in \A}N_{a_{1}a_{2}}^{a_{3}}W^{a_{3}}$,
we obtain
$$\sum_{k=1}^{N_{a_{3}a_{4}}^{a_{4}}}
\lambda_{a_{3}, a_{4}, k, l}
\Y_{a_{3}a_{4}; k}^{a_{4}; (3)}=0$$
for $a_{3}, a_{4}\in \A$ and $l=1, \dots, N_{a_{1}a_{2}}^{a_{3}}$. 
Since $\Y_{a_{3}a_{4}; k}^{a_{4}; (3)}$ for 
$k=1, \dots, N_{a_{3}a_{4}}^{a_{4}}$ are linearly independent, we 
obtain $\lambda_{a_{3}, a_{4}, k, l}=0$ 
for $a_{3}, a_{4}\in \A$, $k=1, \dots, N_{a_{3}a_{4}}^{a_{4}}$
and $l=1, \dots, N_{a_{1}a_{2}}^{a_{3}}$. 
\epfv

We now  introduce 
a space $\F_{1; 2}$ spanned by linear maps of the form
\begin{eqnarray*}
\Psi^{k, l}_{a_{1}, a_{2}, a_{3}}:
\coprod_{a\in \A}W^{a}\otimes W^{a'}&\to &
\mathbb{G}^{e}_{1; 2}\nn
w_{a}\otimes w_{a'}&\mapsto & 
\Psi_{a_{1}, a_{2}, a_{3}}(w_{a}, w_{a'}; z_{1}, z_{2}; \tau)
\end{eqnarray*}
for $a_{1}, a_{2}, a_{3}\in \A$, $k=1, \dots, 
N_{a_{3}a_{1}}^{a_{1}}$, $l=1, \dots, N_{a_{2}a_{2}'}^{a_{3}}$,
where 
$$\Psi^{k, l}_{a_{1}, a_{2}, a_{3}}(w_{a}, w_{a'}; 
z_{1}, z_{2}; \tau)=0$$ 
when $a\ne a_{2}$ and 
\begin{eqnarray*}
\lefteqn{\Psi^{k, l}_{a_{1},  a_{2}, a_{3}}(w_{a_{2}}, w_{a_{2}'}; 
z_{1}, z_{2}; \tau)}\nn
&&=
E(\tr_{W^{a_{1}}}
\Y_{a_{3}a_{1}; k}^{a_{1}; (1)}(\mathcal{U}(e^{2\pi iz_{2}})
\Y_{a_{2}a'_{2}; l}^{a_{3}; (2)}
(w_{a_{2}}, z_{1}-z_{2})w_{a'_{2}}, e^{2\pi i z_{2}})
q_{\tau}^{L(0)-\frac{c}{24}}).
\end{eqnarray*}

Let $\F^{e}_{1; 2}$ be the subspace of $\F_{1; 2}$
spanned by maps of the form $\Psi^{1, 1}_{a_{1}, a_{2}, e}$ for 
$a_{1}, a_{2}\in \A$ and let $\F^{ne}_{1; 2}$ be 
the subspace of $F_{1; 2}$
spanned by maps of the form $\Psi^{k, l}_{a_{1}, a_{2}, a_{3}}$ for 
$a_{1}, a_{2}, a_{3}\in \A$, $a_{3}\ne e$, 
$k=1, \dots, 
N_{a_{3}a_{1}}^{a_{1}}$ and $l=1, \dots, N_{a_{2}a_{2}'}^{a_{3}}$. We
have:

\begin{prop}\label{projection}
The intersection of $\F^{e}_{1; 2}$ and $\F^{ne}_{1; 2}$ 
is $0$. In particular, 
$$F_{1; 2}=\F^{e}_{1; 2}\oplus \F^{ne}_{1; 2}$$
and there exist a projection $\pi: \F_{1; 2}\to \F^{e}_{1; 2}$. 
\end{prop}
\pf
By Proposition \ref{1-independence0}, 
$\Psi^{k, l}_{a_{1}, a_{2}, a_{3}}$
for $a_{1}, a_{2}, a_{3}\in \A$, $k=1, \dots, 
N_{a_{3}a_{1}}^{a_{1}}$ and $l=1, \dots, N_{a_{2}a_{2}'}^{a_{3}}$
are linearly independent. 
Thus the intersection of the space spanned by $\Psi^{1, 1}_{a_{1}, a_{2}, e}$ 
for 
$a_{1}, a_{2}\in \A$ and the space spanned 
by $\Psi^{k, l}_{a_{1}, a_{2}, a_{3}}$ for 
$a_{1}, a_{2}, a_{3}\in \A$, $a_{3}\ne e$, 
$k=1, \dots, 
N_{a_{3}a_{1}}^{a_{1}}$ and 
$l=1, \dots, N_{a_{2}a_{2}'}^{a_{3}}$ are $0$.
\epfv

We define $S: \F^{e}_{1;1}\to \F^{e}_{1;1}$ as follows: For 
$a\in \mathcal{A}$, let
\begin{eqnarray*}
(S(\Psi_{a}))
(u; \tau)&=&\Psi_{a}\left(\left(-\frac{1}{\tau}\right)^{L(0)}u; -\frac{1}{\tau}\right)\nn
&=&\tr_{W^{a}}
\Y_{ea; 1}^{a; (1)}\left(\mathcal{U}(e^{-2\pi i\frac{z_{2}}{\tau}})
\left(-\frac{1}{\tau}\right)^{L(0)}u, e^{-2\pi i \frac{z_{2}}{\tau}}\right)
q_{-\frac{1}{\tau}}^{L(0)-\frac{c}{24}}.
\end{eqnarray*}
Here we have used our convention that 
$$\left(-\frac{1}{\tau}\right)^{L(0)}=e^{(\log (-\frac{1}{\tau}))L(0)}.$$
Note that by the modular invariance of genus-one one-point functions
proved in \cite{M} and \cite{H9}, 
$S(\Psi_{a})$ is indeed in $\F^{e}_{1;1}$. Thus we do obtain maps 
$S: \F^{e}_{1;1}\to \F^{e}_{1;1}$.

Now we define an action of the map $S$  on the space $\F_{1; 2}^{e}$ by 
$$(S(\Psi_{a_{1}, a_{2}, e}^{1, 1}))
(w_{a}, w_{a'}; z_{1}, z_{2}; \tau)=0$$
when $a\ne a_{2}$ and 
\begin{eqnarray*}
\lefteqn{(S(\Psi_{a_{1}, a_{2}, e}^{1, 1}))
(w_{a_{2}}, w_{a'_{2}}; z_{1}, z_{2}; \tau)}\nn
&&=E\Bigg(\tr_{W^{a_{1}}}
\Y_{ea_{1}; 1}^{a_{1}; (1)}\Bigg(\mathcal{U}(e^{-2\pi i\frac{z_{2}}{\tau}})
\left(-\frac{1}{\tau}\right)^{L(0)}\cdot\nn
&&\quad\quad\quad\quad\quad\quad\quad\quad\quad\cdot \Y_{a_{2}a'_{2}; 1}^{e; (2)}
(w_{a_{2}}, z_{1}-z_{2})w_{a'_{2}}, e^{-2\pi i \frac{z_{2}}{\tau}}\Bigg)
q_{-\frac{1}{\tau}}^{L(0)-\frac{c}{24}}\Bigg)\nn
&&=E\Bigg(\tr_{W^{a_{1}}}
\Y_{ea_{1}; 1}^{a_{1}; (1)}\Bigg(\mathcal{U}(e^{-2\pi i\frac{z_{2}}{\tau}})
\cdot\nn
&&\quad \cdot \Y_{a_{2}a'_{2}; 1}^{e; (2)}
\left(\left(-\frac{1}{\tau}\right)^{L(0)}w_{a_{2}}, 
-\frac{1}{\tau}z_{1}-\left(-\frac{1}{\tau}z_{2}\right)\Bigg)
\left(-\frac{1}{\tau}\right)^{L(0)}w_{a'_{2}}, e^{-2\pi i \frac{z_{2}}{\tau}}\right)
\cdot\nn
&&\quad\quad\quad\quad\cdot
q_{-\frac{1}{\tau}}^{L(0)-\frac{c}{24}}\Bigg)\nn
&&=\Psi_{a_{1}, a_{2}, e}^{1, 1}\left(\left(-\frac{1}{\tau}\right)^{L(0)}w_{a},
\left(-\frac{1}{\tau}\right)^{L(0)}w_{a'}; -\frac{1}{\tau}z_{1}, -\frac{1}{\tau}z_{2};
-\frac{1}{\tau}\right).
\end{eqnarray*}

We shall also need the following result:

\begin{prop}
For $a\in A$, $u\in V$, we have
\begin{eqnarray}\label{contrag-1}
\lefteqn{\tr_{W^{a}}
Y_{W^{a}}(\mathcal{U}(e^{2\pi iz})u, 
e^{2\pi i z})
q_{\tau}^{L(0)-\frac{c}{24}}}\nn
&&=\tr_{W^{a'_{4}}}
Y_{W^{a'}}(\mathcal{U}(e^{-2\pi iz})e^{\pi i L(0)}u, 
e^{-2\pi i z})
q_{\tau}^{L(0)-\frac{c}{24}}).
\end{eqnarray}
For $a_{3}, a_{4}\in \A$, $i=1, \dots, N_{a_{1}a_{3}}^{a_{4}}$,
$j=1, \dots, N_{a_{2}a_{4}}^{a_{3}}$, $w_{a_{1}}\in W^{a_{1}}$,
$w_{a_{2}}\in W^{a_{2}}$, we have
\begin{eqnarray}\label{contrag-p}
\lefteqn{E(\tr_{W^{a_{4}}}
\sigma_{23}(\Y_{a_{1}a'_{4}; i}^{a'_{3}; (1)})(\mathcal{U}(e^{2\pi iz_{1}})w_{a_{1}}, 
e^{2\pi i z_{1}})\cdot}\nn
&&\quad\quad\quad\quad\quad\quad\quad\quad\quad\quad\quad
\cdot \sigma_{23}(\Y_{a_{2}a'_{3}; j}^{a'_{4}; (2)})
(\mathcal{U}(e^{2\pi iz_{2}})w_{a_{2}},  e^{2\pi i z_{2}})
q_{\tau}^{L(0)-\frac{c}{24}})\nn
&&=e^{-\pi i(h_{a_{1}}+h_{a_{2}})}E(\tr_{W^{a'_{4}}}
\Y_{a_{2}a'_{3}; j}^{a'_{4}; (2)}
(\mathcal{U}(e^{-2\pi iz_{2}})e^{\pi i L(0)}w_{a_{2}},  e^{-2\pi i z_{2}})
\cdot\nn
&&\quad\quad\quad\quad\quad\quad\quad\quad\quad\quad\quad
\cdot\Y_{a_{1}a'_{4}; i}^{a'_{3}; (1)}(\mathcal{U}(e^{-2\pi iz_{1}})
e^{\pi i L(0)}w_{a_{1}}, 
e^{-2\pi i z_{1}})
q_{\tau}^{L(0)-\frac{c}{24}})\nn
&&
\end{eqnarray}
and
\begin{eqnarray}\label{contrag-i}
\lefteqn{E(\tr_{W^{a'_{4}}}
\sigma_{23}(\Y_{a_{3}a_{4}; k}^{a_{4}; (3)})(\mathcal{U}(e^{2\pi iz_{2}})
\Y_{a_{1}a_{2}; l}^{a_{3}; (4)}
(w_{a_{1}},  z_{1}-z_{2})w_{a_{2}}, 
e^{2\pi i z_{2}})
q_{\tau}^{L(0)-\frac{c}{24}})}\nn
&&=e^{\pi ih_{a_{3}}}
E(\tr_{W^{a_{4}}}
\Y_{a_{3}a_{4}; k}^{a_{4}; (3)}(\mathcal{U}(e^{2\pi iz_{2}})
e^{-\pi i L(0)}\cdot\nn
&&\quad\quad\quad\quad\quad\quad\quad\quad\quad\quad\quad
\cdot\Y_{a_{1}a_{2}; l}^{a_{3}; (4)}
(w_{a_{1}},  z_{1}-z_{2})w_{a_{2}}, 
e^{2\pi i z_{2}})
q_{\tau}^{L(0)-\frac{c}{24}}).\nn
&&
\end{eqnarray}
\end{prop}
\pf
These formulas follows immediately from 
(\ref{contrag-m}).
\epfv

\renewcommand{\theequation}{\thesection.\arabic{equation}}
\renewcommand{\thethm}{\thesection.\arabic{thm}}
\setcounter{equation}{0}
\setcounter{thm}{0}

\section{Properties of fusing and braiding matrices}

In the present section, we prove some properties of 
the fusing and braiding matrices. These properties 
play important role in the proofs of Moore-Seiberg formulas
in the next section and in the proof of the symmetry 
of the matrix associated to the modular transformation 
$\tau\mapsto -1/\tau$ in Section 5. 

In this section, for $p=1, 2, 3, 4, 5, 6$ and $a_{1}, a_{2}, a_{3}
\in \mathcal{A}$, 
$\mathcal{Y}_{a_{1}a_{2}; i}^{a_{3}; (p)}$, $i=1, \dots, 
N_{a_{1}a_{2}}^{a_{3}}$, are bases 
of $\mathcal{V}_{a_{1}a_{2}}^{a_{3}}$.

\begin{prop}
The following equality expressing the squares of 
braiding matrices in terms of the fusing matrices and the 
inverses of fusing matrices holds:
\begin{eqnarray}\label{braiding}
\lefteqn{\sum_{a_{7}\in \A}\sum_{k=1}^{N_{a_{1}a_{2}}^{a_{7}}}
\sum_{l=1}^{N_{a_{7}a_{3}}^{a_{4}}}
F(\Y_{a_{1}a_{5}; i}^{a_{4}; (1)}\otimes \Y_{a_{2}a_{3}; j}^{a_{5}; (2)}; 
\Y_{a_{7}a_{3}; l}^{a_{4}; (5)}
\otimes \Y_{a_{1}a_{2}; k}^{a_{7}; (6)})\cdot}\nn
&&\quad\quad\quad\quad\cdot 
e^{2(2r+1)\pi i (h_{a_{7}}-h_{a_{1}}-h_{a_{2}})}
F^{-1}(\Y_{a_{7}a_{3}; l}^{a_{4}; (5)}\otimes \Y_{a_{1}a_{2}; k}^{a_{7}; (6)}; 
\Y_{a_{1}a_{6}; p}^{a_{4}; (3)}\otimes \Y_{a_{2}a_{3}; q}^{a_{6}; (4)}) \nn
&&=(B^{(r)})^{2}(\Y_{a_{1}a_{5}; i}^{a_{4}; (1)}\otimes \Y_{a_{2}a_{3}; j}^{a_{5}; (2)};
\Y_{a_{1}a_{6}; p}^{a_{4}; (3)}\otimes \Y_{a_{2}a_{3}; q}^{a_{6}; (4)}).
\end{eqnarray}
\end{prop}
\pf
Let $w_{a_{1}}\in W^{a_{1}}$,
$w_{a_{2}}\in W^{a_{2}}$, $w_{a_{3}}\in W^{a_{3}}$ and 
$w_{a'_{4}}\in W^{a'_{4}}$. Then we have
\begin{eqnarray*}
\lefteqn{E(\langle w_{a'_{4}}, \Y_{a_{1}a_{5}; i}^{a_{4}; (1)}(w_{a_{1}}, z_{1})
\Y_{a_{2}a_{3}; j}^{a_{5}; (2)}(w_{a_{2}}, z_{2})w_{a_{3}}\rangle)}\nn
&&=\sum_{a_{7}\in \A}\sum_{k=1}^{N_{a_{1}a_{2}}^{a_{7}}}
\sum_{l=1}^{N_{a_{7}a_{3}}^{a_{4}}}
F(\Y_{a_{1}a_{5}; i}^{a_{4}; (1)}\otimes \Y_{a_{2}a_{3}; j}^{a_{5}; (2)}; 
\Y_{a_{7}a_{3}; l}^{a_{4}; (5)}
\otimes \Y_{a_{1}a_{2}; k}^{a_{7}; (6)})\cdot\nn
&&\quad\quad\quad\quad\cdot 
E(\langle w_{a'_{4}}, 
\Y_{a_{7}a_{3}; l}^{a_{4}; (5)}(\Y_{a_{1}a_{2}; k}^{a_{7}; (6)}(w_{a_{1}}, z_{1}-z_{2})
w_{a_{2}}, z_{2})w_{a_{3}}\rangle).
\end{eqnarray*}
Applying $(B^{(r)})^{2}$ to both sides of the above formula, we obtain
\begin{eqnarray}\label{braiding1}
\lefteqn{(B^{(r)})^{2}(E(\langle w_{a'_{4}}, \Y_{a_{1}a_{5}; i}^{a_{4}; (1)}(w_{a_{1}}, z_{1})
\Y_{a_{2}a_{3}; j}^{a_{5}; (2)}(w_{a_{2}}, z_{2})w_{a_{3}}\rangle))}\nn
&&=\sum_{a_{7}\in \A}\sum_{k=1}^{N_{a_{1}a_{2}}^{a_{7}}}
\sum_{l=1}^{N_{a_{7}a_{3}}^{a_{4}}}
F(\Y_{a_{1}a_{5}; i}^{a_{4}; (1)}\otimes \Y_{a_{2}a_{3}; j}^{a_{5}; (2)}; 
\Y_{a_{7}a_{3}; l}^{a_{4}; (5)}
\otimes \Y_{a_{1}a_{2}; k}^{a_{7}; (6)})\cdot\nn
&&\quad\quad\quad\quad\cdot 
(B^{(r)})^{2}(E(\langle w_{a'_{4}}, 
\Y_{a_{7}a_{3}; l}^{a_{4}; (5)}(\Y_{a_{1}a_{2}; k}^{a_{7}; (6)}(w_{a_{1}}, z_{1}-z_{2})
w_{a_{2}}, z_{2})w_{a_{3}}\rangle)).\nn
&&
\end{eqnarray}

In Section 1, we have seen that $(B^{(r)})^{2}$ is the monodromy 
given by 
$$\log (z_{1}-z_{2})\mapsto \log (z_{1}-z_{2})+ 2(2r+1)\pi i.$$
Using this fact and 
\begin{eqnarray*}
\lefteqn{\Y_{a_{1}a_{2}; k}^{a_{7}; (6)}(w_{a_{1}}, x)
|_{x^{n}=e^{n(\log (z_{1}-z_{2})+2(2r+1)\pi i)}, \; n\in \C}}\nn
&&=e^{2(2r+1)\pi i(h_{a_{7}}-h_{a_{1}}-h_{a_{2}})}
\Y_{a_{1}a_{2}; k}^{a_{7}; (6)}(w_{a_{1}}, x)
|_{x^{n}=e^{n\log (z_{1}-z_{2})},\; n\in \C}\;\;,
\end{eqnarray*}
we have 
\begin{eqnarray}\label{braiding2}
\lefteqn{(B^{(r)})^{2}(E(\langle w_{a'_{4}}, 
\Y_{a_{7}a_{3}; l}^{a_{4}; (5)}(\Y_{a_{1}a_{2}; k}^{a_{7}; (6)}(w_{a_{1}}, z_{1}-z_{2})
w_{a_{2}}, z_{2})w_{a_{3}}\rangle))}\nn
&&=e^{2(2r+1)\pi i (h_{a_{7}}-h_{a_{1}}-h_{a_{2}})}E(\langle w_{a'_{4}}, 
\Y_{a_{7}a_{3}; l}^{a_{4}; (5)}(\Y_{a_{1}a_{2}; k}^{a_{7}; (6)}(w_{a_{1}}, z_{1}-z_{2})
w_{a_{2}}, z_{2})w_{a_{3}}\rangle).\nn
&&
\end{eqnarray}
Using (\ref{braiding1}), (\ref{braiding2}) and the associativity 
(\ref{assoc-inv})
expressed in terms of the matrix elements of the inverse of 
the fusing isomorphism, we obtain
\begin{eqnarray}\label{braiding3}
\lefteqn{(B^{(r)})^{2}(E(\langle w_{a'_{4}}, \Y_{a_{1}a_{5}; i}^{a_{4}; (1)}(w_{a_{1}}, z_{1})
\Y_{a_{2}a_{3}; j}^{a_{5}; (2)}(w_{a_{2}}, z_{2})w_{a_{3}}\rangle))}\nn
&&=\sum_{a_{7}\in \A}\sum_{k=1}^{N_{a_{1}a_{2}}^{a_{7}}}
\sum_{l=1}^{N_{a_{7}a_{3}}^{a_{4}}}
F(\Y_{a_{1}a_{5}; i}^{a_{4}; (1)}\otimes \Y_{a_{2}a_{3}; j}^{a_{5}; (2)}; 
\Y_{a_{7}a_{3}; l}^{a_{4}; (5)}
\otimes \Y_{a_{1}a_{2}; k}^{a_{7}; (6)})\cdot\nn
&&\quad\quad\quad\quad\cdot 
e^{2(2r+1)\pi i (h_{a_{7}}-h_{a_{1}}-h_{a_{2}})}\cdot\nn
&&\quad\quad\quad\quad\cdot E(\langle w_{a'_{4}}, 
\Y_{a_{7}a_{3}; l}^{a_{4}; (5)}(\Y_{a_{1}a_{2}; k}^{a_{7}; (6)}(w_{a_{1}}, z_{1}-z_{2})
w_{a_{2}}, z_{2})w_{a_{3}}\rangle)\nn
&&=\sum_{a_{7}\in \A}\sum_{k=1}^{N_{a_{1}a_{2}}^{a_{7}}}
\sum_{l=1}^{N_{a_{7}a_{3}}^{a_{4}}}
F(\Y_{a_{1}a_{5}; i}^{a_{4}; (1)}\otimes \Y_{a_{2}a_{3}; j}^{a_{5}; (2)}; 
\Y_{a_{7}a_{3}; l}^{a_{4}; (5)}
\otimes \Y_{a_{1}a_{2}; k}^{a_{7}; (6)})
\cdot\nn
&&\quad\quad\quad\quad\cdot 
e^{2(2r+1)\pi i (h_{a_{7}}-h_{a_{1}}-h_{a_{2}})}\cdot\nn
&&\quad\quad\quad\quad\cdot \sum_{a_{6}\in \A}\sum_{p=1}^{N_{a_{1}a_{6}}^{a_{4}}}
\sum_{q=1}^{N_{a_{2}a_{3}}^{a_{6}}}
F^{-1}(\Y_{a_{7}a_{3}; l}^{a_{4}; (5)}\otimes \Y_{a_{1}a_{2}; k}^{a_{7}; (6)}; 
\Y_{a_{1}a_{6}; p}^{a_{4}; (3)}\otimes \Y_{a_{2}a_{3}; q}^{a_{6}; (4)})\cdot\nn
&&\quad\quad\quad\quad\cdot 
E(\langle w_{a'_{4}}, 
\Y_{a_{1}a_{6}; p}^{a_{4}; (3)}(w_{a_{1}}, z_{1})\Y_{a_{2}a_{3}; q}^{a_{6}; (4)}
(w_{a_{2}}, z_{2})w_{a_{3}}\rangle).
\end{eqnarray}
Comparing (\ref{braiding3})
with the definition (\ref{b-r-2}) 
of the matrix elements of $(B^{(r)})^{2}$ and 
using Proposition \ref{independence0}, we
obtain (\ref{braiding}). 
\epfv

\begin{prop}
The following equality  expressing the inverses of 
the fusing matrices in terms of the fusing matrices
holds:
\begin{eqnarray}\label{fusing-inv}
\lefteqn{F^{-1}(\Y_{a_{6}a_{3}; l}^{a_{4}; (1)}
\otimes \Y_{a_{1}a_{2}; k}^{a_{6}; (2)};
\Y_{a_{1}a_{5}; i}^{a_{4}; (3)}\otimes \Y_{a_{2}a_{3}; j}^{a_{5}; (4)})}\nn
&&=F(\sigma_{12}(\Y_{a_{6}a_{3}; l}^{a_{4}; (1)})
\otimes \sigma_{12}(\Y_{a_{1}a_{2}; k}^{a_{6}; (2)});
\sigma_{12}(\Y_{a_{1}a_{5}; i}^{a_{4}; (3)})\otimes 
\sigma_{12}(\Y_{a_{2}a_{3}; j}^{a_{5}; (4)})).
\end{eqnarray}
\end{prop}
\pf
Let $w_{a_{1}}\in W^{a_{1}}$,
$w_{a_{2}}\in W^{a_{2}}$, $w_{a_{3}}\in W^{a_{3}}$ and 
$w_{a'_{4}}\in W^{a'_{4}}$. Then using the definition of 
$\sigma_{12}$, the relation $\sigma_{12}^{2}=1$, the associativity
(\ref{assoc}) and the $L(-1)$-conjugation property (see \cite{FHL}), 
we have
\begin{eqnarray}\label{fusing-inv1}
\lefteqn{E(\langle w_{a'_{4}}, 
\Y_{a_{6}a_{3}; l}^{a_{4}; (1)}(\Y_{a_{1}a_{2}; k}^{a_{6}; (2)}(w_{a_{1}}, z_{1}-z_{2})
w_{a_{2}}, z_{2})w_{a_{3}}\rangle)}\nn
&&=E(\langle w_{a'_{4}}, 
\sigma_{12}^{2}(\Y_{a_{6}a_{3}; l}^{a_{4}; (1)})(
\sigma_{12}^{2}(\Y_{a_{1}a_{2}; k}^{a_{6}; (2)})(w_{a_{1}}, z_{1}-z_{2})
w_{a_{2}}, z_{2})w_{a_{3}}\rangle)\nn
&&=e^{\pi i (h_{a_{4}}-h_{a_{6}}-h_{a_{3}})}e^{\pi i (h_{a_{6}}-h_{a_{1}}-h_{a_{2}})}
\cdot\nn
&&\quad\quad\quad\cdot 
E(\langle w_{a'_{4}}, 
e^{z_{2}L(-1)}\sigma_{12}(\Y_{a_{6}a_{3}; l}^{a_{4}; (1)})(w_{a_{3}}, e^{-\pi i}z_{2})
\cdot\nn
&&\quad\quad\quad\quad\quad\quad\quad\quad\quad\cdot 
e^{(z_{1}-z_{2})L(-1)}\sigma_{12}(\Y_{a_{1}a_{2}; k}^{a_{6}; (2)})
(w_{a_{2}}, e^{-\pi i}(z_{1}-z_{2}))
w_{a_{1}}\rangle)\nn
&&=e^{\pi i (h_{a_{4}}-h_{a_{1}}-h_{a_{2}}-h_{a_{3}})}
E(\langle w_{a'_{4}}, 
e^{z_{1}L(-1)}\sigma_{12}(\Y_{a_{6}a_{3}; l}^{a_{4}; (1)})(w_{a_{3}}, e^{-\pi i}z_{1})
\cdot\nn
&&\quad\quad\quad\quad\quad\quad\quad\quad\quad\quad\quad\quad\quad\quad\cdot 
\sigma_{12}(\Y_{a_{1}a_{2}; k}^{a_{6}; (2)})(w_{a_{2}}, e^{-\pi i}(z_{1}-z_{2}))
w_{a_{1}}\rangle)\nn
&&=e^{\pi i (h_{a_{4}}-h_{a_{1}}-h_{a_{2}}-h_{a_{3}})}
\sum_{a_{5}\in \A}\sum_{i=1}^{N_{a_{1}a_{5}}^{a_{4}}}
\sum_{j=1}^{N_{a_{2}a_{3}}^{a_{5}}}\nn
&&\quad\quad\cdot 
F(\sigma_{12}(\Y_{a_{6}a_{3}; l}^{a_{4}; (1)})\otimes 
\sigma_{12}(\Y_{a_{1}a_{2}; k}^{a_{6}; (2)});
\sigma_{12}(\Y_{a_{1}a_{5}; i}^{a_{4}; (3)})\otimes 
\sigma_{12}(\Y_{a_{2}a_{3}; j}^{a_{5}; (4)}))\cdot\nn
&&\quad\quad\cdot 
E(\langle w_{a'_{4}}, 
e^{z_{1}L(-1)}\sigma_{12}(\Y_{a_{1}a_{5}; i}^{a_{4}; (3)})(
\sigma_{12}(\Y_{a_{2}a_{3}; j}^{a_{5}; (4)})(w_{a_{3}}, e^{-\pi i}z_{2})
\cdot\nn
&&\quad\quad\quad\quad\quad\quad\quad\quad\quad\quad\quad\quad\quad\quad\quad\quad
\quad\quad\quad\cdot w_{a_{2}}, 
e^{-\pi i}(z_{1}-z_{2}))
w_{a_{1}}\rangle)\nn
&&=e^{\pi i (h_{a_{4}}-h_{a_{1}}-h_{a_{5}})}
e^{\pi i (h_{a_{5}}-h_{a_{2}}-h_{a_{3}})}
\sum_{a_{5}\in \A}\sum_{i=1}^{N_{a_{1}a_{5}}^{a_{4}}}
\sum_{j=1}^{N_{a_{2}a_{3}}^{a_{5}}}\nn
&&\quad\quad\cdot 
F(\sigma_{12}(\Y_{a_{6}a_{3}; l}^{a_{4}; (1)})\otimes 
\sigma_{12}(\Y_{a_{1}a_{2}; k}^{a_{6}; (2)});
\sigma_{12}(\Y_{a_{1}a_{5}; i}^{a_{4}; (3)})\otimes 
\sigma_{12}(\Y_{a_{2}a_{3}; j}^{a_{5}; (4)}))\cdot\nn
&&\quad\quad\cdot 
E(\langle w_{a'_{4}}, 
e^{z_{1}L(-1)}\sigma_{12}(\Y_{a_{1}a_{5}; i}^{a_{4}; (3)})(e^{z_{2}L(-1)}\cdot\nn
&&\quad\quad\quad\quad\quad\quad\quad\quad\quad\quad\quad\quad\quad\cdot
\sigma_{12}(\Y_{a_{2}a_{3}; j}^{a_{5}; (4)})(w_{a_{3}}, e^{-\pi i}z_{2})w_{a_{2}}, 
e^{-\pi i}z_{1})
w_{a_{1}}\rangle)\nn
&&=\sum_{a_{5}\in \A}\sum_{i=1}^{N_{a_{1}a_{5}}^{a_{4}}}
\sum_{j=1}^{N_{a_{2}a_{3}}^{a_{5}}}\nn
&&\quad\quad
F(\sigma_{12}(\Y_{a_{6}a_{3}; l}^{a_{4}; (1)})\otimes 
\sigma_{12}(\Y_{a_{1}a_{2}; k}^{a_{6}; (2)});
\sigma_{12}(\Y_{a_{1}a_{5}; i}^{a_{4}; (3)})\otimes 
\sigma_{12}(\Y_{a_{2}a_{3}; j}^{a_{5}; (4)}))\cdot\nn
&&\quad\quad\cdot 
E(\langle w_{a'_{4}}, 
\Y_{a_{1}a_{5}; i}^{a_{4}; (3)}(w_{a_{1}}, 
z_{1})\Y_{a_{2}a_{3}; j}^{a_{5}; (4)}(w_{a_{2}}, z_{2})w_{a_{3}}\rangle).
\end{eqnarray}
Comparing (\ref{fusing-inv1}) with (\ref{assoc-inv2}) and using 
Proposition \ref{independence0}, we obtain 
(\ref{fusing-inv}).
\epfv

In the proof of the next property of fusing matrices, 
we need the following lemma:

\begin{lemma}
For any $z_{1}, z_{2}\in \C^{\times}$ satisfying 
$z_{1}\ne z_{2}$ and any $V$-module $W$, 
the following equalities for maps
from $W$ to $\overline{W}$ holds:
\begin{eqnarray}
&e^{z_{2}L(1)}e^{z_{1}^{-1}L(-1)}=e^{(z_{1}-z_{2})^{-1}L(-1)}e^{z_{1}^{-1}z_{2}(z_{1}-z_{2})L(1)}
(z_{1}(z_{1}-z_{2})^{-1})^{2L(0)},&\label{sl-2-1}\\
&e^{z_{1}^{-1}z_{2}(z_{1}-z_{2})L(1)}(z_{1}(z_{1}-z_{2})^{-1})^{2L(0)}
e^{z_{1}L(1)}(e^{-\pi i}z_{1}^{-2})^{L(0)}&\nn
&=e^{(z_{1}-z_{2})L(1)}(e^{-\pi i}(z_{1}-z_{2})^{-2})^{L(0)}.&
\label{sl-2-2}
\end{eqnarray}
\end{lemma}
\pf
From the identity
$$
\frac{1}{x^{-1}+z_{2}}-z_{1}^{-1}=(z_{1}(z_{1}-z_{2})^{-1})^{-2}
\left(\frac{1}{(x-(z_{1}-z_{2})^{-1})^{-1}+z_{1}^{-1}z_{2}(z_{1}-z_{2})}\right),
$$
we obtain 
\begin{equation}\label{sl-2-1-1}
e^{-z_{2}x^{2}\frac{d}{dx}}e^{-z_{1}^{-1}\frac{d}{dx}}x
=e^{-(z_{1}-z_{2})^{-1}\frac{d}{dx}}e^{-z_{1}^{-1}z_{2}(z_{1}-z_{2})x^{2}\frac{d}{dx}}
(z_{1}(z_{1}-z_{2})^{-1})^{-2x\frac{d}{dx}}x.
\end{equation}
Using (\ref{sl-2-1-1}) and the theory developed in 
Chapters 4 and 5 of \cite{H4}, we obtain (\ref{sl-2-1}).

To prove (\ref{sl-2-2}), we note that the weight of $L(1)$ is $-1$.
So we have
\begin{equation}\label{sl-2-2-0}
(z_{1}(z_{1}-z_{2})^{-1})^{2L(0)}L(1)(z_{1}(z_{1}-z_{2})^{-1})^{-2L(0)}
=z^{-2}_{1}(z_{1}-z_{2})^{2}L(1).
\end{equation}
Using (\ref{sl-2-2-0}), we obtain
\begin{eqnarray}\label{sl-2-2-1}
\lefteqn{(z_{1}(z_{1}-z_{2})^{-1})^{2L(0)}
e^{z_{1}L(1)}}\nn
&&=(z_{1}(z_{1}-z_{2})^{-1})^{2L(0)}
e^{z_{1}L(1)}(z_{1}(z_{1}-z_{2})^{-1})^{-2L(0)}(z_{1}(z_{1}-z_{2})^{-1})^{2L(0)}\nn
&&=e^{z_{1}(z_{1}(z_{1}-z_{2})^{-1})^{2L(0)}L(1)
(z_{1}(z_{1}-z_{2})^{-1})^{-2L(0)}}(z_{1}(z_{1}-z_{2})^{-1})^{2L(0)}\nn
&&=e^{z_{1}^{-1}(z_{1}-z_{2})^{2}L(1)}(z_{1}(z_{1}-z_{2})^{-1})^{2L(0)}.
\end{eqnarray}
The equality (\ref{sl-2-2}) follows immediately from 
(\ref{sl-2-2-1}).
\epfv

Let $\sigma_{123}=\sigma_{12}\sigma_{23}$
and $\sigma_{132}=\sigma_{23}\sigma_{12}$.
We have:

\begin{prop}
The following equality between fusing matrices holds:
\begin{eqnarray}\label{fusing}
\lefteqn{F(\Y_{a_{1}a_{5}; i}^{a_{4}; (1)}\otimes \Y_{a_{2}a_{3}; j}^{a_{5}; (2)}; 
\Y_{a_{6}a_{3}; l}^{a_{4}; (3)}
\otimes \Y_{a_{1}a_{2}; k}^{a_{6}; (4)})}\nn
&&=F(\sigma_{132}(\Y_{a_{2}a_{3}; j}^{a_{5}; (2)})\otimes 
\sigma_{123}(\Y_{a_{1}a_{5}; i}^{a_{4}; (1)}); 
\sigma_{123}(\Y_{a_{1}a_{2}; k}^{a_{6}; (4)})
\otimes \sigma_{132}(\Y_{a_{6}a_{3}; l}^{a_{4}; (3)})).\nn
&&
\end{eqnarray}
\end{prop}
\pf
Let $w_{a_{1}}\in W^{a_{1}}$,
$w_{a_{2}}\in W^{a_{2}}$, $w_{a_{3}}\in W^{a_{3}}$ and 
$w_{a'_{4}}\in W^{a'_{4}}$.  Then using the definitions of 
$\sigma_{12}$ and $\sigma_{23}$ and the relations $\sigma_{12}^{2}=\sigma_{23}^{2}=1$,
we obtain
\begin{eqnarray}\label{fusing-1}
\lefteqn{E(\langle w_{a'_{4}}, \Y_{a_{1}a_{5}; i}^{a_{4}; (1)}(w_{a_{1}}, z_{1})
\Y_{a_{2}a_{3}; j}^{a_{5}; (2)}(w_{a_{2}}, z_{2})w_{a_{3}}\rangle)}\nn
&&=E(\langle w_{a'_{4}}, \sigma_{23}^{2}(\Y_{a_{1}a_{5}; i}^{a_{4}; (1)})(w_{a_{1}}, z_{1})
\sigma_{12}^{2}(\Y_{a_{2}a_{3}; j}^{a_{5}; (2)})(w_{a_{2}}, z_{2})w_{a_{3}}\rangle)\nn
&&=e^{\pi i(h_{a_{1}}+h_{a_{5}}-h_{a_{2}}-h_{a_{3}})}\cdot\nn
&&\quad\quad\cdot 
E(\langle \sigma_{23}(\Y_{a_{1}a_{5}; i}^{a_{4}; (1)})(e^{z_{1}L(1)}(
e^{-\pi i}z_{1}^{-2})^{L(0)}w_{a_{1}}, z^{-1}_{1})w_{a'_{4}}, \nn
&&\quad\quad\quad\quad\quad\quad\quad\quad\quad\quad\quad\quad\quad\quad
e^{z_{2}L(-1)}\sigma_{12}(\Y_{a_{2}a_{3}; j}^{a_{5}; (2)})(w_{a_{3}}, 
e^{-\pi i}z_{2})w_{a_{2}}\rangle)\nn
&&=e^{\pi i(h_{a_{1}}+h_{a_{5}}-h_{a_{2}}-h_{a_{3}})}\cdot\nn
&&\quad\quad\cdot 
E(\langle \sigma_{12}^{2}(\sigma_{23}(\Y_{a_{1}a_{5}; i}^{a_{4}; (1)}))(e^{z_{1}L(1)}(
e^{-\pi i}z_{1}^{-2})^{L(0)}w_{a_{1}}, z^{-1}_{1})w_{a'_{4}}, \nn
&&\quad\quad\quad\quad\quad\quad\quad\quad\quad\quad\quad\quad
e^{z_{2}L(-1)}\sigma_{23}^{2}(\sigma_{12}(\Y_{a_{2}a_{3}; j}^{a_{5}; (2)}))(w_{a_{3}}, 
e^{-\pi i}z_{2})w_{a_{2}}\rangle)\nn
&&=e^{\pi i(2h_{a_{5}}-h_{a_{2}}-h_{a_{4}})}\cdot\nn
&&\quad\quad\cdot 
E(\langle \sigma_{132}(\Y_{a_{2}a_{3}; j}^{a_{5}; (2)})
(e^{-z_{2}L(1)}(e^{-\pi i}z_{2}^{-2})^{L(0)}w_{a_{3}}, 
e^{\pi i}z^{-1}_{2})e^{z_{2}L(1)}e^{z_{1}^{-1}L(-1)}\cdot\nn
&&\quad\quad\quad\quad\quad\quad\cdot
\sigma_{123}(
\Y_{a_{1}a_{5}; i}^{a_{4}; (1)})(w_{a'_{4}}, e^{-\pi i}z^{-1}_{1})
e^{z_{1}L(1)}(e^{-\pi i}z_{1}^{-2})^{L(0)}w_{a_{1}}, 
w_{a_{2}}\rangle).\nn
&&
\end{eqnarray}

By (\ref{sl-2-1}), (\ref{sl-2-2}), the $L(-1)$-, $L(0)$- and 
$L(1)$-conjugation formulas (see \cite{FHL}), (\ref{y-delta})
and the associativity, 
the right-hand side of (\ref{fusing1}) is 
equal to
\begin{eqnarray}\label{fusing-1.1}
\lefteqn{e^{\pi i(2h_{a_{5}}-h_{a_{2}}-h_{a_{4}})}\cdot}\nn
&&\quad\quad\cdot 
E(\langle \sigma_{132}(\Y_{a_{2}a_{3}; j}^{a_{5}; (2)})
(e^{-z_{2}L(1)}(e^{-\pi i}z_{2}^{-2})^{L(0)}w_{a_{3}}, 
e^{\pi i}z^{-1}_{2})\cdot\nn
&&\quad\quad\quad\quad\quad\cdot 
e^{(z_{1}-z_{2})^{-1}L(-1)}e^{z_{1}^{-1}z_{2}(z_{1}-z_{2})L(1)}
(z_{1}(z_{1}-z_{2})^{-1})^{2L(0)}\cdot\nn
&&\quad\quad\quad\quad\quad\quad\quad\cdot 
\sigma_{123}(
\Y_{a_{1}a_{5}; i}^{a_{4}; (1)})(w_{a'_{4}}, e^{-\pi i}z^{-1}_{1})
e^{z_{1}L(1)}(e^{-\pi i}z_{1}^{-2})^{L(0)}w_{a_{1}}, 
w_{a_{2}}\rangle)\nn
&&=e^{\pi i(2h_{a_{5}}-h_{a_{2}}-h_{a_{4}})}\cdot\nn
&&\quad\quad\cdot 
E(\langle e^{(z_{1}-z_{2})^{-1}L(-1)}\sigma_{132}(\Y_{a_{2}a_{3}; j}^{a_{5}; (2)})
(e^{-z_{2}L(1)}(e^{-\pi i}z_{2}^{-2})^{L(0)}\cdot\nn
&&\quad\quad\quad\quad\cdot 
w_{a_{3}}, 
e^{\pi i}z_{1}z^{-1}_{2}(z_{1}-z_{2}))e^{z_{1}^{-1}z_{2}(z_{1}-z_{2})L(1)}\cdot\nn
&&\quad\quad\quad\quad\quad\;\cdot 
\sigma_{123}(\Y_{a_{1}a_{5}; i}^{a_{4}; (1)})(
(z_{1}(z_{1}-z_{2})^{-1})^{2L(0)}w_{a'_{4}}, e^{-\pi i}z^{-1}_{1}(z_{1}-z_{2})^{-2})
\cdot\nn
&&\quad\quad\quad\quad\quad\quad\quad\quad\quad\cdot 
(z_{1}(z_{1}-z_{2})^{-1})^{2L(0)}
e^{z_{1}L(1)}(e^{-\pi i}z_{1}^{-2})^{L(0)}w_{a_{1}}, 
w_{a_{2}}\rangle)\nn
&&=e^{\pi i(2h_{a_{5}}-h_{a_{2}}-h_{a_{4}})}\cdot\nn
&&\quad\quad\cdot 
E(\langle e^{(z_{1}-z_{2})^{-1}L(-1)}\sigma_{132}(\Y_{a_{2}a_{3}; j}^{a_{5}; (2)})
(e^{-z_{2}L(1)}(e^{-\pi i}z_{2}^{-2})^{L(0)}\cdot\nn
&&\quad\quad\quad\cdot w_{a_{3}}, 
e^{\pi i}z_{1}z^{-1}_{2}(z_{1}-z_{2}))
\sigma_{123}(\Y_{a_{1}a_{5}; i}^{a_{4}; (1)})(
e^{z_{2}L(1)}w_{a'_{4}}, e^{-\pi i}(z_{1}-z_{2})^{-1})\cdot\nn
&&\quad\quad\quad\quad\cdot
e^{z_{1}^{-1}z_{2}(z_{1}-z_{2})L(1)}(z_{1}(z_{1}-z_{2})^{-1})^{2L(0)}
e^{z_{1}L(1)}(e^{-\pi i}z_{1}^{-2})^{L(0)}w_{a_{1}}, 
w_{a_{2}}\rangle)\nn
&&=e^{\pi i(2h_{a_{5}}-h_{a_{2}}-h_{a_{4}})}\cdot\nn
&&\quad\quad\cdot 
E(\langle e^{(z_{1}-z_{2})^{-1}L(-1)}\sigma_{132}(\Y_{a_{2}a_{3}; j}^{a_{5}; (2)})
(e^{-z_{2}L(1)}(e^{-\pi i}z_{2}^{-2})^{L(0)}\cdot\nn
&&\quad\quad\quad\quad\quad\quad\quad\quad\quad\quad\quad\quad
\quad\quad\quad\quad\quad\quad\cdot w_{a_{3}}, 
e^{2\pi i}e^{-\pi i}z_{1}z^{-1}_{2}(z_{1}-z_{2}))\cdot\nn
&&\quad\quad\quad\cdot
\sigma_{123}(\Y_{a_{1}a_{5}; i}^{a_{4}; (1)})(
e^{z_{2}L(1)}w_{a'_{4}}, e^{-\pi i}(z_{1}-z_{2})^{-1})\cdot\nn
&&\quad\quad\quad\quad\cdot
e^{(z_{1}-z_{2})L(1)}(e^{-\pi i}(z_{1}-z_{2})^{-2})^{L(0)}w_{a_{1}}, 
w_{a_{2}}\rangle)\nn
&&=e^{\pi i(h_{a_{2}}-2h_{a_{3}}-h_{a_{4}})}\cdot\nn
&&\quad\quad\cdot 
E(\langle e^{(z_{1}-z_{2})^{-1}L(-1)}\sigma_{132}(\Y_{a_{2}a_{3}; j}^{a_{5}; (2)})
(e^{-z_{2}L(1)}(e^{-\pi i}z_{2}^{-2})^{L(0)}\cdot\nn
&&\quad\quad\quad\cdot w_{a_{3}}, 
e^{-\pi i}z_{1}z^{-1}_{2}(z_{1}-z_{2}))
\sigma_{123}(\Y_{a_{1}a_{5}; i}^{a_{4}; (1)})(
e^{z_{2}L(1)}w_{a'_{4}}, e^{-\pi i}(z_{1}-z_{2})^{-1})\cdot\nn
&&\quad\quad\quad\quad\cdot
e^{(z_{1}-z_{2})L(1)}(e^{-\pi i}(z_{1}-z_{2})^{-2})^{L(0)}w_{a_{1}}, 
w_{a_{2}}\rangle)\nn
&&=e^{\pi i(h_{a_{2}}-2h_{a_{3}}-h_{a_{4}})}
\sum_{a_{6}\in \A}\sum_{k=1}^{N_{a_{6}a_{3}}^{a_{4}}}
\sum_{l=1}^{N_{a_{1}a_{2}}^{a_{6}}}\nn
&&\quad
F(\sigma_{132}(\Y_{a_{2}a_{3}; j}^{a_{5}; (2)})\otimes 
\sigma_{123}(\Y_{a_{1}a_{5}; i}^{a_{4}; (1)});
\sigma_{123}(\Y_{a_{1}a_{2}; l}^{a_{6}; (4)})\otimes 
\sigma_{132}(\Y_{a_{6}a_{3}; j}^{a_{4}; (3)}))
\cdot\nn
&&\quad\quad\cdot
E(\langle e^{(z_{1}-z_{2})^{-1}L(-1)}\sigma_{123}(\Y_{a_{1}a_{2}; j}^{a_{6}; (4)})
(\sigma_{132}(\Y_{a_{6}a_{3}; i}^{a_{4}; (3)})
(e^{-z_{2}L(1)}(e^{-\pi i}z_{2}^{-2})^{L(0)}\cdot\nn
&&\quad\quad\quad\cdot w_{a_{3}}, 
e^{-\pi i}z^{-1}_{2})e^{z_{2}L(1)}w_{a'_{4}}, e^{-\pi i}(z_{1}-z_{2})^{-1})\cdot\nn
&&\quad\quad\quad\quad\cdot
e^{(z_{1}-z_{2})L(1)}(e^{-\pi i}(z_{1}-z_{2})^{-2})^{L(0)}w_{a_{1}}, 
w_{a_{2}}\rangle).
\end{eqnarray}
Using the definitions of 
$\sigma_{12}$ and $\sigma_{23}$, the relations $\sigma_{12}^{2}=
\sigma_{23}^{2}=1$, $\sigma_{12}\sigma_{123}=\sigma_{23}$, 
$\sigma_{23}\sigma_{132}=\sigma_{12}$, we see that the right-hand side 
of (\ref{fusing-1.1}) is equal to
\begin{eqnarray}\label{fusing-1.2}
\lefteqn{e^{\pi i(h_{a_{6}}+h_{a_{1}}-2h_{a_{3}}-h_{a_{4}})}
\sum_{a_{6}\in \A}\sum_{k=1}^{N_{a_{6}a_{3}}^{a_{4}}}\sum_{l=1}^{N_{a_{1}a_{2}}^{a_{6}}}}\nn
&&
F(\sigma_{132}(\Y_{a_{2}a_{3}; j}^{a_{5}; (2)})\otimes \sigma_{123}(\Y_{a_{1}a_{5}; i}^{a_{4}; (1)});
\sigma_{123}(\Y_{a_{1}a_{2}; l}^{a_{6}; (4)})\otimes \sigma_{132}(\Y_{a_{6}a_{3}; j}^{a_{4}; (3)}))
\cdot\nn
&&\quad\cdot 
E(\langle \sigma_{12}(\sigma_{123}(\Y_{a_{1}a_{2}; j}^{a_{6}; (4)}))
(e^{(z_{1}-z_{2})L(1)}(e^{-\pi i}(z_{1}-z_{2})^{-2})^{L(0)}\cdot\nn
&&\quad\quad\cdot w_{a_{1}}, (z_{1}-z_{2})^{-1})
\sigma_{132}(\Y_{a_{6}a_{3}; i}^{a_{4}; (3)})(e^{-z_{2}L(1)}(e^{-\pi i}z_{2}^{-2})^{L(0)}w_{a_{3}}, 
e^{-2\pi i}e^{\pi i}z^{-1}_{2})\cdot\nn
&&\quad\quad\quad\quad\quad\quad\quad\quad\quad\quad\quad\quad\quad\quad\quad\quad\quad
\quad\quad\quad\quad\quad\quad\quad\cdot e^{z_{2}L(1)}w_{a'_{4}}, 
w_{a_{2}}\rangle)\nn
&&=e^{\pi i(h_{a_{1}}+h_{a_{4}}-h_{a_{6}})}
\sum_{a_{6}\in \A}\sum_{k=1}^{N_{a_{6}a_{3}}^{a_{4}}}\sum_{l=1}^{N_{a_{1}a_{2}}^{a_{6}}}\nn
&&\quad
F(\sigma_{132}(\Y_{a_{2}a_{3}; j}^{a_{5}; (2)})\otimes \sigma_{123}(\Y_{a_{1}a_{5}; i}^{a_{4}; (1)});
\sigma_{123}(\Y_{a_{1}a_{2}; l}^{a_{6}; (4)})\otimes \sigma_{132}(\Y_{a_{6}a_{3}; j}^{a_{4}; (3)}))
\cdot\nn
&&\quad\quad\cdot 
E(\langle \sigma_{23}(\Y_{a_{1}a_{2}; j}^{a_{6}; (4)})
(e^{(z_{1}-z_{2})L(1)}(e^{-\pi i}(z_{1}-z_{2})^{-2})^{L(0)}w_{a_{1}}, (z_{1}-z_{2})^{-1})\cdot\nn
&&\quad\quad\quad\cdot
\sigma_{132}(\Y_{a_{6}a_{3}; i}^{a_{4}; (3)})(e^{-z_{2}L(1)}(e^{-\pi i}z_{2}^{-2})^{L(0)}w_{a_{3}}, 
e^{\pi i}z^{-1}_{2})e^{z_{2}L(1)}w_{a'_{4}}, 
w_{a_{2}}\rangle)\nn
&&=e^{\pi i(h_{a_{4}}-h_{a_{3}}-h_{a_{6}})}
\sum_{a_{6}\in \A}\sum_{k=1}^{N_{a_{6}a_{3}}^{a_{4}}}\sum_{l=1}^{N_{a_{1}a_{2}}^{a_{6}}}\nn
&&\quad
F(\sigma_{132}(\Y_{a_{2}a_{3}; j}^{a_{5}; (2)})\otimes \sigma_{123}(\Y_{a_{1}a_{5}; i}^{a_{4}; (1)});
\sigma_{123}(\Y_{a_{1}a_{2}; l}^{a_{6}; (4)})\otimes \sigma_{132}(\Y_{a_{6}a_{3}; j}^{a_{4}; (3)}))
\cdot\nn
&&\quad\quad\quad\cdot 
E(\langle w_{a'_{4}}, e^{z_{2}L(-1)}
\sigma_{23}(\sigma_{132}(\Y_{a_{6}a_{3}; i}^{a_{4}; (3)}))(w_{a_{3}}, 
e^{-\pi i}z_{2})\cdot\nn
&&\quad\quad\quad\quad\quad\quad\quad\quad\quad\quad\quad\quad\quad\quad
\quad\quad\quad\cdot
\sigma^{2}_{23}(\Y_{a_{1}a_{2}; j}^{a_{6}; (4)})
(w_{a_{1}}, z_{1}-z_{2})w_{a_{2}}\rangle)\nn
&&=e^{\pi i(h_{a_{4}}-h_{a_{3}}-h_{a_{6}})}
\sum_{a_{6}\in \A}\sum_{k=1}^{N_{a_{6}a_{3}}^{a_{4}}}\sum_{l=1}^{N_{a_{1}a_{2}}^{a_{6}}}\nn
&&\quad
F(\sigma_{132}(\Y_{a_{2}a_{3}; j}^{a_{5}; (2)})\otimes \sigma_{123}(\Y_{a_{1}a_{5}; i}^{a_{4}; (1)});
\sigma_{123}(\Y_{a_{1}a_{2}; l}^{a_{6}; (4)})\otimes \sigma_{132}(\Y_{a_{6}a_{3}; j}^{a_{4}; (3)}))
\cdot\nn
&&\quad\quad\quad\cdot 
E(\langle w_{a'_{4}}, e^{z_{2}L(-1)}
\sigma_{12}(\Y_{a_{6}a_{3}; i}^{a_{4}; (3)})(w_{a_{3}}, 
e^{-\pi i}z_{2})\cdot\nn
&&\quad\quad\quad\quad\quad\quad\quad\quad\quad\quad\quad\quad\quad\quad
\quad\quad\quad\cdot
\sigma^{2}_{23}(\Y_{a_{1}a_{2}; j}^{a_{6}; (4)})
(w_{a_{1}}, z_{1}-z_{2})w_{a_{2}}\rangle)\nn
&&=\sum_{a_{6}\in \A}\sum_{k=1}^{N_{a_{6}a_{3}}^{a_{4}}}\sum_{l=1}^{N_{a_{1}a_{2}}^{a_{6}}}\nn
&&\quad
F(\sigma_{132}(\Y_{a_{2}a_{3}; j}^{a_{5}; (2)})\otimes \sigma_{123}(\Y_{a_{1}a_{5}; i}^{a_{4}; (1)});
\sigma_{123}(\Y_{a_{1}a_{2}; l}^{a_{6}; (4)})\otimes \sigma_{132}(\Y_{a_{6}a_{3}; j}^{a_{4}; (3)}))
\cdot\nn
&&\quad\quad\quad\cdot 
E(\langle w_{a'_{4}}, 
\sigma_{12}^{2}(\Y_{a_{6}a_{3}; i}^{a_{4}; (3)})(\sigma^{2}_{23}(\Y_{a_{1}a_{2}; j}^{a_{6}; (4)})
(w_{a_{1}}, z_{1}-z_{2})w_{a_{2}}, 
z_{2})w_{a_{3}}\rangle)\nn
&&=\sum_{a_{6}\in \A}\sum_{k=1}^{N_{a_{6}a_{3}}^{a_{4}}}\sum_{l=1}^{N_{a_{1}a_{2}}^{a_{6}}}\nn
&&\quad
F(\sigma_{132}(\Y_{a_{2}a_{3}; j}^{a_{5}; (2)})\otimes \sigma_{123}(\Y_{a_{1}a_{5}; i}^{a_{4}; (1)});
\sigma_{123}(\Y_{a_{1}a_{2}; l}^{a_{6}; (4)})\otimes \sigma_{132}(\Y_{a_{6}a_{3}; j}^{a_{4}; (3)}))
\cdot\nn
&&\quad\quad\quad\cdot 
E(\langle w_{a'_{4}}, 
\Y_{a_{6}a_{3}; i}^{a_{4}; (3)}(\Y_{a_{1}a_{2}; j}^{a_{6}; (4)}
(w_{a_{1}}, z_{1}-z_{2})w_{a_{2}}, 
z_{2})w_{a_{3}}\rangle).
\end{eqnarray}

Using (\ref{fusing-1})--(\ref{fusing-1.2}), we see that 
the left-hand side of (\ref{fusing-1}) is equal to the 
right-hand side of (\ref{fusing-1.2}). 
Comparing this resulting equality with (\ref{assoc}) and using 
Proposition \ref{independence0}, we obtain 
(\ref{fusing}). 
\epfv

\renewcommand{\theequation}{\thesection.\arabic{equation}}
\renewcommand{\thethm}{\thesection.\arabic{thm}}
\setcounter{equation}{0}
\setcounter{thm}{0}

\section{Moore-Seiberg formulas}

In \cite{MS1}, Moore and Seiberg derived two formulas from 
the (assumed) axioms for rational conformal 
field theories. The Verlinde conjecture is a consequence
of these equations.  In this section, we prove these formulas mathematically
using the results obtained in the representation theory of 
vertex operator algebras, especially those obtained 
in \cite{H8} and \cite{H9} and in the 
preceding sections.

We now want to choose
a basis $\mathcal{Y}_{a_{1}a_{2}; i}^{a_{3}}$, $i=1, \dots, 
N_{a_{1}a_{2}}^{a_{3}}$, 
of $\mathcal{V}_{a_{1}a_{2}}^{a_{3}}$ for each triple 
$a_{1}, a_{2}, a_{3}\in \mathcal{A}$. Note that for each element 
$\sigma\in S_{3}$, $\sigma(\mathcal{Y})_{a_{1}a_{2}; i}^{a_{3}}$, $i=1, \dots, 
N_{a_{1}a_{2}}^{a_{3}}$, is also a basis of $\mathcal{V}_{a_{1}a_{2}}^{a_{3}}$.

For $a\in \A$, we choose $\Y_{ea; 1}^{a}$ to be the the vertex operator 
$Y_{W^{a}}$ defining the module structure on $W^{a}$ and we choose 
$\Y_{ae; 1}^{a}$ to be the intertwining operator defined using 
the action of $\sigma_{12}$, or equivalently the skew-symmetry 
in this case,
\begin{eqnarray*}
\Y_{ae; 1}^{a}(w_{a}, x)u&=&\sigma_{12}(\Y_{ea; 1}^{a})(w_{a}, x)u\nn
&=&e^{xL(-1)}\Y_{ea; 1}^{a}(u, -x)w_{a}\nn
&=&e^{xL(-1)}Y_{W^{a}}(u, -x)w_{a}
\end{eqnarray*}
for $u\in V$ and $w_{a}\in W^{a}$. 
Since $V'$ as a $V$-module is isomorphic to $V$, we have 
$e'=e$. From \cite{FHL}, we know that there is a nondegerate
invariant 
bilinear form $(\cdot, \cdot)$ on $V$ such that $(\mathbf{1}, 
\mathbf{1})=1$. 
We choose $\Y_{aa'; 1}^{e}=\Y_{aa'; 1}^{e'}$
to be the intertwining operator defined using the action of 
$\sigma_{23}$ by
$$\Y_{aa'; 1}^{e'}=\sigma_{23}(\Y_{ae; 1}^{a}),$$
that is,
$$(u, \Y_{aa'; 1}^{e'}(w_{a}, x)w_{a'})
=e^{\pi i h_{a}}\langle \Y_{ae; 1}^{a}(e^{xL(1)}(e^{-\pi i}x^{-2})^{L(0)}w_{a}, x^{-1})u, 
w_{a'}\rangle$$
for $u\in V$, $w_{a}\in W^{a}$ and $w_{a'}\in W^{a'}$. Since the actions of
$\sigma_{12}$
and $\sigma_{23}$ generate the action of $S_{3}$ on $\mathcal{V}$, we have
$$\Y_{a'a; 1}^{e}=\sigma_{12}(\Y_{aa'; 1}^{e})$$
for any $a\in \mathcal{A}$.
When $a_{1}, a_{2}, a_{3}\ne e$, we choose 
$\mathcal{Y}_{a_{1}a_{2}; i}^{a_{3}}$, $i=1, \dots, 
N_{a_{1}a_{2}}^{a_{3}}$, to be an arbitrary basis
of $\mathcal{V}_{a_{1}a_{2}}^{a_{3}}$.

Recall that $\sigma_{123}=\sigma_{12}\sigma_{23}$ and 
$\sigma_{132}=\sigma_{23}\sigma_{12}$. 
The following theorem gives the first Moore-Seiberg formula 
in \cite{MS1}:

\begin{thm}
For $a_{1}, a_{2}, a_{3}\in \A$, we have 
\begin{eqnarray}\label{formula1-2}
\lefteqn{\sum_{i=1}^{N_{a_{1}a_{2}}^{a_{3}}}\sum_{k=1}^{N_{a'_{1}a_{3}}^{a_{2}}}
F(\Y_{a_{2}e; 1}^{a_{2}}\otimes \Y_{a'_{3}a_{3}; 1}^{e}; 
\Y_{a'_{1}a_{3}; k}^{a_{2}}\otimes \Y_{a_{2}a'_{3}; i}^{a'_{1}})\cdot}\nn
&&\quad\quad\quad\quad\quad\cdot 
F(\Y_{a'_{1}a_{3}; k}^{a_{2}}\otimes \sigma_{123}(\Y_{a_{2}a'_{3}; i}^{a'_{1}});
\Y_{ea_{2}; 1}^{a_{2}}\otimes \Y_{a'_{1}a_{1}; 1}^{e})\nn
&&=N_{a_{1}a_{2}}^{a_{3}}
F(\Y_{a_{2}e; 1}^{a_{2}} \otimes \Y_{a'_{2}a_{2}; 1}^{e};
\Y_{ea_{2}; 1}^{a_{2}}\otimes \Y_{a_{2}a'_{2}; 1}^{e}).
\end{eqnarray}
\end{thm}
\pf 
For $a_{1}, a_{2}, a_{3}\in \A$, $w_{a_{1}}\in W^{a_{1}}$, 
$w_{a'_{3}}\in W^{a'_{3}}$, $w^{1}_{a_{2}}, w^{2}_{a_{2}}\in W^{a_{2}}$,  
$w_{a'_{2}}\in W^{a'_{2}}$, $i=1, \dots, N_{a_{1}a_{2}}^{a_{3}}$,
by the associativity (\ref{assoc}), we have
\begin{eqnarray}\label{pentagon1}
\lefteqn{E(\langle w_{a'_{2}}, \Y_{a_{2}e; 1}^{a_{2}}(w^{1}_{a_{2}}, z_{1})
\Y_{a'_{3}a_{3}; 1}^{e}(w_{a'_{3}}, z_{2})
\sigma_{123}(\Y_{a_{2}a'_{3}; i}^{a'_{1}})(w_{a_{1}}, z_{3})w^{2}_{a_{2}}\rangle)}\nn
&&=\sum_{a_{4}\in \A}\sum_{j=1}^{N_{a_{2}a'_{3}}^{a_{4}}}
\sum_{k=1}^{N_{a_{4}a_{3}}^{a_{2}}}
F(\Y_{a_{2}e; 1}^{a_{2}}\otimes \Y_{a'_{3}a_{3}; 1}^{e}; 
\Y_{a_{4}a_{3}; k}^{a_{2}}\otimes \Y_{a_{2}a'_{3}; j}^{a_{4}})\cdot\nn
&&\quad\cdot 
E(\langle w_{a'_{2}}, 
\Y_{a_{4}a_{3}; k}^{a_{2}}(\Y_{a_{2}a'_{3}; j}^{a_{4}}(w^{1}_{a_{2}}, 
z_{1}-z_{2})w_{a'_{3}}, z_{2})
\sigma_{123}(\Y_{a_{2}a'_{3}; i}^{a'_{1}})(w_{a_{1}}, z_{3})w^{2}_{a_{2}}\rangle)\nn
&&=\sum_{a_{4}\in \A}\sum_{j=1}^{N_{a_{2}a'_{3}}^{a_{4}}}
\sum_{k=1}^{N_{a_{4}a_{3}}^{a_{2}}}
F(\Y_{a_{2}e; 1}^{a_{2}}\otimes \Y_{a'_{3}a_{3}; 1}^{e}; 
\Y_{a_{4}a_{3}; k}^{a_{2}}\otimes \Y_{a_{2}a'_{3}; j}^{a_{4}})\cdot\nn
&&\quad\cdot\sum_{a_{5}\in \A}\sum_{l=1}^{N_{a_{4}a_{1}}^{a_{5}}}
\sum_{m=1}^{N_{a_{5}a_{2}}^{a_{2}}}
F(\Y_{a_{4}a_{3}; k}^{a_{2}}\otimes \sigma_{123}(\Y_{a_{2}a'_{3}; i}^{a'_{1}});
\Y_{a_{5}a_{2}; m}^{a_{2}}\otimes \Y_{a_{4}a_{1}; l}^{a_{5}})\cdot\nn
&&\quad\cdot 
E(\langle w_{a'_{2}}, 
\Y_{a_{5}a_{2}; m}^{a_{2}}(\Y_{a_{4}a_{1}; l}^{a_{5}}
(\Y_{a_{2}a'_{3}; j}^{a_{4}}(w^{1}_{a_{2}}, 
z_{1}-z_{2})w_{a'_{3}}, z_{2}-z_{3})
w_{a_{1}}, z_{3})w^{2}_{a_{2}}\rangle)\nn
&&=\sum_{a_{4}, a_{5}\in \A}
\sum_{j=1}^{N_{a_{2}a'_{3}}^{a_{4}}}
\sum_{k=1}^{N_{a_{4}a_{3}}^{a_{2}}}
\sum_{l=1}^{N_{a_{4}a_{1}}^{a_{5}}}
\sum_{m=1}^{N_{a_{5}a_{2}}^{a_{2}}}
F(\Y_{a_{2}e; 1}^{a_{2}}\otimes \Y_{a'_{3}a_{3}; 1}^{e}; 
\Y_{a_{4}a_{3}; k}^{a_{2}}\otimes \Y_{a_{2}a'_{3}; j}^{a_{4}})\cdot\nn
&&\quad\cdot 
F(\Y_{a_{4}a_{3}; k}^{a_{2}}\otimes \sigma_{123}(\Y_{a_{2}a'_{3}; i}^{a'_{1}});
\Y_{a_{5}a_{2}; m}^{a_{2}}\otimes \Y_{a_{4}a_{1}; l}^{a_{5}})\cdot\nn
&&\quad\cdot 
E(\langle w_{a'_{2}}, 
\Y_{a_{5}a_{2}; m}^{a_{2}}(\Y_{a_{4}a_{1}; l}^{a_{5}}
(\Y_{a_{2}a'_{3}; j}^{a_{4}}(w^{1}_{a_{2}}, 
z_{1}-z_{2})w_{a'_{3}}, z_{2}-z_{3})
w_{a_{1}}, z_{3})w^{2}_{a_{2}}\rangle).\nn
&&
\end{eqnarray}

On the other hand, also by the associativity (\ref{assoc}),
we have
\begin{eqnarray}\label{pentagon2}
\lefteqn{E(\langle w_{a'_{2}}, \Y_{a_{2}e; 1}^{a_{2}}(w^{1}_{a_{2}}, z_{1})
\Y_{a'_{3}a_{3}; 1}^{e}(w_{a'_{3}}, z_{2})
\sigma_{123}(\Y_{a_{2}a'_{3}; i}^{a'_{1}})(w_{a_{1}}, z_{3})w^{2}_{a_{2}}\rangle)}\nn
&&=\sum_{n=1}^{N_{a'_{3}a_{1}}^{a'_{2}}}
F(\Y_{a'_{3}a_{3}; 1}^{e}\otimes \sigma_{123}(\Y_{a_{2}a'_{3}; i}^{a'_{1}});
\Y_{a'_{2}a_{2}; 1}^{e}\otimes \sigma^{2}_{123}(\Y_{a_{2}a'_{3}; n}^{a'_{1}}))\cdot\nn
&&\quad\cdot 
E(\langle w_{a'_{2}}, \Y_{a_{2}e; 1}^{a_{2}}(w^{1}_{a_{2}}, z_{1})
\Y_{a'_{2}a_{2}; 1}^{e}(\sigma^{2}_{123}(
\Y_{a_{2}a'_{3}; n}^{a'_{1}})(w_{a'_{3}}, z_{2}-z_{3})
w_{a_{1}}, z_{3})w^{2}_{a_{2}}\rangle)\nn
&&=\sum_{n=1}^{N_{a'_{3}a_{1}}^{a'_{2}}}
F(\Y_{a'_{3}a_{3}; 1}^{e}\otimes \sigma_{123}(\Y_{a_{2}a'_{3}; i}^{a'_{1}});
\Y_{a'_{2}a_{2}; 1}^{e}\otimes \sigma^{2}_{123}(\Y_{a_{2}a'_{3}; n}^{a'_{1}}))\cdot\nn
&&\quad\cdot 
\sum_{a_{6}\in \A}\sum_{r=1}^{N_{a_{2}a'_{2}}^{a_{6}}}
\sum_{s=1}^{N_{a_{6}a_{2}}^{a_{2}}}
F(\Y_{a_{2}e; 1}^{a_{2}} \otimes \Y_{a'_{2}a_{2}; 1}^{e};
\Y_{a_{6}a_{2}; s}^{a_{2}}\otimes \Y_{a_{2}a'_{2}; r}^{a_{6}})\cdot\nn
&&\quad\quad\cdot 
E(\langle w_{a'_{2}}, \Y_{a_{6}a_{2}; s}^{a_{2}}(
\Y_{a_{2}a'_{2}; r}^{a_{6}}(w^{1}_{a_{2}}, z_{1}-z_{3})\cdot\nn
&&\quad\quad\quad\quad\quad\quad\quad\quad\quad\quad\quad\quad\quad
\cdot 
\sigma^{2}_{123}(\Y_{a_{2}a'_{3}; n}^{a'_{1}})(w_{a'_{3}}, z_{2}-z_{3})
w_{a_{1}}, z_{3})w^{2}_{a_{2}}\rangle)\nn
&&=\sum_{n=1}^{N_{a'_{3}a_{1}}^{a'_{2}}}
F(\Y_{a'_{3}a_{3}; 1}^{e}\otimes \sigma_{123}(\Y_{a_{2}a'_{3}; i}^{a'_{1}});
\Y_{a'_{2}a_{2}; 1}^{e}\otimes \sigma^{2}_{123}(\Y_{a_{2}a'_{3}; n}^{a'_{1}}))\cdot\nn
&&\quad\cdot 
\sum_{a_{6}\in \A}\sum_{r=1}^{N_{a_{2}a'_{2}}^{a_{6}}}
\sum_{s=1}^{N_{a_{6}a_{2}}^{a_{2}}}
F(\Y_{a_{2}e; 1}^{a_{2}} \otimes \Y_{a'_{2}a_{2}; 1}^{e};
\Y_{a_{6}a_{2}; s}^{a_{2}}\otimes \Y_{a_{2}a'_{2}; r}^{a_{6}})\cdot\nn
&&\quad\cdot 
\sum_{a_{7}\in \A}\sum_{p=1}^{N_{a_{2}a'_{3}}^{a_{7}}}
\sum_{q=1}^{N_{a_{7}a_{1}}^{a_{6}}}
F(\Y_{a_{2}a'_{2}; r}^{a_{6}}\otimes \sigma^{2}_{123}(\Y_{a_{2}a'_{3}; n}^{a'_{1}});
\Y_{a_{7}a_{1}; q}^{a_{6}}\otimes  \Y_{a_{2}a'_{3}; p}^{a_{7}})\cdot\nn
&&\quad\cdot 
E(\langle w_{a'_{2}}, \Y_{a_{6}a_{2}; s}^{a_{2}}(
\Y_{a_{7}a_{1}; q}^{a_{6}}(\Y_{a_{2}a'_{3}; p}^{a_{7}}(
w^{1}_{a_{2}}, z_{1}-z_{2})
w_{a'_{3}}, z_{2}-z_{3})
w_{a_{1}}, z_{3})w^{2}_{a_{2}}\rangle)\nn
&&=\sum_{a_{6}, a_{7}\in \A}\sum_{n=1}^{N_{a'_{3}a_{1}}^{a'_{2}}}
\sum_{r=1}^{N_{a_{2}a'_{2}}^{a_{6}}}
\sum_{s=1}^{N_{a_{6}a_{2}}^{a_{2}}}
\sum_{p=1}^{N_{a_{2}a'_{3}}^{a_{7}}}
\sum_{q=1}^{N_{a_{7}a_{1}}^{a_{6}}}\nn
&&\quad \; F(\Y_{a'_{3}a_{3}; 1}^{e}\otimes 
\sigma_{123}(\Y_{a_{2}a'_{3}; i}^{a'_{1}});
\Y_{a'_{2}a_{2}; 1}^{e}\otimes \sigma^{2}_{123}(\Y_{a_{2}a'_{3}; n}^{a'_{1}}))\cdot\nn
&&\quad\cdot 
F(\Y_{a_{2}e; 1}^{a_{2}} \otimes \Y_{a'_{2}a_{2}; 1}^{e};
\Y_{a_{6}a_{2}; s}^{a_{2}}\otimes \Y_{a_{2}a'_{2}; r}^{a_{6}})\cdot\nn
&&\quad\cdot 
F(\Y_{a_{2}a'_{2}; r}^{a_{6}}\otimes \sigma^{2}_{123}(\Y_{a_{2}a'_{3}; n}^{a'_{1}});
\Y_{a_{7}a_{1}; q}^{a_{6}}\otimes  \Y_{a_{2}a'_{3}; p}^{a_{7}})\cdot\nn
&&\quad\cdot 
E(\langle w_{a'_{2}}, \Y_{a_{6}a_{2}; s}^{a_{2}}(
\Y_{a_{7}a_{1}; q}^{a_{6}}(\Y_{a_{2}a'_{3}; p}^{a_{7}}(
w^{1}_{a_{2}}, z_{1}-z_{2})
w_{a'_{3}}, z_{2}-z_{3})
w_{a_{1}}, z_{3})w^{2}_{a_{2}}\rangle).\nn
&&
\end{eqnarray}

Using (\ref{pentagon1})--(\ref{pentagon2}) and Proposition 
\ref{independence-3p}, we obtain
\begin{eqnarray*}
\lefteqn{\sum_{k=1}^{N_{a_{4}a_{3}}^{a_{2}}}
F(\Y_{a_{2}e; 1}^{a_{2}}\otimes \Y_{a'_{3}a_{3}; 1}^{e}; 
\Y_{a_{4}a_{3}; k}^{a_{2}}\otimes \Y_{a_{2}a'_{3}; j}^{a_{4}})\cdot}\nn
&&\quad\quad\quad\quad\quad\cdot 
F(\Y_{a_{4}a_{3}; k}^{a_{2}}\otimes \sigma_{123}(\Y_{a_{2}a'_{3}; i}^{a'_{1}});
\Y_{a_{5}a_{2}; m}^{a_{2}}\otimes \Y_{a_{4}a_{1}; l}^{a_{5}})\nn
&&=\sum_{n=1}^{N_{a'_{3}a_{1}}^{a'_{2}}}
\sum_{r=1}^{N_{a_{2}a'_{2}}^{a_{5}}}
F(\Y_{a'_{3}a_{3}; 1}^{e}\otimes \sigma_{123}(\Y_{a_{2}a'_{3}; i}^{a'_{1}});
\Y_{a'_{2}a_{2}; 1}^{e}\otimes \sigma^{2}_{123}(\Y_{a_{2}a'_{3}; n}^{a'_{1}}))\cdot\nn
&&\quad\quad\quad\quad\quad\cdot 
F(\Y_{a_{2}e; 1}^{a_{2}} \otimes \Y_{a'_{2}a_{2}; 1}^{e};
\Y_{a_{5}a_{2}; m}^{a_{2}}\otimes \Y_{a_{2}a'_{2}; r}^{a_{5}})\cdot\nn
&&\quad\quad\quad\quad\quad\cdot 
F(\Y_{a_{2}a'_{2}; r}^{a_{5}}\otimes \sigma^{2}_{123}(\Y_{a_{2}a'_{3}; n}^{a'_{1}});
\Y_{a_{4}a_{1}; l}^{a_{5}}\otimes  \Y_{a_{2}a'_{3}; j}^{a_{4}}).
\end{eqnarray*}
In particular, for $a_{5}=e$ and $a_{4}=a'_{1}$,
we have
\begin{eqnarray}\label{pentagon3}
\lefteqn{\sum_{k=1}^{N_{a'_{1}a_{3}}^{a_{2}}}
F(\Y_{a_{2}e; 1}^{a_{2}}\otimes \Y_{a'_{3}a_{3}; 1}^{e}; 
\Y_{a'_{1}a_{3}; k}^{a_{2}}\otimes \Y_{a_{2}a'_{3}; j}^{a'_{1}})\cdot}\nn
&&\quad\quad\quad\quad\quad\cdot 
F(\Y_{a_{4}a_{3}; k}^{a_{2}}\otimes \sigma_{123}(\Y_{a_{2}a'_{3}; i}^{a'_{1}});
\Y_{ea_{2}; m}^{a_{2}}\otimes \Y_{a'_{1}a_{1}; l}^{e})\nn
&&=\sum_{n=1}^{N_{a'_{3}a_{1}}^{a'_{2}}}
F(\Y_{a'_{3}a_{3}; 1}^{e}\otimes \sigma_{123}(\Y_{a_{2}a'_{3}; i}^{a'_{1}});
\Y_{a'_{2}a_{2}; 1}^{e}\otimes \sigma^{2}_{123}(\Y_{a_{2}a'_{3}; n}^{a'_{1}}))\cdot\nn
&&\quad\quad\quad\quad\quad\cdot 
F(\Y_{a_{2}e; 1}^{a_{2}} \otimes \Y_{a'_{2}a_{2}; 1}^{e};
\Y_{ea_{2}; 1}^{a_{2}}\otimes \Y_{a_{2}a'_{2}; 1}^{e})\cdot\nn
&&\quad\quad\quad\quad\quad\cdot 
F(\Y_{a_{2}a'_{2}; 1}^{e}\otimes \sigma^{2}_{123}(\Y_{a_{2}a'_{3}; n}^{a'_{1}});
\Y_{a'_{1}a_{1}; 1}^{e}\otimes  \Y_{a_{2}a'_{3}; j}^{a'_{1}})\nn
&&=\Bigg(\sum_{n=1}^{N_{a'_{3}a_{1}}^{a'_{2}}}
F(\Y_{a'_{3}a_{3}; 1}^{e}\otimes \sigma_{123}(\Y_{a_{2}a'_{3}; i}^{a'_{1}});
\Y_{a'_{2}a_{2}; 1}^{e}\otimes \sigma^{2}_{123}(\Y_{a_{2}a'_{3}; n}^{a'_{1}}))\cdot\nn
&&\quad\quad\quad\quad\quad\cdot 
F(\Y_{a_{2}a'_{2}; 1}^{e}\otimes \sigma^{2}_{123}(\Y_{a_{2}a'_{3}; n}^{a'_{1}});
\Y_{a'_{1}a_{1}; 1}^{e}\otimes  \Y_{a_{2}a'_{3}; j}^{a'_{1}})\Bigg)\cdot\nn
&&\quad\quad\quad\quad\quad\cdot 
F(\Y_{a_{2}e; 1}^{a_{2}} \otimes \Y_{a'_{2}a_{2}; 1}^{e};
\Y_{ea_{2}; 1}^{a_{2}}\otimes \Y_{a_{2}a'_{2}; 1}^{e}).
\end{eqnarray}

On the other hand, by the definition of $\sigma_{12}$ and $\sigma_{23}$,
the relations $\sigma_{12}^{2}=\sigma_{23}^{2}=1$ and the choices
of the bases $\Y_{ea; 1}^{a}$, $\Y_{ae; 1}^{a}$, $\Y_{aa'}^{e}$
and $\Y_{a'a;1}^{e}$ for $a\in \A$, we have
\begin{eqnarray}\label{fusing1}
\lefteqn{\langle u, \Y_{a'_{3}a_{3}; 1}^{e}(w_{a'_{3}}, x_{1})
\sigma_{123}(\Y_{a_{2}a'_{3}; i}^{a'_{1}})(w_{a_{1}}, x_{2})w_{a_{2}}\rangle}\nn
&&=\langle u, \sigma_{23}^{2}(\Y_{a'_{3}a_{3}; 1}^{e})(w_{a'_{3}}, x_{1})
\sigma_{23}^{2}(\sigma_{123}(
\Y_{a_{2}a'_{3}; i}^{a'_{1}}))(w_{a_{1}}, x_{2})w_{a_{2}}\rangle\nn
&&=e^{\pi i(h_{a_{1}}+h_{a_{3}})}
\langle \sigma_{23}(\sigma_{123}(
\Y_{a_{2}a'_{3}; i}^{a'_{1}}))(e^{x_{2}L(1)}(e^{-\pi i}x_{2}^{-2})^{L(0)}
w_{a_{1}}, x_{2}^{-1})\cdot\nn
&&\quad\quad\quad\quad\quad\quad\quad\quad\cdot 
\sigma_{23}(\Y_{a'_{3}a_{3}; 1}^{e})(e^{x_{1}L(1)}(e^{-\pi i}x_{1}^{-2})^{L(0)}w_{a'_{3}}, 
x^{-1}_{1})u, w_{a_{2}}\rangle\nn
&&=e^{\pi i(h_{a_{1}}+h_{a_{3}})}
\langle \sigma_{23}(\sigma_{123}(
\Y_{a_{2}a'_{3}; i}^{a'_{1}}))(e^{x_{2}L(1)}(e^{-\pi i}x_{2}^{-2})^{L(0)}
w_{a_{1}}, x_{2}^{-1})\cdot\nn
&&\quad\quad\quad\quad\quad\quad\quad\quad\cdot 
\Y_{a'_{3}e; 1}^{a'_{3}}(e^{x_{1}L(1)}(e^{-\pi i}x_{1}^{-2})^{L(0)}w_{a'_{3}}, 
x^{-1}_{1})u, w_{a_{2}}\rangle\nn
&&=e^{\pi i(h_{a_{1}}+h_{a_{3}})}
\langle \sigma_{12}(\sigma_{123}(
\Y_{a_{2}a'_{3}; i}^{a'_{1}}))(e^{x_{2}L(1)}((e^{-\pi i}x_{2}^{-2})^{L(0)}
w_{a_{1}}, x_{2}^{-1})\cdot\nn
&&\quad\quad\quad\quad\quad\quad\quad\quad\cdot e^{x_{1}^{-1}L(-1)}
\Y_{ea'_{3}; 1}^{a'_{3}}(u, 
-x^{-1}_{1})e^{x_{1}L(1)}((e^{-\pi i}x_{1}^{-2})^{L(0)}w_{a'_{3}}, w_{a_{2}}\rangle\nn
&&=e^{\pi i(h_{a_{1}}+h_{a_{3}})}\cdot\nn
&&\quad\quad\cdot\langle e^{x_{1}^{-1}L(-1)}
\sigma_{23}(\sigma_{123}(
\Y_{a_{2}a'_{3}; i}^{a'_{1}}))(e^{x_{2}L(1)}(e^{-\pi i}x_{2}^{-2})^{L(0)}
w_{a_{1}}, x_{2}^{-1}-x_{1}^{-1})\cdot\nn
&&\quad\quad\quad\quad\quad\quad\quad\quad\cdot
Y_{W^{a'_{3}}}(u, 
-x^{-1}_{1})e^{x_{1}L(1)}(e^{-\pi i}x_{1}^{-2})^{L(0)}w_{a'_{3}}, w_{a_{2}}\rangle
\end{eqnarray}
By the locality between vertex operators on $V$-modules and intertwining 
operators and the definition of $\sigma_{23}$, 
there exists a positive integer $N$ such that 
\begin{eqnarray}\label{fusing2}
\lefteqn{x_{2}^{-N}e^{\pi i(h_{a_{1}}+h_{a_{3}})}}\nn
&&\quad\cdot
\langle e^{x_{1}^{-1}L(-1)}
\sigma_{23}(\sigma_{123}(
\Y_{a_{2}a'_{3}; i}^{a'_{1}}))(e^{x_{2}L(1)}(e^{-\pi i}x_{2}^{-2})^{L(0)}
w_{a_{1}}, x_{2}^{-1}-x_{1}^{-1})\cdot\nn
&&\quad\quad\quad\quad\quad\quad\quad\quad\cdot
Y_{W^{a'_{3}}}(u, 
-x^{-1}_{1})e^{x_{1}L(1)}(e^{-\pi i}x_{1}^{-2})^{L(0)}w_{a'_{3}}, w_{a_{2}}\rangle\nn
&&=((x_{2}^{-1}-x_{1}^{-1})-(-x_{1}^{-1})^{N}e^{\pi i(h_{a_{1}}+h_{a_{3}})}\cdot\nn
&&\quad\cdot
\langle e^{x_{1}^{-1}L(-1)}
\sigma_{23}(\sigma_{123}(
\Y_{a_{2}a'_{3}; i}^{a'_{1}}))(e^{x_{2}L(1)}(e^{-\pi i}x_{2}^{-2})^{L(0)}
w_{a_{1}}, x_{2}^{-1}-x_{1}^{-1})\cdot\nn
&&\quad\quad\quad\quad\quad\quad\quad\quad\quad\quad\quad\quad\cdot
Y_{W^{a'_{3}}}(u, 
-x^{-1}_{1})e^{x_{1}L(1)}(e^{-\pi i}x_{1}^{-2})^{L(0)}w_{a'_{3}}, w_{a_{2}}\rangle\nn
&&=((x_{2}^{-1}-x_{1}^{-1})-(-x_{1}^{-1})^{N}e^{\pi i(h_{a_{1}}+h_{a_{3}})}\cdot\nn
&&\quad\quad \cdot
\langle e^{x_{1}^{-1}L(-1)}
Y_{W_{a'_{2}}}(u, 
-x^{-1}_{1})\cdot\nn
&&\quad\quad\quad\quad\quad\quad\quad\cdot 
\sigma_{23}(\sigma_{123}(
\Y_{a_{2}a'_{3}; i}^{a'_{1}}))(e^{x_{2}L(1)}(e^{-\pi i}x_{2}^{-2})^{L(0)}
w_{a_{1}}, z_{2}^{-1}-z_{1}^{-1})\cdot\nn
&&\quad\quad\quad\quad\quad\quad\quad\quad\quad\quad\quad\quad\cdot 
e^{z_{1}L(1)}(e^{-\pi i}z_{1}^{-2})^{L(0)}w_{a'_{3}}, w_{a_{2}}\rangle\nn
&&=x_{2}^{-N}e^{\pi i(h_{a_{1}}+h_{a_{3}})}\cdot\nn
&&\quad\cdot\langle 
\sigma_{12}(Y_{W_{a'_{2}}})(
\sigma_{23}(\sigma_{123}(
\Y_{a_{2}a'_{3}; i}^{a'_{1}}))(e^{x_{2}L(1)}(e^{-\pi i}x_{2}^{-2})^{L(0)}
w_{a_{1}}, x_{2}^{-1}-x_{1}^{-1})\cdot\nn
&&\quad\quad\quad\quad\quad\quad\quad\quad\quad\quad\quad\quad\cdot 
e^{x_{1}L(1)}(e^{-\pi i}x_{1}^{-2})^{L(0)}w_{a'_{3}}, 
x^{-1}_{1})u, w_{a_{2}}\rangle\nn
&&=x_{2}^{-N}e^{\pi i(h_{a_{1}}+h_{a_{3}}-h_{a_{2}})}\langle 
e^{\pi ih_{a_{1}}}\Y_{a'_{2}e; 1}^{a'_{2}}(e^{x_{1}L(1)}(e^{-\pi i}x_{1}^{-2})^{L(0)}
(e^{\pi i}x_{1}^{2})^{L(0)}e^{-x_{1}L(1)}\cdot\nn
&&\quad\quad\quad\quad\quad\quad\quad\cdot 
\sigma_{23}(\sigma_{123}(
\Y_{a_{2}a'_{3}; i}^{a'_{1}}))(e^{x_{2}L(1)}(e^{-\pi i}x_{2}^{-2})^{L(0)}
w_{a_{1}}, x_{2}^{-1}-x_{1}^{-1})\cdot\nn
&&\quad\quad\quad\quad\quad\quad\quad\quad\quad\quad\quad\quad\cdot 
e^{x_{1}L(1)}(e^{-\pi i}x_{1}^{-2})^{L(0)}w_{a'_{3}}, 
x^{-1}_{1})u, w_{a_{2}}\rangle\nn
&&=x_{2}^{-N}e^{\pi i(h_{a_{1}}+h_{a_{3}}-h_{a_{2}})}
\langle u, \sigma_{23}(
\Y_{a'_{2}e; 1}^{a'_{2}})((e^{\pi i}x_{1}^{2})^{L(0)}e^{-x_{1}L(1)}\cdot\nn
&&\quad\quad\quad\quad\quad\quad\quad\cdot 
\sigma_{23}(\sigma_{123}(
\Y_{a_{2}a'_{3}; i}^{a'_{1}}))(e^{x_{2}L(1)}(e^{-\pi i}x_{2}^{-2})^{L(0)}
w_{a_{1}}, x_{2}^{-1}-x_{1}^{-1})\cdot\nn
&&\quad\quad\quad\quad\quad\quad\quad\quad\quad\quad\quad\quad\cdot 
e^{x_{1}L(1)}(e^{-\pi i}x_{1}^{-2})^{L(0)}w_{a'_{3}}, 
x_{1})w_{a_{2}}\rangle\nn
&&=x_{2}^{-N}e^{\pi i(h_{a_{1}}+h_{a_{3}}-h_{a_{2}})}
\langle u, \Y_{a'_{2}a_{2}; 1}^{e}((e^{-\pi i}x_{1}^{2})^{L(0)}e^{-x_{1}L(1)}\cdot\nn
&&\quad\quad\quad\quad\quad\quad\quad\cdot 
\sigma_{23}(\sigma_{123}(
\Y_{a_{2}a'_{3}; i}^{a'_{1}}))(e^{x_{2}L(1)}(e^{-\pi i}x_{2}^{-2})^{L(0)}
w_{a_{1}}, x_{2}^{-1}-x_{1}^{-1})\cdot\nn
&&\quad\quad\quad\quad\quad\quad\quad\quad\quad\quad\quad\quad\cdot 
e^{x_{1}L(1)}(e^{-\pi i}x_{1}^{-2})^{L(0)}w_{a'_{3}}, 
x_{1})w_{a_{2}}\rangle\nn
&&=x_{2}^{-N}e^{\pi i(h_{a_{1}}+h_{a_{3}}-h_{a_{2}})}\cdot\nn
&&\quad\quad\cdot 
\langle u, \Y_{a'_{2}a_{2}; 1}^{e}(
\sigma_{23}(\sigma_{123}(\Y_{a_{2}a'_{3}; i}^{a'_{1}}))(
w_{a_{1}}, e^{\pi i}(x_{1}-x_{2}))w_{a'_{3}}, 
x_{1})w_{a_{2}}\rangle\nn
&&=x_{2}^{-N}
\langle u, \Y_{a'_{2}a_{2}; 1}^{e}(e^{(x_{2}-x_{1})L(-1)}
e^{-\pi i \Delta(\sigma_{23}(\sigma_{123}(
\Y_{a_{2}a'_{3}; i}^{a'_{1}})))}e^{(x_{1}-x_{2})L(-1)}\cdot\nn
&&\quad\quad\quad\quad\quad\quad\quad\quad\quad\quad\cdot 
\sigma_{23}(\sigma_{123}(\Y_{a_{2}a'_{3}; i}^{a'_{1}}))(
w_{a_{1}}, e^{\pi i}(x_{1}-x_{2}))w_{a'_{3}}, 
z_{1})w_{a_{2}}\rangle\nn
&&=x_{2}^{-N}\langle u, \Y_{a'_{2}a_{2}; 1}^{e}(e^{(x_{2}-x_{1})L(-1)}\cdot\nn
&&\quad\quad\quad\quad\quad\quad\quad\quad\quad\quad\cdot 
\sigma_{12}(\sigma_{23}(\sigma_{123}(\Y_{a_{2}a'_{3}; i}^{a'_{1}})))(
w_{a_{1}}, x_{1}-x_{2})w_{a'_{3}}, 
x_{1})w_{a_{2}}\rangle.\nn
&&
\end{eqnarray}
Using (\ref{fusing1}), (\ref{fusing2}), the $L(-1)$-derivative 
property for $\sigma_{12}(\sigma_{23}(\sigma_{123}(\Y_{a_{2}a'_{3}; i}^{a'_{1}})))$
and the equality $\sigma_{12}\sigma_{23}=\sigma_{123}$,
we obtain
\begin{eqnarray*}
\lefteqn{E(\langle u, \Y_{a'_{3}a_{3}; 1}^{e}(w_{a'_{3}}, z_{1})
\sigma_{123}(\Y_{a_{2}a'_{3}; i}^{a'_{1}})(w_{a_{1}}, z_{2})w_{a_{2}}\rangle)}\nn
&&=E(\langle u, \Y_{a'_{2}a_{2}; 1}^{e}(e^{(z_{2}-z_{1})L(-1)}\cdot\nn
&&\quad\quad\quad\quad\quad\quad\quad\quad\quad\quad\cdot 
\sigma_{12}(\sigma_{23}(\sigma_{123}(\Y_{a_{2}a'_{3}; i}^{a'_{1}})))(
w_{a_{1}}, z_{1}-z_{2})w_{a'_{3}}, 
z_{1})w_{a_{2}}\rangle)\nn
&&=E(\langle u, \Y_{a'_{2}a_{2}; 1}^{e}(
\sigma^{2}_{123}(\Y_{a_{2}a'_{3}; i}^{a'_{1}})(
w_{a_{1}}, z_{1}-z_{2})w_{a'_{3}}, 
z_{2})w_{a_{2}}\rangle)
\end{eqnarray*}
Thus we obtain
\begin{equation}\label{fusing3}
F(\Y_{a'_{3}a_{3}; 1}^{e}\otimes \sigma_{123}(\Y_{a_{2}a'_{3}; i}^{a'_{1}});
\Y_{a'_{2}a_{2}; 1}^{e}\otimes \sigma^{2}_{123}(\Y_{a_{2}a'_{3}; n}^{a'_{1}}))
=\delta_{in}.
\end{equation}
Similarly, we can prove 
\begin{equation}\label{fusing4}
F(\Y_{a_{2}a'_{2}; 1}^{e}\otimes \sigma^{2}_{123}(\Y_{a_{2}a'_{3}; n}^{a'_{1}});
\Y_{a'_{1}a_{1}; 1}^{e}\otimes  \Y_{a_{2}a'_{3}; j}^{a'_{1}})=\delta_{nj}
\end{equation}

Using (\ref{fusing3}) and (\ref{fusing4}), we see that (\ref{pentagon3})
becomes
\begin{eqnarray}\label{formula1}
\lefteqn{\sum_{k=1}^{N_{a'_{1}a_{3}}^{a_{2}}}
F(\Y_{a_{2}e; 1}^{a_{2}}\otimes \Y_{a'_{3}a_{3}; 1}^{e}; 
\Y_{a'_{1}a_{3}; k}^{a_{2}}\otimes \Y_{a_{2}a'_{3}; j}^{a'_{1}})\cdot}\nn
&&\quad\quad\quad\quad\quad\cdot 
F(\Y_{a_{1}'a_{3}; k}^{a_{2}}\otimes \sigma_{123}(\Y_{a_{2}a'_{3}; i}^{a'_{1}});
\Y_{ea_{2}; 1}^{a_{2}}\otimes \Y_{a'_{1}a_{1}; 1}^{e})\nn
&&=\delta_{ij}
F(\Y_{a_{2}e; 1}^{a_{2}} \otimes \Y_{a'_{2}a_{2}; 1}^{e};
\Y_{ea_{2}; 1}^{a_{2}}\otimes \Y_{a_{2}a'_{2}; 1}^{e}).
\end{eqnarray}
Summing over $i=1, \dots, \N_{a_{1}a_{3}}^{a_{2}}$ 
on both sides in the special case $j=i$ of 
(\ref{formula1}), we obtain (\ref{formula1-2}).
\epfv

We now prepare to prove the second formula. Recall from Section 2 the 
maps $\Psi_{a_{1}, a_{2}, e}^{1, 1}: 
\coprod_{a\in \A}W^{a}\otimes W^{a'}\to \mathbb{G}^{e}_{1;2}$
for $a_{1}, a_{2}\in \A$ and the projection 
$\pi: \mathcal{F}_{1;2}\to \mathcal{F}_{1;2}^{e}$. 
For any $f\in \mathcal{F}_{1;2}$, we shall,
for convenience,
denote $(\pi(f))(w_{a}\otimes w_{a'})$
by $\pi(f(w_{a}\otimes w_{a'}))$.

We need the following lemma:

\begin{lemma}
For $a_{1}, a_{2}\in \A$, $w_{a_{2}}\in W^{a_{2}}$ 
and $w_{a'_{2}}\in W^{a'_{2}}$, we have
\begin{eqnarray}\label{alpha0}
\lefteqn{\Psi_{a_{1}, a_{2}, e}^{1, 1}
(w_{a_{2}}, w_{a'_{2}}; z_{1}, z_{2}-1; \tau)}\nn
&&=\sum_{a_{3}\in \mathcal{A}}\sum_{i=1}^{N_{a_{2}a_{3}}^{a_{1}}}
\sum_{j=1}^{N_{a'_{2}a_{1}}^{a_{3}}}
\sum_{a_{4}\in \mathcal{A}}\sum_{k=1}^{N_{a_{4}a_{1}}^{a_{1}}}
\sum_{l=1}^{N_{a_{2}a'_{2}}^{a_{4}}}e^{-2\pi i(h_{a_{3}}-h_{a_{1}})}\cdot\nn
&&\quad\quad \quad \cdot F^{-1}(\Y_{ea_{1}; 1}^{a_{1}}\otimes \Y_{a_{2}a'_{2}; 1}^{e};
\Y_{a_{2}a_{3}; i}^{a_{1}}\otimes \Y_{a'_{2}a_{1}; j}^{a_{3}})
\cdot \nn
&&\quad\quad \quad  \cdot
F(\Y_{a_{2}a_{3}; i}^{a_{1}}\otimes \Y_{a'_{2}a_{1}; j}^{a_{3}};
\Y_{a_{4}a_{1}; k}^{a_{1}}\otimes \Y_{a_{2}a'_{2}; l}^{a_{4}})\cdot \nn
&&\quad\quad\quad \cdot
E\bigg(\tr_{W^{a_{1}}}
\Y_{a_{4}a_{1}; k}^{a_{1}}(\mathcal{U}(e^{2\pi iz_{2}})\cdot \nn
&&\quad\quad\quad \quad\quad\quad\quad\quad\quad\quad\cdot
\Y_{a_{2}a'_{2}; l}^{a_{4}}(w_{a_{2}},z_{1}-z_{2})w_{a'_{2}},  e^{2\pi i z_{2}})
q_{\tau}^{L(0)-\frac{c}{24}}\bigg)\nn
&&
\end{eqnarray}
and 
\begin{eqnarray}\label{beta0}
\lefteqn{\Psi_{a_{1}, a_{2}, e}^{1, 1}
(w_{a_{2}}, w_{a'_{2}}; z_{1}, z_{2}+\tau; \tau)}\nn
&&=\sum_{a_{3}\in \mathcal{A}}\sum_{i=1}^{N_{a_{2}a_{3}}^{a_{1}}}
\sum_{j=1}^{N_{a'_{2}a_{1}}^{a_{3}}}
\sum_{a_{4}\in \mathcal{A}}\sum_{k=1}^{N_{a_{4}a_{3}}^{a_{3}}}
\sum_{l=1}^{N_{a_{2}a'_{2}}^{a_{4}}}e^{\pi i (-2h_{a_{2}}+h_{a_{4}})}\nn
&&\quad \quad F^{-1}(\Y_{ea_{1}; 1}^{a_{1}}\otimes \Y_{a_{2}a'_{2}; 1}^{e};
\sigma_{23}(\Y_{a_{2}a'_{1}; i}^{a'_{3}}) \otimes 
\sigma_{13}(\Y_{a'_{3}a_{1}; j}^{a_{2}}))\cdot \nn
&&\quad\quad \cdot 
F(\Y_{a_{2}a'_{1}; i}^{a'_{3}}\otimes \sigma_{123}(\Y_{a'_{3}a_{1}; j}^{a_{2}});
\Y_{a_{4}a'_{3}; k}^{a'_{3}}\otimes \Y_{a_{2}a'_{2}; l}^{a_{4}})\cdot \nn
&&\quad\quad \cdot 
E\bigg(\tr_{W^{a_{3}}}
\Y_{a_{4}a_{3}; k}^{a_{3}}(\mathcal{U}(e^{2\pi iz_{2}})
\Y_{a_{2}a'_{2}; l}^{a_{4}}(
w_{a_{2}}, z_{1}-z_{2})
w_{a'_{2}}, e^{2\pi iz_{2}})
q_{\tau}^{L(0)-\frac{c}{24}}
\bigg).\nn
&&
\end{eqnarray}
In particular,
for any $a_{1}, a_{2}\in \A$, the maps from 
$\coprod_{a\in \A}W^{a}\otimes W^{a'}$ to the space of
single-valued analytic functions on $\widetilde{M}_{1}^{2}$
given by 
\begin{eqnarray*}
w_{a}\otimes w_{a'}&\mapsto & 
\Psi_{a_{1}, a_{2}, e}^{1,1}
(w_{a}, w_{a'}; z_{1}, z_{2}-1; \tau),\\
w_{a}\otimes w_{a'}&\mapsto & 
\Psi_{a_{1}, a_{2}, e}^{1,1}
(w_{a}, w_{a'}; z_{1}, z_{2}+\tau; \tau)
\end{eqnarray*}
for $w_{a}\in W^{a}$ and $w_{a'}\in W^{a'}$ are
in $\mathcal{F}_{1;2}$.
\end{lemma}
\pf
Using the definition of $\Psi_{a_{1}, a_{2}, e}^{1,1}$, 
the associativity properties (\ref{assoc}), (\ref{assoc-inv}) and
(\ref{y-delta}), we have
\begin{eqnarray*}
\lefteqn{\Psi_{a_{1}, a_{2}, e}^{1, 1}
(w_{a_{2}}, w_{a'_{2}}; z_{1}, z_{2}-1; \tau)}\nn
&&=E\bigg(\tr_{W^{a_{1}}}
\Y_{ea_{1}; 1}^{a_{1}}(\mathcal{U}(e^{2\pi i(z_{2}-1)})\cdot\nn
&&\quad\quad\quad \quad\quad\quad\quad\cdot
\Y_{a_{2}a'_{2}; 1}^{e}
(w_{a_{2}}, z_{1}-(z_{2}-1))w_{a'_{2}}, e^{2\pi i (z_{2}-1)})
q_{\tau}^{L(0)-\frac{c}{24}}\bigg)\nn
&&=\sum_{a_{3}\in \mathcal{A}}
\sum_{i=1}^{N_{a_{2}a_{3}}^{a_{1}}}\sum_{j=1}^{N_{a'_{2}a_{1}}^{a_{3}}}
F^{-1}(\Y_{ea_{1}; 1}^{a_{1}}\otimes \Y_{a_{2}a'_{2}; 1}^{e};
\Y_{a_{2}a_{3}; i}^{a_{1}}\otimes \Y_{a'_{2}a_{1}; j}^{a_{3}})\cdot \nn
&&\quad\quad\quad \cdot 
E\bigg(\tr_{W^{a_{1}}}
\Y_{a_{2}a_{3}; i}^{a_{1}}(\mathcal{U}(e^{2\pi iz_{1}})w_{a_{2}}, e^{2\pi i z_{1}})\cdot \nn
&&\quad\quad\quad\quad\quad\quad\quad\quad\quad\quad\cdot 
\Y_{a'_{2}a_{1}; j}^{a_{3}}(\mathcal{U}(e^{2\pi i(z_{2}-1)})w_{a'_{2}}, e^{2\pi i (z_{2}-1)})
q_{\tau}^{L(0)-\frac{c}{24}}\bigg)\nn
&&=\sum_{a_{3}\in \mathcal{A}}\sum_{i=1}^{N_{a_{2}a_{3}}^{a_{1}}}
\sum_{j=1}^{N_{a'_{2}a_{1}}^{a_{3}}}
F^{-1}(\Y_{ea_{1}; 1}^{a_{1}}\otimes \Y_{a_{2}a'_{2}; 1}^{e};
\Y_{a_{2}a_{3}; i}^{a_{1}}\otimes \Y_{a'_{2}a_{1}; j}^{a_{3}})\cdot \nn
&&\quad\quad\quad \cdot 
E\bigg(\tr_{W^{a_{1}}}
\Y_{a_{2}a_{3}; i}^{a_{1}}(\mathcal{U}(e^{2\pi iz_{1}})w_{a_{2}}, e^{2\pi i z_{1}})\cdot \nn
&&\quad\quad\quad \quad\quad\quad\quad\quad\quad\quad\cdot
\Y_{a'_{2}a_{1}; j}^{a_{3}}(\mathcal{U}(e^{-2\pi i}
e^{2\pi iz_{2}})w_{a'_{2}},  e^{-2\pi i}e^{2\pi i z_{2}})
q_{\tau}^{L(0)-\frac{c}{24}}\bigg)\nn
&&=\sum_{a_{3}\in \mathcal{A}}\sum_{i=1}^{N_{a_{2}a_{3}}^{a_{1}}}
\sum_{j=1}^{N_{a'_{2}a_{1}}^{a_{3}}}e^{-2\pi i(h_{a_{3}}-h_{a_{1}})}
F^{-1}(\Y_{ea_{1}; 1}^{a_{1}}\otimes \Y_{a_{2}a'_{2}; 1}^{e};
\Y_{a_{2}a_{3}; i}^{a_{1}}\otimes \Y_{a'_{2}a_{1}; j}^{a_{3}})
\cdot \nn
&&\quad\quad\quad \cdot
E\bigg(\tr_{W^{a_{1}}}
\Y_{a_{2}a_{3}; i}^{a_{1}}(\mathcal{U}(e^{2\pi iz_{1}})w_{a_{2}}, e^{2\pi i z_{1}})\cdot \nn
&&\quad\quad\quad \quad\quad\quad\quad\quad\quad\quad\cdot
\Y_{a'_{2}a_{1}; j}^{a_{3}}(\mathcal{U}(
e^{2\pi iz_{2}})w_{a'_{2}},  e^{2\pi i z_{2}})
q_{\tau}^{L(0)-\frac{c}{24}}\bigg)\nn
&&=\sum_{a_{3}\in \mathcal{A}}\sum_{i=1}^{N_{a_{2}a_{3}}^{a_{1}}}
\sum_{j=1}^{N_{a'_{2}a_{1}}^{a_{3}}}e^{-2\pi i(h_{a_{3}}-h_{a_{1}})}
F^{-1}(\Y_{ea_{1}; 1}^{a_{1}}\otimes \Y_{a_{2}a'_{2}; 1}^{e};
\Y_{a_{2}a_{3}; i}^{a_{1}}\otimes \Y_{a'_{2}a_{1}; j}^{a_{3}})
\cdot \nn
&&\quad \cdot
\sum_{a_{4}\in \mathcal{A}}\sum_{k=1}^{N_{a_{4}a_{1}}^{a_{1}}}
\sum_{l=1}^{N_{a_{2}a'_{2}}^{a_{4}}}
F(\Y_{a_{2}a_{3}; i}^{a_{1}}\otimes \Y_{a'_{2}a_{1}; j}^{a_{3}};
\Y_{a_{4}a_{1}; k}^{a_{1}}\otimes \Y_{a_{2}a'_{2}; l}^{a_{4}})\cdot \nn
&&\quad\quad\quad \cdot
E\bigg(\tr_{W^{a_{1}}}
\Y_{a_{4}a_{1}; k}^{a_{1}}(\mathcal{U}(e^{2\pi iz_{2}})\cdot \nn
&&\quad\quad\quad \quad\quad\quad\quad\quad\quad\quad\cdot
\Y_{a_{2}a'_{2}; l}^{a_{4}}(w_{a_{2}},z_{1}-z_{2})w_{a'_{2}},  e^{2\pi i z_{2}})
q_{\tau}^{L(0)-\frac{c}{24}}\bigg),
\end{eqnarray*}
proving (\ref{alpha0}).

The proof of (\ref{beta0}) is more complicated. 
Using the definition of $\Psi_{a_{1}, a_{2}, e}^{1,1}$, 
the associativity properties (\ref{assoc}), (\ref{assoc-inv}),
the $L(0)$-conjugation property, the property of traces,  and
(\ref{y-delta}), 
\begin{eqnarray}\label{beta1}
\lefteqn{\Psi_{a_{1}, a_{2}, e}^{1, 1}
(w_{a_{2}}, w_{a'_{2}}; z_{1}, z_{2}+\tau; \tau)}\nn
&&=E\bigg(\tr_{W^{a_{1}}}
\Y_{ea_{1}; 1}^{a_{1}}(\mathcal{U}(e^{2\pi i(z_{2}+\tau)})\cdot\nn
&&\quad\quad\quad \quad\quad\quad\quad\cdot
\Y_{a_{2}a'_{2}; 1}^{e}
(w_{a_{2}}, z_{1}-(z_{2}+\tau))w_{a'_{2}}, e^{2\pi i (z_{2}+\tau)})
q_{\tau}^{L(0)-\frac{c}{24}}\bigg)\nn
&&=\sum_{a_{3}\in \mathcal{A}}\sum_{i=1}^{N_{a_{2}a_{3}}^{a_{1}}}
\sum_{j=1}^{N_{a'_{2}a_{1}}^{a_{3}}}
F^{-1}(\Y_{ea_{1}; 1}^{a_{1}}\otimes \Y_{a_{2}a'_{2}; 1}^{e};
\sigma_{23}(\Y_{a_{2}a'_{1}; i}^{a'_{3}})\otimes 
\sigma_{13}(\Y_{a'_{3}a_{1}; j}^{a_{2}}))\cdot \nn
&&\quad\quad\quad \cdot 
E\bigg(\tr_{W^{a_{1}}}
\sigma_{23}(\Y_{a_{2}a'_{1}; i}^{a'_{3}})(
\mathcal{U}(e^{2\pi iz_{1}})w_{a_{2}}, e^{2\pi i z_{1}})\cdot \nn
&&\quad\quad\quad\quad\quad\quad\quad\quad\quad\quad\cdot 
\sigma_{13}(\Y_{a'_{3}a_{1}; j}^{a_{2}})
(\mathcal{U}(e^{2\pi i(z_{2}+\tau)})w_{a'_{2}}, e^{2\pi i (z_{2}+\tau)})
q_{\tau}^{L(0)-\frac{c}{24}}\bigg)\nn
&&=\sum_{a_{3}\in \mathcal{A}}\sum_{i=1}^{N_{a_{2}a_{3}}^{a_{1}}}
\sum_{j=1}^{N_{a'_{2}a_{1}}^{a_{3}}}
F^{-1}(\Y_{ea_{1}; 1}^{a_{1}}\otimes \Y_{a_{2}a'_{2}; 1}^{e};
\sigma_{23}(\Y_{a_{2}a'_{1}; i}^{a'_{3}}) \otimes 
\sigma_{13}(\Y_{a'_{3}a_{1}; j}^{a_{2}}))\cdot \nn
&&\quad\quad\quad \cdot 
E\bigg(\tr_{W^{a_{1}}}
\sigma_{23}(\Y_{a_{2}a'_{1}; i}^{a'_{3}})(
\mathcal{U}(e^{2\pi iz_{1}})w_{a_{2}}, e^{2\pi i z_{1}})\cdot \nn
&&\quad\quad\quad\quad\quad\quad\quad\quad\quad\quad\cdot 
\sigma_{13}(\Y_{a'_{3}a_{1}; j}^{a_{2}})
(\mathcal{U}(q_{\tau}e^{2\pi iz_{2}})w_{a'_{2}}, q_{\tau}e^{2\pi iz_{2}})
q_{\tau}^{L(0)-\frac{c}{24}}\bigg)\nn
&&=\sum_{a_{3}\in \mathcal{A}}\sum_{i=1}^{N_{a_{2}a_{3}}^{a_{1}}}
\sum_{j=1}^{N_{a'_{2}a_{1}}^{a_{3}}}
F^{-1}(\Y_{ea_{1}; 1}^{a_{1}}\otimes \Y_{a_{2}a'_{2}; 1}^{e};
\sigma_{23}(\Y_{a_{2}a'_{1}; i}^{a'_{3}}) \otimes 
\sigma_{13}(\Y_{a'_{3}a_{1}; j}^{a_{2}}))\cdot \nn
&&\quad\quad\quad \cdot 
E\bigg(\tr_{W^{a_{1}}}
\sigma_{23}(\Y_{a_{2}a'_{1}; i}^{a'_{3}})(
\mathcal{U}(e^{2\pi iz_{1}})w_{a_{2}}, e^{2\pi i z_{1}})\cdot \nn
&&\quad\quad\quad\quad\quad\quad\quad\quad\quad\quad\cdot 
q_{\tau}^{L(0)-\frac{c}{24}}
\sigma_{13}(\Y_{a'_{3}a_{1}; j}^{a_{2}})(
\mathcal{U}(e^{2\pi iz_{2}})w_{a'_{2}}, e^{2\pi iz_{2}})
\bigg)\nn
&&=\sum_{a_{3}\in \mathcal{A}}\sum_{i=1}^{N_{a_{2}a_{3}}^{a_{1}}}
\sum_{j=1}^{N_{a'_{2}a_{1}}^{a_{3}}}
F^{-1}(\Y_{ea_{1}; 1}^{a_{1}}\otimes \Y_{a_{2}a'_{2}; 1}^{e};
\sigma_{23}(\Y_{a_{2}a'_{1}; i}^{a'_{3}}) \otimes 
\sigma_{13}(\Y_{a'_{3}a_{1}; j}^{a_{2}}))\cdot \nn
&&\quad\quad\quad \cdot 
E\bigg(\tr_{W^{a_{3}}}
\sigma_{13}(\Y_{a'_{3}a_{1}; j}^{a_{2}})(
\mathcal{U}(e^{2\pi iz_{2}})w_{a'_{2}}, e^{2\pi iz_{2}})\cdot \nn
&&\quad\quad\quad\quad\quad\quad\quad\quad\quad\quad\cdot 
\sigma_{23}(\Y_{a_{2}a'_{1}; i}^{a'_{3}})(
\mathcal{U}(e^{2\pi iz_{1}})w_{a_{2}}, e^{2\pi i z_{1}})
q_{\tau}^{L(0)-\frac{c}{24}}
\bigg).
\end{eqnarray}

Using (\ref{contrag-p}), the relations $\sigma_{23}^{2}=1$,
$\sigma_{23}\sigma_{13}=\sigma_{123}$ and the genus-one
associativity, we have
\begin{eqnarray}\label{beta2}
\lefteqn{E\bigg(\tr_{W^{a_{3}}}
\sigma_{13}(\Y_{a'_{3}a_{1}; j}^{a_{2}})(
\mathcal{U}(e^{2\pi iz_{2}})w_{a'_{2}}, e^{2\pi iz_{2}})\cdot} \nn
&&\quad\quad\quad\quad\quad\quad\quad\quad\quad\quad\cdot 
\sigma_{23}(\Y_{a_{2}a'_{1}; i}^{a'_{3}})(
\mathcal{U}(e^{2\pi iz_{1}})w_{a_{2}}, e^{2\pi i z_{1}})
q_{\tau}^{L(0)-\frac{c}{24}}
\bigg)\nn
&&=E\bigg(\tr_{W^{a_{3}}}
\sigma_{23}^{2}(\sigma_{13}(\Y_{a'_{3}a_{1}; j}^{a_{2}}))(
\mathcal{U}(e^{2\pi iz_{2}})w_{a'_{2}}, e^{2\pi iz_{2}})\cdot \nn
&&\quad\quad\quad\quad\quad\quad\quad\quad\quad\quad\cdot 
\sigma_{23}(\Y_{a_{2}a'_{1}; i}^{a'_{3}})(
\mathcal{U}(e^{2\pi iz_{1}})w_{a_{2}}, e^{2\pi i z_{1}})
q_{\tau}^{L(0)-\frac{c}{24}}
\bigg)\nn
&&=e^{-2\pi ih_{a_{2}}}
E\bigg(\tr_{W^{a'_{3}}}
\Y_{a_{2}a'_{1}; i}^{a'_{3}}\bigg(
\mathcal{U}(e^{2\pi iz_{1}})e^{\pi iL(0)}
w_{a_{2}}, e^{-2\pi i z_{1}}\bigg)\cdot \nn
&&\quad\quad\quad\quad\quad\quad\quad\cdot 
\sigma_{23}(\sigma_{13}(\Y_{a'_{3}a_{1}; j}^{a_{2}}))
\bigg(\mathcal{U}(e^{2\pi iz_{2}})e^{\pi iL(0)}w_{a'_{2}}, e^{-2\pi iz_{2}}\bigg)
q_{\tau}^{L(0)-\frac{c}{24}}
\bigg)\nn
&&=e^{-2\pi ih_{a_{2}}}
E\bigg(\tr_{W^{a'_{3}}}
\Y_{a_{2}a'_{1}; i}^{a'_{3}}\bigg(\mathcal{U}(e^{-2\pi iz_{1}})e^{\pi i L(0)}
w_{a_{2}}, e^{-2\pi i z_{1}}\bigg)\cdot \nn
&&\quad\quad\quad\quad\quad\quad\quad\cdot 
\sigma_{123}(\Y_{a'_{3}a_{1}; j}^{a_{2}})\bigg(\mathcal{U}(e^{-2\pi iz_{2}})e^{\pi i L(0)}
w_{a'_{2}}, e^{-2\pi iz_{2}}\bigg)
q_{\tau}^{L(0)-\frac{c}{24}}
\bigg)\nn
&&
\end{eqnarray}

We now prove 
\begin{eqnarray}\label{beta3}
\lefteqn{E\bigg(\tr_{W^{a'_{3}}}
\Y_{a_{2}a'_{1}; i}^{a'_{3}}\bigg(\mathcal{U}(e^{-2\pi iz_{1}})e^{\pi i L(0)}
w_{a_{2}}, e^{-2\pi i z_{1}}\bigg)\cdot} \nn
&&\quad\quad\quad\quad\quad\quad\quad\cdot 
\sigma_{123}(\Y_{a'_{3}a_{1}; j}^{a_{2}})\bigg(\mathcal{U}(e^{-2\pi iz_{2}})e^{\pi i L(0)}
w_{a'_{2}}, e^{-2\pi iz_{2}}\bigg)
q_{\tau}^{L(0)-\frac{c}{24}}
\bigg)\nn
&&=\sum_{a_{4}\in \mathcal{A}}\sum_{k=1}^{N_{a_{4}a_{3}}^{a_{3}}}
\sum_{l=1}^{N_{a_{2}a'_{2}}^{a_{4}}}
F(\Y_{a_{2}a'_{1}; i}^{a'_{3}}\otimes \sigma_{123}(\Y_{a'_{3}a_{1}; j}^{a_{2}});
\Y_{a_{4}a'_{3}; k}^{a'_{3}}\otimes \Y_{a_{2}a'_{2}; l}^{a_{4}})\cdot \nn
&&\quad\quad\quad\cdot E\bigg(\tr_{W^{a'_{3}}}
\Y_{a_{4}a'_{3}; k}^{a'_{3}}\bigg(\mathcal{U}(e^{-2\pi iz_{2}})
\Y_{a_{2}a'_{2}; l}^{a_{4}}\bigg(e^{\pi i L(0)}
w_{a_{2}}, e^{\pi i}(z_{1}-z_{2})\bigg)\cdot \nn
&&\quad\quad\quad\quad\quad\quad\quad\quad\quad
\quad\quad\quad\quad\quad\quad\quad\quad\quad\cdot 
e^{\pi i L(0)}
w_{a'_{2}}, e^{-2\pi iz_{2}}\bigg)
q_{\tau}^{L(0)-\frac{c}{24}}
\bigg).\nn
&&
\end{eqnarray}
To prove (\ref{beta3}), we need only prove 
their restrictions to a subregion of 
$\widetilde{M}_{1}^{2}$ are equal. 
So we need only prove that
\begin{eqnarray}\label{beta3.5}
\lefteqn{\tr_{W^{a'_{3}}}
\Y_{a_{2}a'_{1}; i}^{a'_{3}}\bigg(\mathcal{U}(e^{-2\pi iz_{1}})e^{\pi i L(0)}
w_{a_{2}}, e^{-2\pi i z_{1}}\bigg)\cdot} \nn
&&\quad\quad\quad\quad\quad\quad\quad\cdot 
\sigma_{123}(\Y_{a'_{3}a_{1}; j}^{a_{2}})\bigg(\mathcal{U}(e^{-2\pi iz_{2}})e^{\pi i L(0)}
w_{a'_{2}}, e^{-2\pi iz_{2}}\bigg)
q_{\tau}^{L(0)-\frac{c}{24}}\nn
&&=\sum_{a_{4}\in \mathcal{A}}\sum_{k=1}^{N_{a_{4}a_{3}}^{a_{3}}}
\sum_{l=1}^{N_{a_{2}a'_{2}}^{a_{4}}}
F(\Y_{a_{2}a'_{1}; i}^{a'_{3}}\otimes \sigma_{123}(\Y_{a'_{3}a_{1}; j}^{a_{2}});
\Y_{a_{4}a'_{3}; k}^{a'_{3}}\otimes \Y_{a_{2}a'_{2}; l}^{a_{4}})\cdot \nn
&&\quad\quad\quad\cdot \tr_{W^{a'_{3}}}
\Y_{a_{4}a'_{3}; k}^{a'_{3}}\bigg(\mathcal{U}(e^{-2\pi iz_{2}})
\Y_{a_{2}a'_{2}; l}^{a_{4}}\bigg(e^{\pi i L(0)}
w_{a_{2}}, e^{\pi i}(z_{1}-z_{2})\bigg)\cdot \nn
&&\quad\quad\quad\quad\quad\quad\quad\quad\quad
\quad\quad\quad\quad\quad\quad\quad\quad\quad\cdot 
e^{\pi i L(0)}
w_{a'_{2}}, e^{-2\pi iz_{2}}\bigg)
q_{\tau}^{L(0)-\frac{c}{24}}
\bigg)\nn
&&
\end{eqnarray}
holds when $|q_{\tau}|<|e^{-2\pi iz_{2}}|<|e^{-2\pi iz_{1}}|<1$
and $0<|e^{2\pi i(-z_{1}+z_{2})}-1|<1$. 
From (\ref{1-assoc}), we see  that in this region the left-hand side of 
(\ref{beta3.5}) is equal to 
\begin{eqnarray}\label{beta3.7}
\lefteqn{\sum_{a_{4}\in \mathcal{A}}\sum_{k=1}^{N_{a_{4}a_{3}}^{a_{3}}}
\sum_{l=1}^{N_{a_{2}a'_{2}}^{a_{4}}}
F(\Y_{a_{2}a'_{1}; i}^{a'_{3}}\otimes \sigma_{123}(\Y_{a'_{3}a_{1}; j}^{a_{2}});
\Y_{a_{4}a'_{3}; k}^{a'_{3}}\otimes \Y_{a_{2}a'_{2}; l}^{a_{4}})\cdot} \nn
&&\quad\quad\quad\cdot \tr_{W^{a'_{3}}}
\Y_{a_{4}a'_{3}; k}^{a'_{3}}\bigg(\mathcal{U}(e^{-2\pi iz_{2}})
\Y_{a_{2}a'_{2}; l}^{a_{4}}\bigg(e^{\pi i L(0)}
w_{a_{2}}, (-z_{1}+z_{2})\bigg)\cdot \nn
&&\quad\quad\quad\quad\quad\quad\quad\quad\quad
\quad\quad\quad\quad\quad\quad\quad\quad\quad\cdot 
e^{\pi i L(0)}
w_{a'_{2}}, e^{-2\pi iz_{2}}\bigg)
q_{\tau}^{L(0)-\frac{c}{24}}
\bigg)\nn
&&
\end{eqnarray}
Now in this region, because 
$|e^{-2\pi iz_{2}}|<|e^{-2\pi iz_{1}}|$, 
the imaginary part of $z_{1}$
must be bigger than the imaginary part of $z_{2}$. 
Thus $z_{1}-z_{2}$ is in the upper half plane. 
This means that  $\arg (z_{1}-z_{2})<\pi$ and 
$\arg (z_{1}-z_{2}) +\pi <2\pi$. So we have
$\arg (-(z_{1}-z_{2}))= \arg (z_{1}-z_{2})+\pi$.
Now for any $n\in \C$, by our convention,
\begin{eqnarray*}
(-z_{1}+z_{2})^{n}&=&e^{n\log (-z_{1}+z_{2})}\nn
&=&e^{n\log (-(z_{1}-z_{2}))}\nn
&=&e^{n(\log |-(z_{1}-z_{2})|+i\arg (-(z_{1}-z_{2})))}\nn
&=&e^{n(\log |(z_{1}-z_{2})|+i\arg (z_{1}-z_{2})+\pi i)}\nn
&=&e^{n(\log (z_{1}-z_{2})+\pi i)}\nn
&=&(e^{\pi i}(z_{1}-z_{2}))^{n}.
\end{eqnarray*}
This shows that indeed when 
$|q_{\tau}|<|e^{-2\pi iz_{2}}|<|e^{-2\pi iz_{1}}|<1$
and $0<|e^{2\pi i(-z_{1}+z_{2})}-1|<1$, 
(\ref{beta3.7}) is equal to the right-hand side of 
(\ref{beta3.5}) and  (\ref{beta3.5}) holds. 
Consequently, we obtain (\ref{beta3}).

Using (\ref{contrag-i}) and the $L(0)$-conjugation formula, 
we have
\begin{eqnarray}\label{beta4}
\lefteqn{E\bigg(\tr_{W^{a'_{3}}}
\Y_{a_{4}a'_{3}; k}^{a'_{3}}\bigg(\mathcal{U}(e^{-2\pi iz_{2}})
\Y_{a_{2}a'_{2}; l}^{a_{4}}\bigg(e^{\pi i L(0)}
w_{a_{2}}, e^{\pi i}(z_{1}-z_{2})\bigg)\cdot} \nn
&&\quad\quad\quad\quad\quad\quad\quad\quad\quad
\quad\quad\quad\quad\quad\quad\quad\quad\quad\cdot 
e^{\pi i L(0)}
w_{a'_{2}}, e^{-2\pi iz_{2}}\bigg)
q_{\tau}^{L(0)-\frac{c}{24}}
\bigg)\nn
&&=e^{\pi h_{a_{4}}}E\bigg(\tr_{W^{a_{3}}}
\Y_{a_{4}a_{3}; k}^{a_{3}}\bigg(\mathcal{U}(e^{2\pi iz_{2}})e^{-\pi iL(0)}\cdot \nn
&&\quad\quad\quad\quad\quad
\cdot 
\Y_{a_{2}a'_{2}; l}^{a_{4}}\bigg(e^{\pi i L(0)}
w_{a_{2}}, e^{\pi i}(z_{1}-z_{2})\bigg)e^{\pi i L(0)}
w_{a'_{2}}, e^{2\pi iz_{2}}\bigg)
q_{\tau}^{L(0)-\frac{c}{24}}
\bigg)\nn
&&=e^{\pi h_{a_{4}}}
E\bigg(\tr_{W^{a_{3}}}
\Y_{a_{4}a_{3}; k}^{a_{3}}(\mathcal{U}(e^{2\pi iz_{2}})
\Y_{a_{2}a'_{2}; l}^{a_{4}}(
w_{a_{2}}, z_{1}-z_{2})
w_{a'_{2}}, e^{2\pi iz_{2}})
q_{\tau}^{L(0)-\frac{c}{24}}
\bigg).\nn
&&
\end{eqnarray}
Combining (\ref{beta1}), (\ref{beta2}), (\ref{beta3}) and 
(\ref{beta4}), we obtain (\ref{beta0}).
\epfv

For $a_{1}, a_{2}\in \A$, we define $\alpha(\Psi_{a_{1}, a_{2}, e}^{1, 1})$ and 
$\beta(\Psi_{a_{1}, a_{2}, e}^{1, 1})$ by
\begin{eqnarray*}
(\alpha(\Psi_{a_{1}, a_{2}, e}^{1, 1}))
(w_{a}, w_{a'}; z_{1}, z_{2}; \tau)
&=&\pi(\Psi_{a_{1}, a_{2}, e}^{1, 1}
(w_{a}, w_{a'}; z_{1}, z_{2}-1; \tau)),\\
(\beta(\Psi_{a_{1}, a_{2}, e}^{1, 1}))
(w_{a}, w_{a'}; z_{1}, z_{2}; \tau)
&=&\pi(\Psi_{a_{1}, a_{2}, e}^{1, 1}
(w_{a}, w_{a'}; z_{1}, z_{2}+\tau; \tau))
\end{eqnarray*}
for $a\in \A$, $w_{a}\in W^{a}$ and $w_{a'}\in W^{a'}$. 

\begin{prop}
For $a_{1}, a_{2}\in \A$, we have 
\begin{equation}\label{alpha}
\alpha(\Psi_{a_{1}, a_{2}, e}^{1, 1})=e^{-2\pi ih_{a_{2}}}
(B^{(-1)})^{2}(\Y_{a_{1}e; 1}^{a_{1}}\otimes \Y_{a'_{2}a_{2}; 1}^{e};
\Y_{a_{1}e; 1}^{a_{1}}\otimes \Y_{a'_{2}a_{2}; 1}^{e})\Psi_{a_{1}, a_{2}, e}^{1, 1}
\end{equation}
and 
\begin{eqnarray}\label{beta}
\beta(\Psi_{a_{1}, a_{2}, e}^{1, 1})&=&e^{-2\pi i h_{a_{2}}}
\sum_{a_{3}\in \mathcal{A}}\sum_{i=1}^{N_{a_{2}a_{3}}^{a_{1}}}
\sum_{j=1}^{N_{a'_{2}a_{1}}^{a_{3}}}
F(\Y_{a_{2}e; 1}^{a_{2}}\otimes \Y_{a'_{1}a_{1}; 1}^{e};
\Y_{a'_{3}a_{1}; j}^{a_{2}} \otimes 
\Y_{a_{2}a'_{1}; i}^{a'_{3}})\cdot \nn
&&\quad\cdot 
F(\Y_{a'_{3}a_{1}; j}^{a_{2}}\otimes 
\sigma_{123}(\Y_{a_{2}a'_{1}; i}^{a'_{3}});
\Y_{ea_{2}; l}^{a_{2}}\otimes \Y_{a'_{3}a_{3}; k}^{e})
\Psi^{1, 1}_{a_{3}, a_{2}, e}.\nn
&&
\end{eqnarray}
\end{prop}
\pf
Using the definitions of $\alpha$, $\pi$
and $\Psi_{a_{1}, a_{2}, e}^{1, 1}
(w_{a_{2}}, w_{a'_{2}}; z_{1}, z_{2}; \tau)$,
(\ref{alpha0}) and (\ref{braiding}), we have
\begin{eqnarray*}
\lefteqn{(\alpha(\Psi_{a_{1}, a_{2}, e}^{1, 1}))
(w_{a_{2}}, w_{a'_{2}}; z_{1}, z_{2}; \tau)}\nn
&&=\pi(\Psi_{a_{1}, a_{2}, e}^{1, 1}
(w_{a_{2}}, w_{a'_{2}}; z_{1}, z_{2}-1; \tau))\nn
&&=\sum_{a_{3}\in \mathcal{A}}\sum_{i=1}^{N_{a_{2}a_{3}}^{a_{1}}}
\sum_{j=1}^{N_{a'_{2}a_{1}}^{a_{3}}}
\sum_{a_{4}\in \mathcal{A}}\sum_{k=1}^{N_{a_{4}a_{1}}^{a_{1}}}
\sum_{l=1}^{N_{a_{2}a'_{2}}^{a_{4}}}e^{-2\pi i(h_{a_{3}}-h_{a_{1}})}\cdot\nn
&&\quad\quad\quad \cdot F^{-1}(\Y_{ea_{1}; 1}^{a_{1}}\otimes \Y_{a_{2}a'_{2}; 1}^{e};
\Y_{a_{2}a_{3}; i}^{a_{1}}\otimes \Y_{a'_{2}a_{1}; j}^{a_{3}})
\cdot \nn
&&\quad\quad \quad  \cdot
F(\Y_{a_{2}a_{3}; i}^{a_{1}}\otimes \Y_{a'_{2}a_{1}; j}^{a_{3}};
\Y_{a_{4}a_{1}; k}^{a_{1}}\otimes \Y_{a_{2}a'_{2}; l}^{a_{4}})\cdot \nn
&&\quad\quad\quad \cdot
\pi\bigg(E\bigg(\tr_{W^{a_{1}}}
\Y_{a_{4}a_{1}; k}^{a_{1}}(\mathcal{U}(e^{2\pi iz_{2}})\cdot \nn
&&\quad\quad\quad \quad\quad\quad\quad\quad\quad\quad\cdot
\Y_{a_{2}a'_{2}; l}^{a_{4}}(w_{a_{2}},z_{1}-z_{2})w_{a'_{2}},  e^{2\pi i z_{2}})
q_{\tau}^{L(0)-\frac{c}{24}}\bigg)\bigg)\nn
&&=\sum_{a_{3}\in \mathcal{A}}\sum_{i=1}^{N_{a_{2}a_{3}}^{a_{1}}}
\sum_{j=1}^{N_{a'_{2}a_{1}}^{a_{3}}}e^{-2\pi i(h_{a_{3}}-h_{a_{1}})}
F^{-1}(\Y_{ea_{1}; 1}^{a_{1}}\otimes \Y_{a_{2}a'_{2}; 1}^{e};
\Y_{a_{2}a_{3}; i}^{a_{1}}\otimes \Y_{a'_{2}a_{1}; j}^{a_{3}})
\cdot \nn
&&\quad\quad\quad \cdot
F(\Y_{a_{2}a_{3}; i}^{a_{1}}\otimes \Y_{a'_{2}a_{1}; j}^{a_{3}};
\Y_{ea_{1}; 1}^{a_{1}}\otimes \Y_{a_{2}a'_{2}; 1}^{e})\cdot \nn
&&\quad\quad\quad \cdot
E\bigg(\tr_{W^{a_{1}}}
\Y_{ea_{1}; 1}^{a_{1}}(\mathcal{U}(e^{2\pi iz_{2}})\cdot \nn
&&\quad\quad\quad \quad\quad\quad\quad\quad\quad\quad\cdot
\Y_{a_{2}a'_{2}; 1}^{e}(w_{a_{2}},z_{1}-z_{2})w_{a'_{2}},  e^{2\pi i z_{2}})
q_{\tau}^{L(0)-\frac{c}{24}}\bigg)\nn
&&=e^{-2\pi ih_{a_{2}}}
\sum_{a_{3}\in \mathcal{A}}\sum_{i=1}^{N_{a_{2}a_{3}}^{a_{1}}}
\sum_{j=1}^{N_{a'_{2}a_{1}}^{a_{3}}}
F(\Y_{a_{1}e; 1}^{a_{1}}\otimes \Y_{a'_{2}a_{2}; 1}^{e};
\Y_{a_{3}a_{2}; i}^{a_{1}}\otimes \Y_{a_{1}a'_{2}; j}^{a_{3}})
\cdot \nn
&&\quad\quad\quad \cdot e^{-2\pi i(h_{a_{3}}-h_{a_{1}}-h_{a_{2}})}
F^{-1}(\Y_{a_{3}a_{2}; i}^{a_{1}}\otimes \Y_{a_{1}a'_{2}; j}^{a_{3}};
\Y_{a_{1}e; 1}^{a_{1}}\otimes \Y_{a'_{2}a_{2}; 1}^{e})\cdot \nn
&&\quad\quad\quad \cdot
E\bigg(\tr_{W^{a_{1}}}
\Y_{ea_{1}; 1}^{a_{1}}(\mathcal{U}(e^{2\pi iz_{2}})\cdot \nn
&&\quad\quad\quad \quad\quad\quad\quad\quad\quad\quad\cdot
\Y_{a_{2}a'_{2}; 1}^{e}(w_{a_{2}},z_{1}-z_{2})w_{a'_{2}},  e^{2\pi i z_{2}})
q_{\tau}^{L(0)-\frac{c}{24}}\bigg)\nn
&&=e^{-2\pi ih_{a_{2}}}(B^{(-1)})^{2}(\Y_{a_{1}e; 1}^{a_{1}}\otimes \Y_{a'_{2}a_{2}; 1}^{e};
\Y_{a_{1}e; 1}^{a_{1}}\otimes \Y_{a'_{2}a_{2}; 1}^{e})\cdot \nn
&&\quad\quad\quad \quad\quad\quad\quad\quad\quad\quad\quad\quad\quad\quad
\quad\quad\cdot
\Psi_{a_{1}, a_{2}, e}^{1, 1}
(w_{a_{2}}, w_{a'_{2}}; z_{1}, z_{2}; \tau),
\end{eqnarray*}
proving (\ref{alpha}).

Using the definitions of $\beta$, $\pi$ and 
$\Psi_{a_{1}, a_{2}, e}^{1, 1}
(w_{a_{2}}, w_{a'_{2}}; z_{1}, z_{2}; \tau)$, (\ref{beta0}), 
we have
\begin{eqnarray*}
\lefteqn{(\beta(\Psi_{a_{1}, a_{2}, e}^{1, 1}))
(w_{a_{2}}, w_{a'_{2}}; z_{1}, z_{2}; \tau)}\nn
&&=\pi(\Psi_{a_{1}, a_{2}, e}^{1, 1}
(w_{a_{2}}, w_{a'_{2}}; z_{1}, z_{2}+\tau; \tau))\nn
&&=\sum_{a_{3}\in \mathcal{A}}\sum_{i=1}^{N_{a_{2}a_{3}}^{a_{1}}}
\sum_{j=1}^{N_{a'_{2}a_{1}}^{a_{3}}}
\sum_{a_{4}\in \mathcal{A}}\sum_{k=1}^{N_{a_{4}a_{3}}^{a_{3}}}
\sum_{l=1}^{N_{a_{2}a'_{2}}^{a_{4}}}e^{\pi i (-2h_{a_{2}}+h_{a_{4}})}\cdot\nn
&&\quad  \cdot F^{-1}(\Y_{ea_{1}; 1}^{a_{1}}\otimes \Y_{a_{2}a'_{2}; 1}^{e};
\sigma_{23}(\Y_{a_{2}a'_{1}; i}^{a'_{3}}) \otimes 
\sigma_{13}(\Y_{a'_{3}a_{1}; j}^{a_{2}}))\cdot \nn
&&\quad \cdot 
F(\Y_{a_{2}a'_{1}; i}^{a'_{3}}\otimes \sigma_{123}(\Y_{a'_{3}a_{1}; j}^{a_{2}});
\Y_{a_{4}a'_{3}; k}^{a'_{3}}\otimes \Y_{a_{2}a'_{2}; l}^{a_{4}})\cdot \nn
&&\quad\cdot 
\pi\bigg(E\bigg(\tr_{W^{a_{3}}}
\Y_{a_{4}a_{3}; k}^{a_{3}}(\mathcal{U}(e^{2\pi iz_{2}})
\Y_{a_{2}a'_{2}; l}^{a_{4}}(
w_{a_{2}}, z_{1}-z_{2})
w_{a'_{2}}, e^{2\pi iz_{2}})
q_{\tau}^{L(0)-\frac{c}{24}}
\bigg)\bigg)\nn
&&=e^{-2\pi i h_{a_{2}}}
\sum_{a_{3}\in \mathcal{A}}\sum_{i=1}^{N_{a_{2}a_{3}}^{a_{1}}}
\sum_{j=1}^{N_{a'_{2}a_{1}}^{a_{3}}}
F^{-1}(\Y_{ea_{1}; 1}^{a_{1}}\otimes \Y_{a_{2}a'_{2}; 1}^{e};
\sigma_{23}(\Y_{a_{2}a'_{1}; i}^{a'_{3}}) \otimes 
\sigma_{13}(\Y_{a'_{3}a_{1}; j}^{a_{2}}))\cdot \nn
&&\quad\quad \cdot 
F(\Y_{a_{2}a'_{1}; i}^{a'_{3}}\otimes \sigma_{123}(\Y_{a'_{3}a_{1}; j}^{a_{2}});
\Y_{ea'_{3}; 1}^{a'_{3}}\otimes \Y_{a_{2}a'_{2}; 1}^{e})\cdot \nn
&&\quad\quad \cdot 
E\bigg(\tr_{W^{a_{3}}}
\Y_{ea_{3}; 1}^{a_{3}}(\mathcal{U}(e^{2\pi iz_{2}})
\Y_{a_{2}a'_{2}; 1}^{e}(
w_{a_{2}}, z_{1}-z_{2})
w_{a'_{2}}, e^{2\pi iz_{2}})
q_{\tau}^{L(0)-\frac{c}{24}}
\bigg)\nn
&&=e^{-2\pi i h_{a_{2}}}
\sum_{a_{3}\in \mathcal{A}}\sum_{i=1}^{N_{a_{2}a_{3}}^{a_{1}}}
\sum_{j=1}^{N_{a'_{2}a_{1}}^{a_{3}}}
F^{-1}(\Y_{ea_{1}; 1}^{a_{1}}\otimes \Y_{a_{2}a'_{2}; 1}^{e};
\sigma_{23}(\Y_{a_{2}a'_{1}; i}^{a'_{3}}) \otimes 
\sigma_{13}(\Y_{a'_{3}a_{1}; j}^{a_{2}}))\cdot \nn
&&\quad\quad \cdot 
F(\Y_{a_{2}a'_{1}; i}^{a'_{3}}\otimes \sigma_{123}(\Y_{a'_{3}a_{1}; j}^{a_{2}});
\Y_{ea'_{3}; k}^{a'_{3}}\otimes \Y_{a_{2}a'_{2}; l}^{e})\cdot \nn
&&\quad\quad \quad\quad\cdot 
\Psi^{1, 1}_{a_{3}, a_{2}, e}
(w_{a_{2}}, w_{a'_{2}}; z_{1}, z_{2}; \tau).
\end{eqnarray*}
Thus we obtain 
\begin{eqnarray}\label{formula2-4}
\lefteqn{\beta(\Psi_{a_{1}, a_{2}, e}^{1, 1})}\nn
&&=e^{-2\pi i h_{a_{2}}}
\sum_{a_{3}\in \mathcal{A}}\sum_{i=1}^{N_{a_{2}a_{3}}^{a_{1}}}
\sum_{j=1}^{N_{a'_{2}a_{1}}^{a_{3}}}
F^{-1}(\Y_{ea_{1}; 1}^{a_{1}}\otimes \Y_{a_{2}a'_{2}; 1}^{e};
\sigma_{23}(\Y_{a_{2}a'_{1}; i}^{a'_{3}}) \otimes 
\sigma_{13}(\Y_{a'_{3}a_{1}; j}^{a_{2}}))\cdot \nn
&&\quad\quad \quad\quad\quad\quad\cdot 
F(\Y_{a_{2}a'_{1}; i}^{a'_{3}}\otimes \sigma_{123}(\Y_{a'_{3}a_{1}; j}^{a_{2}});
\Y_{ea'_{3}; k}^{a'_{3}}\otimes \Y_{a_{2}a'_{2}; l}^{e})
\Psi^{1, 1}_{a_{3}, a_{2}, e}.
\end{eqnarray}

Now using (\ref{fusing-inv}), (\ref{fusing}) and the 
relations $\sigma_{12}\sigma_{23}=\sigma_{123}$,
$\sigma_{12}\sigma_{13}=\sigma_{132}$ and 
$\sigma_{123}\sigma_{132}=\sigma_{132}\sigma_{123}=1$, we have 
\begin{eqnarray}\label{formula2-5}
\lefteqn{F^{-1}(\Y_{ea_{1}; 1}^{a_{1}}\otimes \Y_{a_{2}a'_{2}; 1}^{e};
\sigma_{23}(\Y_{a_{2}a'_{1}; i}^{a'_{3}}) \otimes 
\sigma_{13}(\Y_{a'_{3}a_{1}; j}^{a_{2}}))}\nn
&&=F(\sigma_{12}(\Y_{ea_{1}; 1}^{a_{1}})\otimes 
\sigma_{12}(\Y_{a_{2}a'_{2}; 1}^{e});
\sigma_{12}(\sigma_{23}(\Y_{a_{2}a'_{1}; i}^{a'_{3}})) \otimes 
\sigma_{12}(\sigma_{13}(\Y_{a'_{3}a_{1}; j}^{a_{2}})))\nn
&&=F(\Y_{a_{1}e; 1}^{a_{1}}\otimes \Y_{a'_{2}a_{2}; 1}^{e};
\sigma_{12}(\sigma_{23}(\Y_{a_{2}a'_{1}; i}^{a'_{3}})) \otimes 
\sigma_{12}(\sigma_{13}(\Y_{a'_{3}a_{1}; j}^{a_{2}})))\nn
&&=F(\Y_{a_{2}e; 1}^{a_{2}}\otimes \Y_{a'_{1}a_{1}; 1}^{e};
\sigma_{123}(\sigma_{132}(\Y_{a'_{3}a_{1}; j}^{a_{2}})) \otimes 
\sigma_{132}(\sigma_{123}(\Y_{a_{2}a'_{1}; i}^{a'_{3}})))\nn
&&=F(\Y_{a_{2}e; 1}^{a_{2}}\otimes \Y_{a'_{1}a_{1}; 1}^{e};
\Y_{a'_{3}a_{1}; j}^{a_{2}} \otimes 
\Y_{a_{2}a'_{1}; i}^{a'_{3}})
\end{eqnarray}
and 
\begin{eqnarray}\label{formula2-6}
\lefteqn{F(\Y_{a_{2}a'_{1}; i}^{a'_{3}}\otimes \sigma_{123}(\Y_{a'_{3}a_{1}; j}^{a_{2}});
\Y_{ea'_{3}; k}^{a'_{3}}\otimes \Y_{a_{2}a'_{2}; l}^{e})}\nn
&&=F(\sigma_{132}(\sigma_{123}(\Y_{a'_{3}a_{1}; j}^{a_{2}}))\otimes 
\sigma_{123}(\Y_{a_{2}a'_{1}; i}^{a'_{3}});
\sigma_{123}(\Y_{a_{2}a'_{2}; l}^{e}) \otimes 
\sigma_{132}(\Y_{ea'_{3}; k}^{a'_{3}}))\nn
&&=F(\Y_{a'_{3}a_{1}; j}^{a_{2}}\otimes 
\sigma_{123}(\Y_{a_{2}a'_{1}; i}^{a'_{3}});
\Y_{ea_{2}; l}^{a_{2}}\otimes \Y_{a'_{3}a_{3}; k}^{e}).
\end{eqnarray}
From (\ref{formula2-4})--(\ref{formula2-6}), we obtain (\ref{beta}). 
\epfv

By Proposition \ref{1-independence}, $\Psi_{a_{1}, a_{2}, e}^{1, 1}$,
$a_{1}\in \mathcal{A}$ form a basis of $\mathcal{F}_{1;2}^{e}$. For
fixed $a_{2}\in \mathcal{A}$, we use $\alpha_{a_{1}}^{a_{3}}(a_{2})$
and $\beta_{a_{1}}^{a_{3}}(a_{2})$, $a_{1}, a_{3}\in \mathcal{A}$, to
denote the matrix elements of $\alpha$ and $\beta$, respectively,
under the basis $\Psi_{a_{1}, a_{2}, e}^{1, 1}$, $a_{1}\in
\mathcal{A}$.

\begin{cor}
The matrix elements $\alpha_{a_{1}}^{a_{3}}(a_{2})$ and 
$\beta_{a_{1}}^{a_{3}}(a_{2})$, $a_{1}, a_{3}\in \mathcal{A}$,
are given by
\begin{equation}\label{alpha-1}
\alpha_{a_{1}}^{a_{3}}(a_{2})=\delta_{a_{1}a_{3}}
(B^{(-1)})^{2}(\Y_{a_{1}e; 1}^{a_{1}}\otimes \Y_{a'_{2}a_{2}; 1}^{e};
\Y_{a_{1}e; 1}^{a_{1}}\otimes \Y_{a'_{2}a_{2}; 1}^{e})
\end{equation}
and
\begin{eqnarray}\label{beta-1}
\beta_{a_{1}}^{a_{3}}(a_{2})&=&e^{-2\pi i h_{a_{2}}}
\sum_{i=1}^{N_{a_{2}a_{3}}^{a_{1}}}\sum_{j=1}^{N_{a'_{2}a_{1}}^{a_{3}}}
F(\Y_{a_{2}e; 1}^{a_{2}}\otimes \Y_{a'_{1}a_{1}; 1}^{e};
\Y_{a'_{3}a_{1}; j}^{a_{2}} \otimes 
\Y_{a_{2}a'_{1}; i}^{a'_{3}})\cdot \nn
&&\quad\quad \quad\quad\cdot 
F(\Y_{a'_{3}a_{1}; j}^{a_{2}}\otimes 
\sigma_{123}(\Y_{a_{2}a'_{1}; i}^{a'_{3}});
\Y_{ea_{2}; l}^{a_{2}}\otimes \Y_{a'_{3}a_{3}; k}^{e}).\nn
&&
\end{eqnarray}
\end{cor}
\pf
This corollary follows directly from the definition of 
$\alpha_{a_{1}}^{a_{3}}(a_{2})$ and 
$\beta_{a_{1}}^{a_{3}}(a_{2})$, $a_{1}, a_{3}\in \mathcal{A}$,
(\ref{alpha}) and (\ref{beta}). 
\epfv

It is also easy to establish the relationship between 
$\alpha$ and $\beta$:

\begin{prop}
We have the following formula:
\begin{equation}\label{formula2-3}
S\alpha S^{-1}=\beta.
\end{equation}
\end{prop}
\pf
We have 
\begin{eqnarray*}
\lefteqn{(\beta(S(\Psi_{a_{1}, a_{2}, e}^{1, 1})))
(w_{a_{2}}, w_{a_{2}'}; z_{1}, z_{2}; \tau)}\nn
&&=\pi(S(\Psi_{a_{1}, a_{2}, e}^{1, 1})
(w_{a_{2}}, w_{a_{2}'}; z_{1}, z_{2}+\tau; \tau))\nn
&&=\pi\bigg(E\Bigg(\tr_{W^{a_{1}}}
\Y_{ea_{1}; 1}^{a_{1}}(\mathcal{U}(e^{-2\pi i\frac{z_{2}+\tau}{\tau}})
\left(-\frac{1}{\tau}\right)^{L(0)}\cdot\nn
&&\quad\quad\quad\quad\quad\cdot \Y_{a_{2}a'_{2}; 1}^{e}
(w_{a_{2}}, z_{1}-z_{2}+\tau)w_{a'_{2}}, 
e^{-2\pi i\frac{z_{2}+\tau}{\tau}})
q_{-\frac{1}{\tau}}^{L(0)-\frac{c}{24}}\Bigg)\bigg)\nn
&&=\pi\bigg(E\Bigg(\tr_{W^{a_{1}}}
\Y_{ea_{1}; 1}^{a_{1}}\Bigg(\mathcal{U}(e^{2\pi i (-\frac{z_{2}}{\tau}-1)})\cdot\nn
&&\quad\quad\quad \quad\cdot
\Y_{a_{2}a'_{2}; 1}^{e}
\Bigg(\left(-\frac{1}{\tau}\right)^{L(0)}w_{a_{2}}, 
-\frac{1}{\tau}z_{1}-(-\frac{1}{\tau}z_{2}-1)\Bigg)\cdot\nn
&&\quad\quad\quad \quad\cdot \left(-\frac{1}{\tau}\right)^{L(0)}
w_{a'_{2}}, e^{2\pi i (-\frac{z_{2}}{\tau}-1)}\Bigg)
q_{-\frac{1}{\tau}}^{L(0)-\frac{c}{24}}\Bigg)\bigg)\nn
&&=\pi\left(\Psi_{a_{1}, a_{2}, e}^{1, 1}
\left(\left(-\frac{1}{\tau}\right)^{L(0)}w_{a_{2}}, 
\left(-\frac{1}{\tau}\right)^{L(0)}w_{a_{2}'}; -\frac{1}{\tau}z_{1}, -\frac{1}{\tau}z_{2}-1; 
-\frac{1}{\tau}\right)\right)\nn
&&=(\alpha(\Psi_{a_{1}, a_{2}, e}^{1, 1}))
\left(\left(-\frac{1}{\tau}\right)^{L(0)}w_{a_{2}}, 
\left(-\frac{1}{\tau}\right)^{L(0)}w_{a_{2}'}; -\frac{1}{\tau}z_{1}, -\frac{1}{\tau}z_{2}; 
-\frac{1}{\tau}\right)\nn
&&=(S(\alpha(\Psi_{a_{1}, a_{2}, e}^{1, 1})))
(w_{a_{2}}, 
w_{a_{2}'}; z_{1}, z_{2}; \tau).
\end{eqnarray*}
Thus we obtain
$$\beta S=S\alpha$$
or equivalently
(\ref{formula2-3}).
\epfv

With the bases of the spaces of intertwining operators 
we choose in the beginning of this section, 
we have
\begin{eqnarray*}
\Psi_{a_{1}}(u; \tau)&=&\tr_{W^{a_{1}}}
\Y_{ea_{1}; 1}^{a_{1}}(\mathcal{U}(e^{2\pi iz})u, e^{2\pi iz})
q_{\tau}^{L(0)-\frac{c}{24}}\nn
&=&\tr_{W^{a_{1}}}
Y_{W^{a_{1}}}(\mathcal{U}(e^{2\pi iz})u, e^{2\pi iz})
q_{\tau}^{L(0)-\frac{c}{24}}
\end{eqnarray*}
for $a_{1}\in \A$, $u\in V$ and 
\begin{eqnarray*}
\lefteqn{\Psi^{1, 1}_{a_{1},  a_{2}, e}(w_{a_{2}}, w_{a_{2}'}; 
z_{1}, z_{2}; \tau)}\nn
&&=
E(\tr_{W^{a_{1}}}
\Y_{ea_{1}; 1}^{a_{1}}(\mathcal{U}(e^{2\pi iz_{2}})\Y_{a_{2}a'_{2}; 1}^{e}
(w_{a_{2}}, z_{1}-z_{2})w_{a'_{2}}, e^{2\pi i z_{2}})
q_{\tau}^{L(0)-\frac{c}{24}})\nn
&&=
E(\tr_{W^{a_{1}}}
Y_{W^{a_{1}}}(\mathcal{U}(e^{2\pi iz_{2}})\Y_{a_{2}a'_{2}; 1}^{e}
(w_{a_{2}}, z_{1}-z_{2})w_{a'_{2}}, e^{2\pi i z_{2}})
q_{\tau}^{L(0)-\frac{c}{24}})
\end{eqnarray*}
for $a_{1}, a_{2}\in \A$, $w_{a_{2}}\in W^{a_{2}}$ and 
$w_{a'_{2}}\in W^{a'_{2}}$. 

Since $\Psi_{a_{1}}$ for $a_{1}\in \A$
are linear independent, they form a basis of $\F^{e}_{1;1}$. Thus
we know that there exist unique 
$S^{a_{3}}_{a_{1}}\in \C$ for $a_{1}, a_{3}\in \A$
such that 
\begin{eqnarray}
(S(\Psi_{a_{1}}))
(u; \tau)
&=&\sum_{a_{3}\in \A}S^{a_{1}}_{a_{3}}
\Psi_{a_{3}}
(u; \tau),\label{formula2-1}\\
S(\Psi_{a_{1}})
&=&\sum_{a_{3}\in \A}S^{a_{1}}_{a_{3}}
\Psi_{a_{3}}.\label{formula2-2}
\end{eqnarray}
Clearly we also have 
$$
S(\Psi_{a_{1}, a_{2}, e}^{1, 1})
=\sum_{a_{3}\in \A}S^{a_{1}}_{a_{3}}
\Psi_{a_{3}, a_{2}, e}^{1,1}
$$
for $a_{1}, a_{2}\in \A$.

We can identify the space spanned by $\Psi_{a}$ for $a\in \A$
with the vector space $\coprod_{a\in \A}\C[W^{a}]$
where $[W^{a}]$ is the equivalence class of $W^{a}$ and is 
actually equal to $a$. Then we can view $S$ as a 
linear operator on this vector space spanned by $\A$. In terms of the 
basis $[W^{a}]$, $a\in \A$, we have 
$$S([W^{a_{1}}])=\sum_{a_{2}\in \A}S^{a_{1}}_{a_{2}}[W^{a_{2}}].$$

Let $u=\mathbf{1}$ in (\ref{formula2-1}). Then we see that $S_{a_{1}}^{a_{2}}$, 
$a_{1}, a_{2}\in \A$, gives an
action of the modular transformation $\tau\mapsto -1/\tau$ on the space 
spanned by shifted graded dimensions (vacuum characters) of irreducible 
$V$-modules. 

The following theorem gives the second Moore-Seiberg formula 
in \cite{MS1}:

\begin{thm}
For $a_{1}, a_{2}, a_{3}\in \A$, we have
\begin{eqnarray}\label{formula2}
\lefteqn{\sum_{a_{4}\in \mathcal{A}}S_{a_{1}}^{a_{4}}
(B^{(-1)})^{2}(\Y_{a_{4}e; 1}^{a_{4}}\otimes \Y_{a'_{2}a_{2}; 1}^{e};
\Y_{a_{4}e; 1}^{a_{4}}\otimes \Y_{a'_{2}a_{2}; 1}^{e})
(S^{-1})_{a_{4}}^{a_{3}}}\nn
&&=\sum_{i=1}^{N_{a_{1}a_{2}}^{a_{3}}}\sum_{k=1}^{N_{a'_{1}a_{3}}^{a_{2}}}
F(\Y_{a_{2}e; 1}^{a_{2}}\otimes \Y_{a'_{3}a_{3}; 1}^{e}; 
\Y_{a'_{1}a_{3}; k}^{a_{2}}\otimes \Y_{a_{2}a'_{3}; i}^{a'_{1}})\cdot\nn
&&\quad\quad\quad\quad\quad\cdot 
F(\Y_{a'_{1}a_{3}; k}^{a_{2}}\otimes \sigma_{123}(\Y_{a_{2}a'_{3}; i}^{a'_{1}});
\Y_{ea_{2}; 1}^{a_{2}}\otimes \Y_{a'_{1}a_{1}; 1}^{e})
\end{eqnarray}
\end{thm}
\pf
This follows from (\ref{formula2-3}), (\ref{alpha-1}) and (\ref{beta-1})
immediately.
\epfv

\begin{cor}
For $a_{1}, a_{2}, a_{3}\in \A$, we have
\begin{eqnarray}\label{diag1}
&{\displaystyle \sum_{a_{4}\in \mathcal{A}}S_{a_{1}}^{a_{4}}
(B^{(-1)})^{2}(\Y_{a_{4}e; 1}^{a_{4}}\otimes \Y_{a'_{2}a_{2}; 1}^{e};
\Y_{a_{4}e; 1}^{a_{4}}\otimes \Y_{a'_{2}a_{2}; 1}^{e})
(S^{-1})_{a_{4}}^{a_{3}}}&\nn
&=N_{a_{1}a_{2}}^{a_{3}}
F(\Y_{a_{2}e; 1}^{a_{2}} \otimes \Y_{a'_{2}a_{2}; 1}^{e};
\Y_{ea_{2}; 1}^{a_{2}}\otimes \Y_{a_{2}a'_{2}; 1}^{e}).&
\end{eqnarray}
\end{cor}
\pf 
This follows immediately from  (\ref{formula1}) and (\ref{formula2}).
\epfv

\renewcommand{\theequation}{\thesection.\arabic{equation}}
\renewcommand{\thethm}{\thesection.\arabic{thm}}
\setcounter{equation}{0}
\setcounter{thm}{0}

\section{The main theorem and the Verlinde formula}

In this section, using the results obtained in the preceding section,
we prove the main theorem, the Verlinde conjecture, of the present 
paper, derive the Verlinde formula for fusion rules and prove that 
$(S_{a_{1}}^{a_{2}})$ is symmetric. 

First, we have the following:

\begin{prop}\label{nonzero}
For $a_{2}\in \A$, 
$$F(\Y_{a_{2}e; 1}^{a_{2}} \otimes \Y_{a'_{2}a_{2}; 1}^{e};
\Y_{ea_{2}; 1}^{a_{2}}\otimes \Y_{a_{2}a'_{2}; 1}^{e})\ne 0.$$
\end{prop}
\pf
If 
$$F(\Y_{a_{2}e; 1}^{a_{2}} \otimes \Y_{a'_{2}a_{2}; 1}^{e};
\Y_{ea_{2}; 1}^{a_{2}}\otimes \Y_{a_{2}a'_{2}; 1}^{e})= 0,$$
then by (\ref{diag1}), 
$$(B^{(-1)})^{2}(\Y_{a_{4}e; 1}^{a_{4}}\otimes \Y_{a'_{2}a_{2}; 1}^{e};
\Y_{a_{4}e; 1}^{a_{4}}\otimes \Y_{a'_{2}a_{2}; 1}^{e})=0$$
for $a_{4}\in \mathcal{A}$. But we know that 
$$(B^{(-1)})^{2}(\Y_{ee; 1}^{e}\otimes \Y_{a'_{2}a_{2}; 1}^{e};
\Y_{ee; 1}^{e}\otimes \Y_{a'_{2}a_{2}; 1}^{e})=1.$$
Contradiction. 
\epfv

For $a_{2}\in \A$, let $\mathbb{N}(a_{2})$
be the matrices whose entries are $N_{a_{1}a_{2}}^{a_{3}}=N_{a_{2}a_{1}}^{a_{3}}$
for $a_{1}, a_{3}\in \A$, that is, 
$$\mathcal{N}(a_{2})=(N_{a_{1}a_{2}}^{a_{3}})=(N_{a_{2}a_{1}}^{a_{3}}).$$
Then we have the  main result of the present
paper:

\begin{thm}\label{main}
Let $V$ be a simple vertex operator algebra satisfying the conditions 
in Section 1. Then we have
\begin{equation}\label{diag2}
\sum_{a_{1}, a_{3}\in \mathcal{A}}(S^{-1})_{a_{4}}^{a_{1}}
N_{a_{1}a_{2}}^{a_{3}}
S_{a_{3}}^{a_{5}}=\delta_{a_{4}}^{a_{5}}
\frac{(B^{(-1)})^{2}(\Y_{a_{4}e; 1}^{a_{4}}\otimes \Y_{a'_{2}a_{2}; 1}^{e};
\Y_{a_{4}e; 1}^{a_{4}}\otimes \Y_{a'_{2}a_{2}; 1}^{e})}
{F(\Y_{a_{2}e; 1}^{a_{2}} \otimes \Y_{a'_{2}a_{2}; 1}^{e};
\Y_{ea_{2}; 1}^{a_{2}}\otimes \Y_{a_{2}a'_{2}; 1}^{e})}.
\end{equation}
In particular, the matrix $S$ diagonalizes the
matrices $\mathcal{N}(a_{2})$ for all  $a_{2}\in \A$.
\end{thm}
\pf
By Proposition \ref{nonzero}, we can rewrite (\ref{diag1}) as 
(\ref{diag2}). 
Since the right-hand side of (\ref{diag2}) are entries
of diagonal matrices, $S$ diagonalizes the
$\mathcal{N}(a_{2})$ for all  $a_{2}\in \A$.
\epfv

We now prove the 
Verlinde formula for fusion rules. We first need the following:

\begin{prop}
The square $S^{2}$ viewed as a linear operator on the vector 
space spanned by $\A$  is equal to the linear operator 
obtained from the map $': \A\to \A$.
\end{prop}
\pf
By definition, we have 
\begin{eqnarray}\label{s-2}
(S^{2}(\Psi_{a}))
(u; \tau)&=&
(S(\Psi_{a})\left(\left(-\frac{1}{\tau}\right)^{L(0)}u; -\frac{1}{\tau}\right)\nn
&=&\Psi_{a}\left(\left(-\frac{1}{-\frac{1}{\tau}}\right)^{L(0)}
\left(-\frac{1}{\tau}\right)^{L(0)}u; 
-\frac{1}{-\frac{1}{\tau}}\right)\nn
&=&\Psi_{a}\left(\tau^{L(0)}
\left(-\frac{1}{\tau}\right)^{L(0)}u; 
\tau\right)\nn
&=&\Psi_{a}\left(e^{(\log \tau+\log (-\frac{1}{\tau}))L(0)}u; 
\tau\right).
\end{eqnarray}
Note that both $\tau$ and $-\frac{1}{\tau}$ are in the upper
half plane. So $0<\arg \tau, \arg (-\frac{1}{\tau})<\pi$. 
Thus by our convention, 
$$\arg \tau+ \arg (-\frac{1}{\tau})=\arg (-1)=\pi.$$
So we have
\begin{eqnarray}\label{sum-log}
\log \tau+\log \left(-\frac{1}{\tau}\right)&=&\log |\tau|+i \arg \tau
\log \left|-\frac{1}{\tau}\right|+ i\arg \left(-\frac{1}{\tau}\right)\nn
&=&\pi i.
\end{eqnarray}
Using (\ref{sum-log}) and (\ref{contrag-1}), we 
see that the right-hand side of (\ref{s-2})
is equal to 
\begin{eqnarray}\label{s-2-1}
\Psi_{a}(e^{\pi i L(0)}u; 
\tau)&=&\tr_{W^{a}}
Y_{W^{a}}(\mathcal{U}(e^{-2\pi iz})e^{\pi i L(0)}u, 
e^{-2\pi i z})
q_{\tau}^{L(0)-\frac{c}{24}}\nn
&=&\tr_{W^{a'}}
Y_{W^{a'}}(\mathcal{U}(e^{2\pi iz_{2}})u, 
e^{2\pi i z})
q_{\tau}^{L(0)-\frac{c}{24}})\nn
&=&\Psi_{a'}\left(u; 
\tau\right)
\end{eqnarray} 
Combining (\ref{s-2}) and (\ref{s-2-1}), 
we obtain 
$$S^{2}(\Psi_{a})=\Psi_{a'},$$
proving the conclusion. 
\epfv 

An immediate consequence of the proposition above is the 
following:

\begin{cor}
The inverse $S^{-1}$ of $S$ is equal to $S\circ \; '='\circ S$.
In particular, we have 
\begin{eqnarray}\label{s-inv}
(S^{-1})_{a_{1}}^{a_{2}}&=&S_{a_{1}}^{a'_{2}}\nn
&=&S_{a'_{1}}^{a_{2}}
\end{eqnarray}
for $a_{1}, a_{2}\in \A$.
\end{cor}

Now we have the following Verlinde formula for fusion rules:

\begin{thm}
Let $V$ be a vertex operator algebra satisfying the conditions
in Section 1. Then we have $S_{e}^{a}\ne 0$ for $a\in \mathcal{A}$ and
\begin{equation}\label{v-form}
N_{a_{1}a_{2}}^{a_{3}}=
\sum_{a_{4}\in \A}\frac{S_{a_{1}}^{a_{4}}S_{a_{2}}^{a_{4}}S_{a_{4}}^{a'_{3}}}
{S_{e}^{a_{4}}}.
\end{equation}
\end{thm}
\pf
Let 
\begin{equation}\label{lambda}
\lambda_{a_{2}}^{a_{4}}=\frac{(B^{(-1)})^{2}(\Y_{a_{4}e; 1}^{a_{4}}\otimes \Y_{a'_{2}a_{2}; 1}^{e};
\Y_{a_{4}e; 1}^{a_{4}}\otimes \Y_{a'_{2}a_{2}; 1}^{e})}
{F(\Y_{a_{2}e; 1}^{a_{2}} \otimes \Y_{a'_{2}a_{2}; 1}^{e};
\Y_{ea_{2}; 1}^{a_{2}}\otimes \Y_{a_{2}a'_{2}; 1}^{e})}
\end{equation}
for $a_{2}, a_{4}\in \A$. Then by (\ref{diag2}), we have 
$$\sum_{a_{1}, a_{3}\in \mathcal{A}}(S^{-1})_{a_{4}}^{a_{1}}
N_{a_{1}a_{2}}^{a_{3}}
S_{a_{3}}^{a_{5}}=\delta_{a_{4}}^{a_{5}}\lambda_{a_{2}}^{a_{4}},$$
or equivalently
\begin{equation}\label{v-form-1}
N_{a_{1}a_{2}}^{a_{3}}
=\sum_{a_{4}\in\A}
S_{a_{1}}^{a_{4}}\lambda_{a_{2}}^{a_{4}}(S^{-1})_{a_{4}}^{a_{3}}.
\end{equation}
Using (\ref{s-inv}), we see that (\ref{v-form-1}) becomes 
\begin{equation}\label{v-form-2}
N_{a_{1}a_{2}}^{a_{3}}
=\sum_{a_{4}\in\A}
S_{a_{1}}^{a_{4}}\lambda_{a_{2}}^{a_{4}}S_{a_{4}}^{a'_{3}}.
\end{equation}

We know that $N_{ea_{2}}^{a_{3}}=\delta_{a_{2}}^{a_{3}}$.
Combining this fact with (\ref{v-form-1}), we obtain
$$\delta_{a_{2}}^{a_{3}}
=\sum_{a_{4}\in\A}
S_{e}^{a_{4}}\lambda_{a_{2}}^{a_{4}}(S^{-1})_{a_{4}}^{a_{3}}.$$
Thus we have
\begin{equation}\label{v-form-2.5}
S_{e}^{a_{4}}\lambda_{a_{2}}^{a_{4}}=S_{a_{2}}^{a_{4}}.
\end{equation}
From (\ref{v-form-2.5}) we see that if $S_{e}^{a_{4}}=0$ 
for some $a_{4}\in \mathcal{A}$, then there is one column 
of the matrix $S$ is $0$. Contradictary to the fact that 
$S$ is invertible. So $S_{e}^{a_{4}}\ne 0$ 
for $a_{4}\in \mathcal{A}$. Rewrite (\ref{v-form-2.5}) as
\begin{equation}\label{v-form-3}
\lambda_{a_{2}}^{a_{4}}=\frac{S_{a_{2}}^{a_{4}}}{S_{e}^{a_{4}}}.
\end{equation}
Substituting (\ref{v-form-3}) into (\ref{v-form-2}), we obtain
(\ref{v-form}). 
\epfv

Finally we have:

\begin{thm}
The matrix $(S_{a_{1}}^{a_{2}})$ is symmetric. 
\end{thm}
\pf 
Rewriting (\ref{v-form-1}) as 
\begin{equation}\label{v-form-4}
\sum_{a_{1}}(S^{-1})_{a_{4}}^{a_{1}}N_{a_{1}a_{2}}^{a_{3}}
=\lambda_{a_{2}}^{a_{4}}(S^{-1})_{a_{4}}^{a_{3}},
\end{equation}
and then letting $a_{3}=a_{4}=e$ in (\ref{v-form-4}) and using
$N_{a_{1}a_{2}}^{e}=\delta_{a_{1}}^{a'_{2}}$,
we obtain 
\begin{equation}\label{s-form-1}
(S^{-1})_{e}^{a'_{2}}=\lambda_{a_{2}}^{e}(S^{-1})_{e}^{e}.
\end{equation}
Using  (\ref{s-inv}), (\ref{s-form-1}) can be written as 
\begin{equation}\label{s-form-2}
S_{e}^{a_{2}}=\lambda_{a_{2}}^{e}S_{e}^{e}.
\end{equation}
By (\ref{v-form-3}), (\ref{s-form-2})  and 
(\ref{lambda}), 
\begin{eqnarray}\label{s-form-3}
\lefteqn{S_{a_{2}}^{a_{4}}
=\lambda_{a'_{4}}^{e}S_{e}^{e}\lambda_{a_{2}}^{a_{4}}}\nn
&&=\frac{S_{e}^{e}((B^{(-1)})^{2}(\Y_{a_{4}e; 1}^{a_{4}}\otimes \Y_{a'_{2}a_{2}; 1}^{e};
\Y_{a_{4}e; 1}^{a_{4}}\otimes \Y_{a'_{2}a_{2}; 1}^{e}))}
{F(\Y_{a_{2}e; 1}^{a_{2}} \otimes \Y_{a'_{2}a_{2}; 1}^{e};
\Y_{ea_{2}; 1}^{a_{2}}\otimes \Y_{a_{2}a'_{2}; 1}^{e})
F(\Y_{a_{4}e; 1}^{a_{4}} \otimes \Y_{a'_{4}a_{4}; 1}^{e};
\Y_{ea_{4}; 1}^{a_{4}}\otimes \Y_{a_{4}a'_{4}; 1}^{e})}.\nn
&&
\end{eqnarray}

By (\ref{s-form-3}), we obtain
\begin{equation}\label{s-form-4}
S_{a_{2}}^{a'_{4}}
=\frac{S_{e}^{e}((B^{(-1)})^{2}(\Y_{a'_{4}e; 1}^{a'_{4}}\otimes \Y_{a'_{2}a_{2}; 1}^{e};
\Y_{a'_{4}e; 1}^{a'_{4}}\otimes \Y_{a'_{2}a_{2}; 1}^{e}))}
{F(\Y_{a_{2}e; 1}^{a_{2}} \otimes \Y_{a'_{2}a_{2}; 1}^{e};
\Y_{ea_{2}; 1}^{a_{2}}\otimes \Y_{a_{2}a'_{2}; 1}^{e})
F(\Y_{a'_{4}e; 1}^{a'_{4}} \otimes \Y_{a_{4}a'_{4}; 1}^{e};
\Y_{ea'_{4}; 1}^{a'_{4}}\otimes \Y_{a'_{4}a_{4}; 1}^{e})}.
\end{equation}
From (\ref{braiding}), (\ref{fusing}), (\ref{fusing-inv}) and 
the choice of the bases of the spaces of intertwining operators
when some of the modules involved are $V$,
we have
\begin{eqnarray}\label{b2-sym}
\lefteqn{(B^{(-1)})^{2}(\Y_{a'_{4}e; 1}^{a'_{4}}\otimes \Y_{a'_{2}a_{2}; 1}^{e};
\Y_{a'_{4}e; 1}^{a'_{4}}\otimes \Y_{a'_{2}a_{2}; 1}^{e})}\nn
&&=(B^{(-1)})^{2}(\Y_{a_{2}e; 1}^{a_{2}}\otimes \Y_{a_{4}a'_{4}; 1}^{e};
\Y_{a_{2}e; 1}^{a_{2}}\otimes \Y_{a_{4}a'_{4}; 1}^{e}).
\end{eqnarray}
Using (\ref{b2-sym}) and (\ref{s-form-3}), we see that the right-hand 
side of (\ref{s-form-4}) is equal to 
\begin{equation}\label{s-form-5}
\displaystyle 
\frac{S_{e}^{e}((B^{(-1)})^{2}(\Y_{a_{2}e; 1}^{a_{2}}\otimes \Y_{a_{4}a'_{4}; 1}^{e};
\Y_{a_{2}e; 1}^{a_{2}}\otimes \Y_{a_{4}a'_{4}; 1}^{e}))}
{F(\Y_{a'_{2}e; 1}^{a'_{2}} \otimes \Y_{a_{2}a'_{2}; 1}^{e};
\Y_{ea'_{2}; 1}^{a'_{2}}\otimes \Y_{a'_{2}a_{2}; 1}^{e})
F(\Y_{a_{4}e; 1}^{a_{4}} \otimes \Y_{a'_{4}a_{4}; 1}^{e};
\Y_{ea'_{4}; 1}^{a'_{4}}\otimes \Y_{a_{4}a'_{4}; 1}^{e})}
=S_{a'_{4}}^{a_{2}}.
\end{equation}
The formulas (\ref{s-form-4}) and (\ref{s-form-5}) 
gives
$$S_{a_{2}}^{a'_{4}}=S_{a'_{4}}^{a_{2}},$$
proving that $(S_{a_{1}}^{a_{2}})$ is symmetric.
\epfv

\noindent {\small \sc Max Planck Institute of Mathematics,
Vivatsgasse 7, D-53111 Bonn, Germany}
\vspace{1em}

\noindent {\small \sc Institute of Mathematics, Fudan University,
Shanghai, China}
\vspace{1em}

\noindent {\small \sc Department of Mathematics, Rutgers University,
110 Frelinghuysen Rd., Piscataway, NJ 08854-8019 (permanent address)}

\noindent {\em E-mail address}: yzhuang@math.rutgers.edu

\end{document}